\documentclass[leqno,a4paper,12pt]{article}
\usepackage{amsmath,amsthm,bbm}
\usepackage[sans]{dsfont}

\newtheorem{Thm}{\indent Theorem}[section]
\newtheorem{MainThm}[Thm]{\indent Main Theorem}
\newtheorem{Prop}[Thm]{\indent Proposition}
\newtheorem{Lem}[Thm]{\indent Lemma}
\newtheorem{Cor}[Thm]{\indent Corollary}
\newtheorem{Conj}[Thm]{\indent Conjecture}
\newtheorem{Sum}[Thm]{\indent Summary}

\theoremstyle{definition}
\newtheorem{Def}[Thm]{\indent Definition}
\newtheorem{Rem}[Thm]{\indent Remark}
\newtheorem{Ex}[Thm]{\indent Example}
\newtheorem{Exs}[Thm]{\indent Examples}
\newtheorem{Ass}[Thm]{\indent Assumption}

\def\qed{{\hskip0pt\unskip\unskip\nobreak\hfil\penalty50
          \hskip1em\hbox{}\nobreak\hfil
          {\bf q.e.d.}%
          \parfillskip=0pt\finalhyphendemerits=0
          \par}\medskip}

\newenvironment{Proof}
               {{\it Proof.}\quad}
               {\qed}

\newenvironment{Proofof}[1]
               {{\it Proof of #1.}\quad}
               {\qed}

\newcommand{\Prime}{\kern3\fontdimen1\font$'$\kern-7\fontdimen1\font}

\long\def\forget#1{}

\long\def\beginSIDEREMARK#1\endSIDEREMARK
    {{\par\bigskip\advance\leftskip by 2cm
                  \advance\rightskip by -2cm\noindent
      {\bf Our own side remark:} #1
      \par\bigskip\noindent}}

\long\def\beginFORGET#1\endFORGET{#1}
\long\def\beginFORGET#1\endFORGET{}

\def\?{\ ???\ \immediate\write16{}%
\immediate\write16{Warning: There was still a question mark . . . }%
\immediate\write16{}}

\usepackage{amsmath}

\usepackage{exscale}
\usepackage{amssymb}

\input cyracc.def
\font\tencyr=wncyr6
\def\cyr{\tencyr\cyracc}
\newcommand{\cyrb}{{\cyr B}}

\newcommand{\A}{{\rm{A}}}
\newcommand{\B}{{\rm{B}}}

\newcommand{\BA}{{\mathbb{A}}}
\newcommand{\BB}{{\mathbb{B}}}
\newcommand{\BC}{{\mathbb{C}}}
\newcommand{\BD}{{\mathbb{D}}}

\newcommand{\BF}{{\mathbb{F}}}

\newcommand{\BQ}{{\mathbb{Q}}}
\newcommand{\BR}{{\mathbb{R}}}

\newcommand{\BZ}{{\mathbb{Z}}}

\newcommand{\Fg}{{\mathfrak{g}}}

\newcommand{\Fn}{{\mathfrak{n}}}
\newcommand{\Fo}{{\mathfrak{o}}}

\newcommand{\FA}{{\mathfrak{A}}}
\newcommand{\FB}{{\mathfrak{B}}}

\newcommand{\FH}{{\mathfrak{H}}}
\newcommand{\FI}{{\mathfrak{I}}}

\newcommand{\FS}{{\mathfrak{S}}}

\newcommand{\FX}{{\mathfrak{X}}}

\newcommand{\CC}{{\cal C}}
\newcommand{\CD}{{\cal D}}

\newcommand{\CF}{{\cal F}}

\newcommand{\CH}{{\cal H}}

\newcommand{\CK}{{\cal K}}

\newcommand{\CN}{{\cal N}}

\newcommand{\ch}{\mathop{\rm CH}\nolimits}
\newcommand{\Spec}{\mathop{{\bf Spec}}\nolimits}

\newcommand{\rad}{\mathop{{\rm rad}}\nolimits}

\newcommand{\imm}{\mathop{{\rm im}}\nolimits}

\newcommand{\End}{\mathop{\rm End}\nolimits}

\newcommand{\Gr}{\mathop{\rm Gr}\nolimits}
\newcommand{\Hom}{\mathop{\rm Hom}\nolimits}
\newcommand{\Stab}{\mathop{\rm Stab}\nolimits}
\newcommand{\Sym}{\mathop{\rm Sym}\nolimits}

\newcommand{\kernel}{\mathop{\rm ker}\nolimits}
\newcommand{\loccit}{[loc.$\;$cit.]}

\def\tei{\, | \,}
\def\halb{\frac{1}{2}}

\def\id{{\rm id}}

\newbox\mybox
\def\arrover#1{\mathrel{
       \setbox\mybox=\hbox spread 1.4em{\hfil$\scriptstyle#1$\hfil}
       \vbox{\offinterlineskip\copy\mybox
             \hbox to\wd\mybox{\rightarrowfill}}}}
\def\larrover#1{\mathrel{
       \setbox\mybox=\hbox spread 1.4em{\hfil$\scriptstyle#1$\hfil}
       \vbox{\offinterlineskip\copy\mybox
             \hbox to\wd\mybox{\leftarrowfill}}}}

\def\ontoover#1{\mathrel{
       \setbox\mybox=\hbox spread 1.4em{\hfil$\scriptstyle#1$\hfil}
       \vbox{\offinterlineskip\copy\mybox
             \hbox to\wd\mybox{\rightarrowfill\hskip-2.8mm
                               $\rightarrow$}}}}
\def\leftontoover#1{\mathrel{
       \setbox\mybox=\hbox spread 1.4em{\hfil$\scriptstyle#1$\hfil}
       \vbox{\offinterlineskip\copy\mybox
             \hbox to\wd\mybox{$\leftarrow$\hskip-2.8mm
                               \leftarrowfill}}}}
\def\longto{\longrightarrow}
\def\into{\hookrightarrow}
\def\onto{\ontoover{\ }}
\def\longonto{\ontoover{\ }}
\def\isoto{\arrover{\sim}}

\def\longinto{\lhook\joinrel\longrightarrow}

\usepackage[curve,matrix,arrow,cmtip]{xy}
\NoComputerModernTips

\def\myxymessage{\def\messagetext
   {Here an xy-pic diagram was omitted to speed up compilation . . . }
   \immediate\write16{\messagetext}
   \hbox{\bf \messagetext}}
\def\filxymatrix#1{\myxymessage}
\def\filxyarray#1{\myxymessage}

\newdir^{ (}{{}*!/-3pt/\dir^{(}}
\newdir_{ (}{{}*!/-3pt/\dir_{(}}
\newdir^{ )}{{}*!/+3pt/\dir^{)}}
\newdir_{ )}{{}*!/+3pt/\dir_{)}}

\def\rscript#1{\hbox to 0pt{$\scriptstyle#1$\hss}}

\let\oldbullet\bullet
\def\bullet{{\mathchoice{\oldbullet}%
                        {\oldbullet}%
                        {\scriptscriptstyle\oldbullet}%
                        {\oldbullet}}}

\newcommand{\argdot}{{\;\bullet\;}}
\newcommand{\argast}{{\;\ast\;}}

\newcommand{\bA}{\mathop{\overline{\FA}}\nolimits}
\newcommand{\uA}{\mathop{\FA^u}\nolimits}
\newcommand{\bB}{\mathop{\overline{\FB}}\nolimits}
\newcommand{\uB}{\mathop{\FB^u}\nolimits}
\newcommand{\bCU}{\mathop{\overline{\CC(U)_{w=0}}}\nolimits}

\newcommand{\bCX}{\mathop{\overline{\CC(X)_{w=0}}}\nolimits}
\newcommand{\uCX}{\mathop{\CC(X)_{w=0}^u}\nolimits}
\newcommand{\bCZ}{\mathop{\overline{\CC(Z)_{w=0}}}\nolimits}
\newcommand{\uCZ}{\mathop{\CC(Z)_{w=0}^u}\nolimits}
\newcommand{\bi}{\mathop{\overline{i}}\nolimits}
\newcommand{\ig}{\mathop{i_g}\nolimits}

\newcommand{\isp}{\mathop{i_{\sigma , p}}\nolimits}

\newcommand{\iat}{\mathop{i_\tau}\nolimits}
\newcommand{\bj}{\mathop{\overline{j}}\nolimits}
\newcommand{\jp}{\mathop{j_\varphi}\nolimits}
\newcommand{\ujast}{\mathop{j_{!*}}\nolimits}

\newcommand{\CHXM}{\mathop{CHM(X)_\BQ}\nolimits}

\newcommand{\CHFALM}{\mathop{CHM(\A_{g,L})_F}\nolimits}
\newcommand{\CHFALastM}{\mathop{CHM(\A_{g,L}^*)_F}\nolimits}
\newcommand{\CHFAtLastM}{\mathop{CHM(\A_{g,L'}^* \otimes \Spec \BZ[\frac{1}{k \cdot k'}])_F}\nolimits}
\newcommand{\CHFMLM}{\mathop{CHM(M^L)_F}\nolimits}
\newcommand{\CHFMLastM}{\mathop{CHM((M^L)^*)_F}\nolimits}
\newcommand{\CHFMtLastM}{\mathop{CHM((M^{L'})^*)_F}\nolimits}
\newcommand{\CHFSSM}{\mathop{CHM(S(\FS))_F}\nolimits}
\newcommand{\CHFUM}{\mathop{CHM(U)_F}\nolimits}
\newcommand{\CHFXM}{\mathop{CHM(X)_F}\nolimits}
\newcommand{\CHMLMs}{\mathop{CHM^s(M^L)_\BQ}\nolimits}
\newcommand{\CHFMLMs}{\mathop{CHM^s(M^L)_F}\nolimits}
\newcommand{\CHFpA}{\mathop{CHM(\partial \A_{g,L}^*)_F}\nolimits}
\newcommand{\CHFpM}{\mathop{CHM(\partial (M^L)^*)_F}\nolimits}
\newcommand{\CHFXMs}{\mathop{CHM^s(X)_F}\nolimits}
\newcommand{\CHFfXMs}{\mathop{CHM^s(X)_F^{fd}}\nolimits}
\newcommand{\CHFXoMs}{\mathop{CHM^s(X')_F}\nolimits}
\newcommand{\CHFfXoMs}{\mathop{CHM^s(X')_F^{fd}}\nolimits}

\newcommand{\CHFYpM}{\mathop{CHM(\Yp)_F}\nolimits}
\newcommand{\CHFYpMs}{\mathop{CHM^s(\Yp)_F}\nolimits}
\newcommand{\CHTFoKM}{\mathop{CHMT_{\FS}(\Spec \Fo_K)_F}\nolimits}
\newcommand{\CHTFSSM}{\mathop{CHMT_{\FS}(S(\FS))_F}\nolimits}

\newcommand{\DBcM}{\mathop{DM_{\text{\cyrb},c}}\nolimits}

\newcommand{\DBcXM}{\mathop{\DBcM(X)}\nolimits}
\newcommand{\DBcYM}{\mathop{\DBcM(Y)}\nolimits}
\newcommand{\DBcFBsM}{\mathop{\DBcM(\Bs)_F}\nolimits}
\newcommand{\DBcFSsM}{\mathop{\DBcM(\Ss)_F}\nolimits}
\newcommand{\DBcFSSM}{\mathop{\DBcM(S(\FS))_F}\nolimits}

\newcommand{\DBcFUM}{\mathop{\DBcM(U)_F}\nolimits}
\newcommand{\DBcFXM}{\mathop{\DBcM(X)_F}\nolimits}
\newcommand{\DBcFYpM}{\mathop{\DBcM(\Yp)_F}\nolimits}
\newcommand{\DBcFYPM}{\mathop{\DBcM(Y(\Phi))_F}\nolimits}

\newcommand{\DBcFYZM}{\mathop{\DBcM(Y(\Phi_Z))_F}\nolimits}
\newcommand{\DBcFZM}{\mathop{\DBcM(Z)_F}\nolimits}
\newcommand{\DFTBsM}{\mathop{DMT(\Bs)_F}\nolimits}
\newcommand{\DFTBnatsM}{\mathop{DMT(\Bs)_F^\natural}\nolimits}
\newcommand{\DFTSsM}{\mathop{DMT(\Ss)_F}\nolimits}
\newcommand{\DFTSnatsM}{\mathop{DMT(\Ss)_F^\natural}\nolimits}
\newcommand{\DFToKM}{\mathop{DMT_{\FS}(\Spec \Fo_K)_F}\nolimits}
\newcommand{\DFTSMM}{\mathop{DMT_{\FS}(M^K(\FS))_F^\natural}\nolimits}
\newcommand{\DFTSUM}{\mathop{DMT_{\FS_U}(\pi^{-1}(U))_F^\natural}\nolimits}
\newcommand{\DFTSSM}{\mathop{DMT_{\FS}(S(\FS))_F}\nolimits}
\newcommand{\DFTSnatSM}{\mathop{DMT_{\FS}(S(\FS))_F^\natural}\nolimits}
\newcommand{\DFTSnSM}{\mathop{DMT_{\FS_\alpha}(S(\FS_\alpha))_F}\nolimits}

\newcommand{\DFTSUnatSM}{\mathop{DMT_{\FS_U}(S(\FS_U))_F^\natural}\nolimits}

\newcommand{\DFTSZnatSM}{\mathop{DMT_{\FS_Z}(S(\FS_Z))_F^\natural}\nolimits}
\newcommand{\DFTXM}{\mathop{DMT(X)_F}\nolimits}
\newcommand{\MHM}{\mathop{\bf MHM}\nolimits}
\newcommand{\MHS}{\mathop{\bf MHS}\nolimits}

\newcommand{\Bs}{\mathop{B_\sigma}\nolimits}
\newcommand{\Bks}{\mathop{B_{\sigma,k}}\nolimits}
\newcommand{\Perv}{\mathop{{\bf Perv}}\nolimits}
\newcommand{\Pes}{\mathop{P_{1,[\sigma]}}\nolimits}
\newcommand{\Ss}{\mathop{S_\sigma}\nolimits}
\newcommand{\St}{\mathop{S_\tau}\nolimits}
\newcommand{\SH}{\mathop{{\bf SH}}\nolimits}
\newcommand{\bSs}{\mathop{\overline{\Ss}}\nolimits}
\newcommand{\Spi}{\mathop{M^{\pi_1 (L_1)}}\nolimits}
\newcommand{\Spti}{\mathop{M^{\pi_1 (L_1')}}\nolimits}
\newcommand{\Sps}{\mathop{M^{\pi_{[\sigma]} (K_1)}}\nolimits}
\newcommand{\bSps}{\mathop{\overline{\Sps}}\nolimits}

\newcommand{\Wp}{\mathop{W_\varphi}\nolimits}
\newcommand{\Ues}{\mathop{U_{1,[\sigma]}}\nolimits}
\newcommand{\Wes}{\mathop{W_{1,[\sigma]}}\nolimits}
\newcommand{\Xes}{\mathop{\FX_{1,[\sigma]}}\nolimits}
\newcommand{\Yp}{\mathop{Y_\varphi}\nolimits}
\newcommand{\Ykp}{\mathop{Y_{\varphi,k}}\nolimits}

\newcommand{\pis}{\mathop{\pi_\sigma}\nolimits}
\newcommand{\pios}{\mathop{\pi'_\sigma}\nolimits}
\newcommand{\pikos}{\mathop{\pi'_{\sigma,k}}\nolimits}
\newcommand{\pits}{\mathop{\pi''_\sigma}\nolimits}

\newcommand{\one}{\mathds{1}}

\begin{document}

%

\hfuzz=3pt
\overfullrule=10pt                   


\setlength{\abovedisplayskip}{6.0pt plus 3.0pt}
\setlength{\belowdisplayskip}{6.0pt plus 3.0pt}
\setlength{\abovedisplayshortskip}{6.0pt plus 3.0pt}
\setlength{\belowdisplayshortskip}{6.0pt plus 3.0pt}

\setlength{\baselineskip}{13.0pt}
\setlength{\lineskip}{0.0pt}
\setlength{\lineskiplimit}{0.0pt}

%
%

\title{Intermediate extension of Chow motives of Abelian type
\forget{
\footnotemark
\footnotetext{To appear in ....}
}
}
\author{\footnotesize by\\ \\
\mbox{\hskip-2cm
\begin{minipage}{6cm} \begin{center} \begin{tabular}{c}
J\"org Wildeshaus \footnote{
Partially supported by the \emph{Agence Nationale de la
Recherche}, project no.\ ANR-07-BLAN-0142 ``M\'ethodes \`a la
Voevodsky, motifs mixtes et G\'eom\'etrie d'Arakelov''. }\\[0.2cm]
\footnotesize Universit\'e Paris 13\\[-3pt]
\footnotesize Sorbonne Paris Cit\'e \\[-3pt]
\footnotesize LAGA, CNRS (UMR~7539)\\[-3pt]
\footnotesize F-93430 Villetaneuse\\[-3pt]
\footnotesize France\\
\footnotesize tel.\ +33-1-49403577\\
\footnotesize fax~+33-1-49403568\\
{\footnotesize \tt wildesh@math.univ-paris13.fr}
\end{tabular} \end{center} \end{minipage}
\hskip-2cm}
\\[2.5cm]
}
\date{August 22, 2016}
\maketitle
\quad \\[-1.7cm]
\begin{abstract}
\noindent 
In this article, we give an unconditional construction
of a motivic analogue of the intermediate extension
in the context of Chow motives of Abelian type.
Our main application concerns intermediate extensions 
of Chow motives associated to Kuga families
to the Baily--Borel compactification of a Shimura variety. \\


\end{abstract} 

\bigskip
\bigskip
\bigskip

\noindent {\footnotesize Math.\ Subj.\ Class.\ (2010) numbers: 14G35
(11F80, 14C25, 14F32, 14F42, 14K10, 18E05, 19E15).
}

\eject
\tableofcontents

\bigskip
\vspace*{0.5cm}

%
%

\setcounter{section}{-1}
\section{Introduction}
\label{Intro}



A profound conjecture concerning the category of motives over some
base $X$ predicts the existence of a $t$-structure, all of  whose
realizations
are compatible with the so-called \emph{perverse} $t$-structure.
This structure would in particular allow for the construction of the
\emph{intermediate extension} to $X$ of any Chow motive over an open
sub-scheme $U$ of $X$, canonically up to (unipotent) automorphisms
restricting to the identity on $U$. \\

Unfortunately, the assumption 
concerning the $t$-structure is extremely hypothetical: for
$X$ equal to the spectrum of a field of characteristic zero, 
it would imply both Grothendieck's 
standard conjectures and Murre's filtration conjecture \cite{Be}. \\

\forget{
This latter construction is carried out in \cite{CH} for quasi-projective
varieties $X$ over $\BC$, assuming Grothendieck's 
standard conjectures and Murre's filtration conjecture. Note that 
for $X$ equal to the spectrum of a field of charcteristic zero, 
these two conjectures would be implied by the existence of a $t$-structure
\cite{Be}. \\
}
The aim of this paper is to establish the theory of intermediate
extensions. Its ingredients are radically different from $t$-structures.
It is based on two key notions: \emph{weight structures}
\`a la Bondarko, and \emph{semi-primary categories} \`a la Andr\'e--Kahn. 
This approach is \emph{formula free}. Above all, it is \emph{unconditional},
once certain geometric conditions are satisfied. \\

In order to illustrate our results, let us discuss their
implications in a context which lends itself to important arithmetic applications. 
Let $M^L$ be a \emph{pure Shimura variety},
and assume that the group $L$ to which it is associated, is \emph{neat},
which implies that $M^L$ is smooth over the \emph{reflex field} $E$.
The variety $M^L$ is the target of proper, smooth morphisms
\[
\pi : M^K \longto M^L \; ,
\]
induced by a morphism of \emph{Shimura data} $(P,\FX) \to (G,\FH)$, 
which identifies $G$ with the maximal reductive quotient
of $P$. We make a mild technical 
assumption on the Shimura data $(G,\FH)$, 
namely that they satisfy \cite[Condition~(3.1.5)]{P2}
(this condition will be recalled later).
The source of $\pi$ is a \emph{Kuga variety}, \emph{i.e.}, 
a \emph{mixed Shimura variety} admitting the structure
of a torsor under an Abelian scheme over $M^L$. (Note that we admit the case
$\pi = \id_{M^L}$, \emph{i.e.}, the Abelian scheme may be of relative dimension
zero.) Fix one such $\pi$. The scheme $M^L$
being regular, and $\pi$ proper and smooth, the motive
\[
\pi_* \one_{M^K}
\]
belongs to the category $\CHMLMs$ of \emph{smooth Chow motives}
over $M^L$ \cite{L2}. Fix an extension $F$ of $\BQ$, and a direct factor
$N$ of $\pi_* \one_{M^K}$, viewed as an object of the category $\CHFMLM$
of \emph{Chow motives}. Denote by $j: M^L \into (M^L)^*$ the 
\emph{Baily--Borel compactification} \cite{AMRT,BaBo}, and by 
$i: \partial (M^L)^* \into (M^L)^*$ 
the closed immersion of the complement of $M^L$.

\begin{Thm} \label{0A}
(a)~There is a Chow motive $\ujast N \in \CHFMLastM$
extending $N$, \emph{i.e.}, $j^* \ujast N \cong N$, and satisfying
the following properties.
\begin{enumerate}
\item[(1)] $\ujast N$ admits no non-zero direct
factor belonging to $i_* \CHFpM$.
\item[(2)] Any element of the kernel of
\[
j^* : \End_{\CHFMLastM}(\ujast N) \longto \End_{\CHFMLM}(N)
\]
is nilpotent.
\end{enumerate}  
(b)~Among the extensions of $N$ to a Chow motive over $(M^L)^*$,
$\ujast N$ is characterized up to isomorphism
by each of the properties~(1) and (2). \\[0.1cm]
(c)~Any extension of $N$ to a Chow motive over $(M^L)^*$ 
is isomorphic to a direct sum
$\ujast N \oplus i_* N_\partial$,
for an object $N_\partial$ of $\CHFpM$. \\[0.1cm]
(d)~Let $f$ be an idempotent endomorphism of $N$.
Then $f$ admits an idempotent extension to $\ujast N$.
\end{Thm}

For certain Shimura varieties, parts of Theorem~\ref{0A}
are known. This concerns 
modular curves \cite{Sch}, Hilbert--Blumenthal
varieties \cite{GHM}, and Shimura varieties
of dimension $3$ \cite{MS}; note that in the latter two cases the results
are obtained only after a base change from the reflex field $E$
to the field $\BC$ of complex numbers. Over $E$ itself, 
apart from \cite{Sch}, there are the rather recent
results from \cite{V} and \cite{NV2} concerning 
Hilbert--Blumen\-thal varieties
and Siegel threefolds. We refer to Remark~\ref{8H} for details. \\

Note that the restriction $j^*$ is full on Chow motives
\cite[Thm.~1.7~(b)]{W5}. Thus, the isomorphisms in (b) and (c)
can be chosen to give the identity on $N$, when restricted to $M^L$.
The ambiguity related to the presence of a kernel of $j^*$
can be controlled by passing to a suitable quotient category of $\CHFMLastM$. 
This point of view will be developed in Section~\ref{2}; it 
also allows to treat functoriality in a satisfactory manner.  
We refer to $\ujast N$ as the \emph{intermediate extension} of $N$.
This terminology is justified thanks to our next result.

\begin{Thm} \label{0B}
The intermediate extension is compatible with the Betti and $\ell$-adic realizations
\cite{Ay2,CD2}. More precisely, denote by 
\[
R (N) \in D^b_c \bigl( M^L , F \bigr) \quad \text{and} \quad 
R (\ujast N) \in D^b_c \bigl( (M^L)^* , F \bigr)
\]
the realizations of $N$ and $\ujast N$, and by 
\[
H^n: D^b_c(\argdot,F) \longto \Perv_c (\argdot,F) \; , 
\; n \in \BZ
\]
the perverse cohomology functors. Then for any integer $n$, there is
a canonical and functorial isomorphism
\[
H^n R (\ujast N) \cong \ujast H^n R (N)
\]
of perverse sheaves on $(M^L)^*$.
\end{Thm} 

Denote by $m$ the structure morphism $(M^L)^* \to \Spec E$.
From the compatibility of the realizations with direct images \cite{Ay2,CD2},
we deduce the following.

\begin{Cor} \label{0Bcor}
Assume that $R(N)$ is concentrated in a single perverse degree $d$. 
Then $R (\ujast N)$ is concentrated in perverse degree $d$, and 
the complex computing (singular)
\emph{intersection cohomology} of $(M^L)^*$ with coefficients
in $R(N)$ is of motivic origin.
More precisely, there is a canonical isomorphism
\[
R (m_* \ujast N) \cong m_* \ujast R (N) \; .
\]
\end{Cor}

According to \cite[Thm.~3.1]{DM}, the motive
$\pi_* \one_{M^K}$ admits a Chow--K\"unneth decomposition.
Any direct factor $N$ contained in a Chow--K\"unneth component 
then satisfies the additional hypothesis of Corollary~\ref{0Bcor}.
Given its statement, 
it appears justified to think of the object $m_* \ujast N$
as being the \emph{intersection motive} of $M^L$ with coefficients in $N$
(with respect to the Baily--Borel compactification). Note that for
certain Shimura varieties, the construction of the intersection motive,
or at least partial information entering that construction,
already appears in the literature. To the cases which result from the
application of $m_*$ to the intermediate extension \cite{Sch, GHM, MS, V, NV2},
one needs to add surfaces with constant coefficients \cite{CM},
Hilbert--Blumenthal varieties with non-constant coefficients \cite{W4},
and Picard surfaces with regular coefficients \cite{W6}. 
Again, we refer to Remark~\ref{8H} for more details. \\

Now let $[\ \cdot h]: M^L \to M^{L'}$ be a finite morphism
associated to change of the ``level'' $L$. It extends to a finite
morphism $(M^L)^* \to (M^{L'})^*$, denoted by the same symbol 
$[\ \cdot h]$. Assume that $L'$ is neat, too.

\begin{Thm} \label{0C}
The intermediate extension is compatible with $[\ \cdot h]_*$. More precisely,
the Chow motive $[\ \cdot h]_* \ujast N \in \CHFMtLastM$ satisfies
the analogues of the properties~(1), (2) from Theorem~\ref{0A}~(a). Furthermore,
the analogues of Theorem~\ref{0A}~(b)--(d) hold for $[\ \cdot h]_* \ujast N$
and $[\ \cdot h]_* N$ instead of $\ujast N$ and $N$. 
\end{Thm}

Finally, let us discuss \emph{Hecke operators}. 
In order to do so, we need to fix a section $i$ 
of the projection $\pi$ from $P$ to $G$,
and to suppose that the group $K$ contains the image of $L$ under $i$.
In other words, $K$ is the semi-direct product of $L$ and the kernel
of the restriction of $\pi$ to $K$. 
Fix $x \in G(\BA_f)$, and consider the
double coset $LxL$, which is one of the generators of the 
\emph{Hecke algebra} $R(L,G(\BA_f))$
associated to $M^L$. The above results will be shown to imply the following.

\begin{Thm} \label{0D}
(a)~The Hecke operator $LxL$ acts on $m_* \ujast N$, in a way 
compatible with its action on $m_* j_! N$ and on $m_* j_* N$. \\[0.1cm]
(b)~Assume that the (Betti or $\ell$-adic) realization
$R(N)$ is concentrated in a single perverse degree. 
Then the endomorphism
\[
R (LxL) : R(m_* \ujast N) \longto R(m_* \ujast N)
\]
coincides with the Hecke operator defined on the complex computing
intersection cohomology \emph{via} the isomorphism
\[
R(m_* \ujast N) \cong m_* \ujast R(N)
\]
from Corollary~\ref{0Bcor}.
\end{Thm}

Using Faltings-Chai's results on \emph{integral models} of Baily--Borel 
and toroi\-dal compactifications of Siegel varieties $M^L$, and on integral models
of compactifications of Kuga families over $M^L$ \cite{FC}, one can show that analogues of Theorems~\ref{0A}--\ref{0D} hold for 
Chow motives occurring as direct factors in the relative motive of 
the integral models of Kuga families (Theorems~\ref{8K}--\ref{8M}).
Given Lan's recent generalization \cite{La1,La2}, it appears likely that such
analogues of Theorems~\ref{0A}--\ref{0D} actually exist for integral models
of arbitrary \emph{PEL-type Shimura varieties}. 
We refer to Remark~\ref{8I} for more precise comments. \\

The above results will be shown to result from a more general formalism.
Here are the rough ideas: (1)~apply the generalization of Kimura's notion of
\emph{finite dimensionality} \cite[D\'ef.~9.1.1]{AK} to the category
of smooth Chow motives over a base, (2)~\emph{glue} certain finite dimensional
Chow motives over strata of the base, 
and study the finiteness conditions which are respected
by the gluing process, (3)~show that the finiteness conditions 
in question are strong enough
to allow for a formulation of a theory of intermediate extensions. 
It turns out that finite dimensionality itself is in general
not preserved by gluing. This has an obvious reason: in general, gluing will not
preserve rigidity (basically because it will not preserve smoothness).
One of the main results from \cite{AK} (\loccit, Thm.~9.2.2; 
see also \cite[Lemma~4.1]{O'S}) states that a rigid category, all of whose objects
are finite dimensional, is semi-primary. It turns out that semi-primality glues
well. Actually, the only abstract ingredient one then needs in order to formulate
a theory of intermediate extension, is the weight structure on motives,
identifying Chow motives as its \emph{heart}. \\

Let us now give a more detailed description of the content of this paper.
Sections~\ref{1} and \ref{2} can be read independently of the rest of this article,
since they are of purely homological nature. From the beginning, we impose the
finiteness condition which turns out to be be ``right'' one, and thus concentrate
exclusively on semi-primary categories, \emph{i.e.}, 
additive categories $\FA$ which are semi-simple
up to ``nilpotency phenomena'', encoded in a condition on the \emph{radical} of $\FA$. 
Our main homological tool can be obtained 
directly from the definitions: according to Theorem~\ref{1C}, any morphism in
a pseudo-Abelian semi-primary category $\FA$ is a direct sum of an isomorphism and a
morphism belonging to the radical of $\FA$.
The proof of Theorem~\ref{1C} is a variation on a classical theme:
elements of a ring $A$ which are idempotent modulo
a nil-ideal $\Fn$, can be modified by an element of $\Fn$ such that
the result is idempotent in $A$ \cite[1.1.28]{R}. 
We then turn to a relative situation, and study a full inclusion $i_*: \FB \into \FA$
of one additive category into another. We anticipate one of the defining features
of the intermediate extension, and therefore
consider morphisms in $\FA$ between objects of $\FB$ on the one hand,
and objects $A$ admitting no non-zero direct factor belonging to $\FB$ on the other.  
In that context, Theorem~\ref{1C} is applied a first time. Its 
Corollary~\ref{1E} states that if $\FB$ is semi-primary, then all such morphisms
are in the radical; furthermore, $A$ can be split off any object of $\FA$
having the same image as $A$ in the quotient category $\FA / \FB$. 
Modulo a mild technical assumption, satisfied in the geometric context 
we shall consider later (and which will therefore be ignored for the purpose of this
introduction), this observation leads to a first approximation of the
intermediate extension: according to Theorem~\ref{1H}, the quotient morphism
$j^* : \FA \to \FA / \FB$ admits a partial right inverse. More precisely,
a sub-ideal $\Fg$ of the radical of $\FA$ can be identified such that $j^*$ factors through
$\FA / \Fg$, and such that the factorization admits a canonical left inverse;
it is that left inverse that will be defined as the intermediate extension. 
Let us note that 
our requirement for the categories to be pseudo-Abelian could
sometimes be dropped, provided that direct factors be replaced by
retracts. We suppose the coefficient ring $F$ of the semi-primary categories
to be a finite direct 
product of fields in order to have the results from \cite{AK}
at our disposal.\\

From Section~\ref{2} onwards, weight structures enter the picture. The
considerations of Section~\ref{1} are applied to the \emph{heart} of a
weight structure. In this context, 
Theorem~\ref{1C} is used a second time: according to Theorem~\ref{2C},
semi-primality of the heart implies the existence of \emph{minimal weight filtrations}.
Following what we consider as our main insight,
the intermediate extension of an object on an open sub-scheme $U$ of a scheme $X$
should correspond to the minimal weight filtration on the ``boundary'' object on $X - U$. 
This insight is quantified in Theorem~\ref{2K}, where we consider an abstract
gluing situation ``$X = U \coprod Z$'', equipped with weight structures. Applying Theorem~\ref{1H},
we see that the intermediate extension exists once we impose semi-primality over $Z$.
Furthermore, we show that semi-primality glues: if it holds over $Z$ and over $U$,
then it holds over $X$. This latter result opens the way to inductive constructions
of weight structures with semi-primary hearts, as they will be performed later on. 
Summary~\ref{2M} then contains all properties of the abstract intermediate extension
which we shall apply in subsequent sections. \\

Although all results in the present work are unconditional, it appeared useful
to insert a section containing two conjectures guiding our vision,
as far as motives are concerned. Section~\ref{3} starts off with a short review 
of weight structures on \emph{constructible Beilinson motives} over a base $X$ \cite{CD}. 
In particular, we recall the main results from \cite{H}, assuring the existence
of such structures. We also recall the definition of the category $CHM(X)$ of Chow motives
over $X$ as the heart of the motivic weight structure on the category
of Beilinson motives over $X$ \cite{W5}. Conjectures~\ref{3C} and \ref{3D} then
affirm that $CHM(X)$ is semi-primary, and that indecomposable Chow motives over $X$
have trivial radical. The first of these should be viewed as a 
``constructible'' analogue in the relative context of a conjecture proposed
independently by Kimura and O'Sullivan, predicting that over a field, all
Chow motives are finite dimensional. Given the theory that was developed
in Sections~\ref{1} and \ref{2}, it would imply the existence of
the intermediate extension for any base $X$. \\

The aim of Sections~\ref{4} and \ref{5} is to establish semi-primality
for certain explicit sub-categories of Chow motives. In order to do so,
we first need to study the behaviour of motives under gluing. 
The geometric context we study from now on is therefore that of a stratified
scheme $Y$, and of stra\-ti\-fied morphisms.
A first non-trivial example that comes to mind concerns Tate motives \cite{L,EL}. 
Theorem~\ref{4E} gives a sufficient criterion for the ca\-te\-gories of Tate motives
over the strata $Y_\varphi$ to glue. Thanks to \emph{absolute purity} \cite[Thm.~14.4.1]{CD},
this criterion is satisfied in particular if the closures of the strata are regular. 
For the sequel, it will be important to generalize Theorem~\ref{4E} to relative situations,
where we consider direct images of Tate motives under proper morphisms. This is
the content of Corollary~\ref{4K}. In each of these settings, the gluing is shown
to respect the weight structures; in particular, the glued categories are all
generated by their hearts. \\

In Section~\ref{5}, we put everything together, and formulate
our Main Theorem~\ref{5A}. Following the rough idea sketched further above,
we aim at gluing sub-categories of motives generated by finite dimensional Chow motives
over the strata $Y_\varphi \,$, to obtain a semi-primary category of Chow motives
over $Y$. Given that in practice, a proper morphism $\pi$ over a stratified scheme $Y$
rarely admits a nice description over the individual strata $Y_\varphi$,
unless the pre-images of the strata are stratified further, we need the assumptions
on the induced morphisms to be modified: basically, we sacrifice
properness of the strata of $\pi$
in order to obtain regularity of the strata of the source. Since Tate
twists do not affect finite dimensionality, we obtain a condition on the strata of $\pi$:
they factorize over schemes $B$ whose motive is finite dimensional over $Y_\varphi$, 
and such that the strata of the source over $B$ give rise to Tate motives. 
According to the results proved before, we get a weight structure on the
glued category, whose heart is semi-primary. The intermediate extension therefore
exists; its main properties are recalled in Corollary~\ref{5Aa}.
We conclude Section~\ref{5} with a list of examples, showing in particular that
our theory is non-empty. In particular (Example~\ref{5E}~(b), (c)), a recent generalization
of \cite[Thm.~(3.3.1)]{Kue} (see \cite[pp.~54--55, lower half of p.~61]{O'S2})
implies that finite-dimensionality is satisfied whenever $B$ is a torsor under an
Abelian scheme over $Y_\varphi$.   \\

Sections~\ref{6} and \ref{7} contain our main compatibility results concerning
the intermediate extension. Theorems~\ref{6G} and \ref{6H} are about direct images under
certain finite morphisms, and inverse images under certain smooth morphisms,
respectively. 
Theorem~\ref{7C} states that the intermediate extension is compatible with
realizations \cite{Ay2,CD2}. It shows in particular that the approach \emph{via}
weight structures is indeed the ``right'' one in the motivic context. 
All three proofs have one point in common:
the respective functor in question (direct image,
inverse image, realization) is \emph{radicial} \cite[D\'ef.~1.4.6]{AK},  
\emph{i.e.}, maps the radical of the source to the radical of the target.
In the context of Section~\ref{6}, this results from the study of the 
si\-tua\-tion
over strata --- here, we were lucky enough to have strong results from \cite{O'S2}
at our disposal --- and a rather general abstract principle:
``radiciality is compatible with gluing'' (we refer to Proposition~\ref{6E}
for a precise statement). In principle, the same strategy is used in the proof
of Theorem~\ref{7C}. However, already over strata we are confronted with a 
rather serious obstacle: it turns out that radiciality of the realization
is directly related to the conjecture ``numerical equivalence equals
homological equivalence''. This explains why Theorem~\ref{7C}
requires our most restrictive hypothesis: finite dimensionality alone does not
seem to formally imply the conjecture, we need to invoke Lieberman's result
on the validity
of the conjecture for Abelian varieties (Theorem~\ref{7K}). Once the situation over strata
is settled, there is another problem, this time related to gluing.  
The principle sketched above (Proposition~\ref{6E}) is valid only if both the
source and the target of the functor in question carry a weight structure.
But this is not the case for the bounded derived category of sheaves. 
In a sense, radiciality of the realization would thus be easier to prove 
if the Hodge theoretical realization $R_{\bf{H}}$
were known to exist, and to satisfy the full formalism of six operations. 
Note that 
\cite{I} provides a Hodge theoretical realization for schemes which
are smooth over $\BC$. For our purposes however,
it is definitely necessary to consider singular schemes. 
Similarly, we need the compatibility of $R_{\bf{H}}$
with the functors $f^*$, $f_*$, $f_!$, $f^!$.
Given the situation, the statement and the proof of 
our relevant gluing result (Corollary~\ref{7I}) is somewhat less elegant
than what one could hope. \\

In Section~\ref{8}, we show how to deduce the results stated further above 
from the general theory developed in the preceding sections. \\

This work was done while I was enjoying a \emph{d\'el\'egation
aupr\`es du CNRS}, a 
\emph{modulation de service pour les porteurs de projets de recherche},
granted by the \emph{Universit{\'e} Paris~13}, and invitations
to the \emph{Laboratoire UMPA} of the \emph{ENS Lyon}
and to the \emph{Fakult\"at f\"ur Mathematik} of the \emph{Universit\"at 
Duisburg--Essen}.
I am grateful to all four institutions.
I wish to thank G.~Ancona, S.~Bloch, M.V.~Bondarko,
D.-C.~Cisinski, F.~D\'eglise, H.~Esnault, U.~G\"ortz,
J.~Heinloth, K.~K\"unnemann, M.~Levine, A.~Nair, V.~Pilloni, 
M.~Schlichting, S.~Schr\"oer, B.~Stroh and J.~Tilouine for useful comments and discussions. 
Particular thanks go to the referee for a thorough reading and thoughtful criticism. \\

{\bf Conventions}: Throughout the article, 
$F$ denotes a finite direct product of fields
(which are supposed to be of characteristic zero from Section~\ref{3} onwards;
in other words, $F$ is then assumed to be 
a commutative semi-simple Noetherian $\BQ$-algebra).
We fix a base scheme
$\BB$, which is of finite type over some excellent scheme 
of dimension at most two. By definition, \emph{schemes} are  
$\BB$-schemes which are separated and 
of finite type (in particular, they are excellent, and
Noetherian of finite dimension), \emph{morphisms} between schemes
are separated morphisms of $\BB$-schemes, and 
a scheme is \emph{nilregular}
if the underlying reduced scheme is regular. 
A sub-category $\FB$ of an additive category $\FA$ is \emph{dense} if 
for any object $B$ of $\FB$, any direct factor of $B$ formed in $\FA$
belongs already to $\FB$. The pseudo-Abelian completion
of an additive category $\FA$ is denoted by $\FA^\natural$.


\bigskip

%
%

\section{Gluing of semi-primary categories}
\label{1}



Let us fix our coefficients $F$, and consider a
category $\FA$ which is $F$-linear \cite[Sect.~1.1]{AK}.
Recall the following definitions.

\begin{Def}[{\cite{Ke}}] \label{1A}
The \emph{radical} of $\FA$ is the ideal $\rad_\FA$
which associates to each pair of objects $A,B$ of $\FA$ the subset
\[
\rad_\FA (A,B) := \{ f \in \Hom_\FA (A,B) \; , \; 
  \forall \, g \in \Hom_\FA (B,A) \; , \; \id_A - gf \; \text{invertible} \}
\]
of $\Hom_\FA (A,B)$.
\end{Def}

In \cite[D\'ef.~1.4.1]{AK}, the radical is referred to as the \emph{Kelly radical} of $\FA$.
It can be checked that $\rad_\FA$ is indeed a two-sided ideal of $\FA$
in the sense of \cite[Sect.~1.3]{AK}, \emph{i.e.}, for each pair of objects 
$A,B$, $\rad_\FA (A,B)$ is an $F$-submodule of $\Hom_\FA (A,B)$, and for each
pair of morphisms $h: A' \to A$ and $k: B \to B'$ in $\FA$, 
\[
k \rad_\FA (A,B) h \subset \rad_\FA (A',B') \; .
\]
By definition (\emph{cf.}\ \cite[first statement of Prop.~1.4.4~b)]{AK}), 
$\rad_\FA$ is maximal among the two-sided ideals $\FI$
of $\FA$ such that the projection onto the quotient category
\[
\FA \longonto \FA / \FI
\]
is conservative.

\begin{Def}[{\cite[D\'ef.~2.3.1]{AK}}] \label{1B}
The $F$-linear category $\FA$ is called \emph{semi-primary} if \\[0.1cm]
(1)~for all objects $A$ of $\FA$, the radical $\rad_\FA (A,A)$ is nilpotent,
\emph{i.e.}, there exists a positive integer $N$ such that 
\[
\rad_\FA (A,A)^N = 0 \subset \End_\FA (A) \; ,
\]
(2)~the $F$-linear quotient category 
\[
\bA := \FA / \rad_\FA
\]
is semi-simple. 
\end{Def}

Among the basic properties of semi-primary categories,
let us mention one which we shall frequently employ: according
to \cite[Prop.~2.3.4~b)]{AK}, the category
$\bA$ is Abelian if (and only if) $\FA$ is pseudo-Abelian. \\

Here is our main homological tool; it is the key for everything to follow.

\begin{Thm} \label{1C}
Assume $\FA$ to be a pseudo-Abelian semi-primary $F$-linear category. Let
\[
f : A \longto B
\]
be a morphism in $\FA$. Then there exist decompositions
\[
A = A^r \oplus A^s \; , \; B = B^r \oplus B^s \; , 
\]  
such that \\[0.1cm]
(a)~the decompositions are respected by $f$:
\[
f = f^r \oplus f^s \in \Hom_\FA (A^r,B^r) \oplus \Hom_\FA (A^s,B^s) 
\subset \Hom_\FA (A,B) \; ,
\]
(b)~the morphism $f^r$ belongs to the radical $\rad_\FA (A^r,B^r)$, \\[0.1cm]
(c)~the morphism $f^s$ is an isomorphism. \\[0.1cm]
The isomorphism classes of $A^r$, $A^s$, $B^r$ and $B^s$ are uniquely
determined by properties (a)--(c).
\end{Thm}

\begin{Proof}
By \cite[Prop.~2.3.4~b)]{AK}, the category $\bA$ is Abelian.
Therefore, the morphism 
\[
\bar{f} := f \mod \rad_\FA
\]
admits both
a kernel and an image. In addition, $\bA$ is semi-simple, hence
direct complements $\bar{A}^s$ to $\kernel(\bar{f})$ and 
$\bar{B}^r$ to $\imm(\bar{f})$
can be chosen:
\[
A = \kernel(\bar{f}) \oplus \bar{A}^s \; , \; 
B = \bar{B}^r \oplus \imm(\bar{f}) \quad \text{in} \quad \bA \; . 
\] 
Fix such choices, and
denote by $\bar{g}$ the associated projection of $B$
onto $\imm(\bar{f})$, followed by the inverse of the restriction
of $\bar{f}$ to $\bar{A}^s$. 
Thus, both $\bar{f}\bar{g}$ and $\bar{g}\bar{f}$ are idempotent;
actually, we even have the relations 
\[
\bar{f} = \bar{f}\bar{g}\bar{f} \quad \text{and} 
\quad \bar{g} = \bar{g}\bar{f}\bar{g} \; .
\]  
Now let us study lifts $g: B \to A$ of $\bar{g}$ to $\FA$,
and show that there is a choice of lift $g$ such that both $fg$ and
$gf$ are idempotent. \emph{A priori}, for an arbitrary fixed lift $g_0$,
the difference
\[
fg_0 - fg_0fg_0
\]
is in $\rad_\FA(B,B)$. By axiom~\ref{1B}~(1),
it is therefore nilpotent. By \cite[proof of Cor.~7.8]{K},
there is a polynomial $P \in \BZ[X]$ satisfying $P(1)=1$, 
and such that the composition
\[
fg_0P(fg_0) : B \longto B
\]
is idempotent. Put
\[
g_1 := g_0P(fg_0) : B \longto A \; .
\]
Since $P(1)=1$ and $\bar{f}\bar{g}$ is idempotent, 
the polynomial expression $P(fg_0)$ lifts 
$r\cdot \id_B + (1-r)\bar{f}\bar{g}$, for some integer $r$.
Since $\bar{g}\bar{f}\bar{g} = \bar{g}$,
the morphism $g_1$ therefore lifts $\bar{g}$. 
Furthermore,
\[
fg_1 = fg_0P(fg_0)
\]    
is idempotent. Put 
\[
g:= g_1fg_1 \; .
\]
Then $g$ continues to lift $\bar{g}$. 
Let us show that
both $fg$ and $gf$ are idempotent:
\[
(fg)^2 = (fg_1)^4 = (fg_1)^2 = fg \; ,
\]
\[
(gf)^2 = (g_1f)^4 = g_1(fg_1)^3f = g_1(fg_1)f = gf \; .
\]
To finish the proof, put
\[
A^r := \kernel(gf) \quad \text{and} \quad A^s := \imm(gf) \; ,
\]
\[
B^r := \kernel(fg) \quad \text{and} \quad B^s := \imm(fg) \; ,
\]
for a choice of $g$ such that $fg$ and $gf$ are idempotent.
The morphism $f$ obviously repects these decompositions: $f = f^r \oplus f^s$. 
Modulo $\rad_\FA$, they yield the decompositions 
\[
A = \kernel(\bar{f}) \oplus \bar{A}^s \; , \; 
B = \bar{B}^r \oplus \imm(\bar{f})   
\] 
fixed further above. On these, $f^r$ induces the zero morphism,
hence $f^r$ is in the radical, while $f^s$ induces an isomorphism,
hence it is an isomorphism (recall that the projection $\FA \onto \bA$
is conservative).

This proves properties (a)--(c) for the decompositions we constructed.
Conversely, any pair of decompositions satisfying (a)--(c)
induces decompositions 
\[
A = \kernel(\bar{f}) \oplus \bar{A}^s \; , \; 
B = \bar{B}^r \oplus \imm(\bar{f}) \quad \text{in} \quad \bA \; .
\]
Thus, the isomorphism classes of the components of the
decompositions are unique in $\bA$. But then the same is true
for their isomorphism classes in $\FA$, thanks to 
fullness and conservativity of the projection $\FA \onto \bA$. 
\end{Proof}

\begin{Rem}
(a)~The decompositions
\[
A = A^r \oplus A^s \; , \; B = B^r \oplus B^s 
\]
themselves, \emph{i.e.}, the associated idempotent endomorphisms of $A$
and of $B$, should in general not be expected to be unique. \\[0.1cm]
(b)~The use of \cite[proof of Cor.~7.8]{K} in the proof of Theorem~\ref{1C} 
should be seen as a quantitative version of \cite[1.1.28]{R}:
let $R$ be a ring, and $\Fn \subset R$ a nil-ideal. Then every
idempotent element of the quotient $R / \Fn$ admits an idempotent
lift to $R$. 
\end{Rem}

Let us now fix an $F$-linear full, dense sub-category $\FB$ of $\FA$.
Denote by
\[
i_* : \FB \longinto \FA
\]
the fully faithful inclusion. Recall that by definition, 
both $\FA$ and $\FB$ being $F$-linear, they admit finite direct sums 
and products, and both notions coincide \cite[Sect.~1.1]{AK};
their formation commutes therefore with $i_*$.
Denote by 
\[
j^* : \FA \longonto \FA / \FB
\]
the projection of $\FA$ onto the quotient category $\FA / \FB$. Thus,
$j^*$ is full and essentially surjective. Two morphisms
$f_1$, $f_2$ in $\FA$ yield identical images $j^*f_1$, $j^*f_2$ 
in $\FA / \FB$ if and only if
their difference $f_2 - f_1$ factors through $i_* B$, for an object 
$B$ of $\FB$.

\begin{Cor} \label{1E}
In the above situation, 
assume $\FA$ to be pseudo-Abelian, and $\FB$ to be semi-primary.
Let $A$ be an object of $\FA$ admitting no non-zero direct factor belonging
to $\FB$. \\[0.1cm]
(a)~For any object $B$ of $\FB$, 
\[
\Hom_\FA (A,i_*B) = \rad_\FA (A,i_*B) \quad \text{and} \quad 
\Hom_\FA (i_*B,A) = \rad_\FA (i_*B,A) \; .
\]
(b)~Let $A'$ be a second object of $\FA$ 
admitting no non-zero direct factor belonging to $\FB$. Then the surjection
\[
j^* : \Hom_\FA (A',A) \longonto \Hom_{\FA / \FB} (j^* A',j^* A)
\]
detects isomorphisms. It also detects elements of the radical;
more precisely, 
\[
(j^*)^{-1} \bigl( \rad_{\FA / \FB} (j^* A', j^* A) \bigr) =
\rad_\FA (A',A) \; .
\]
(c)~Any object $A'$ of $\FA$ satisfying $j^*A' \cong j^*A$
is isomorphic to a direct sum 
$A \oplus i_* B$,
for some object $B$ of $\FB$. More precisely, given an isomorphism
\[
f : j^*A' \isoto j^*A \; ,
\]
the isomorphism 
\[
A' \isoto A \oplus i_* B
\]
can be chosen such that the associated projection onto the first factor
\[
A' \longonto A
\]
is a pre-image of $f$ under $j^*$. \\[0.1cm]
(d)~Let $A'$ be an object of $\FA$ satisfying $j^*A' \cong j^*A$,
and admitting no non-zero direct factor belonging to $\FB$. Then
$A' \cong A$. More precisely, any pre-image under $j^*$ of
an isomorphism 
$j^*A' \isoto j^*A$ is an isomorphism $A' \isoto A$.
\end{Cor}

\begin{Proof}
As for (a), let $f: A \to i_*B$ and $g: i_*B \to A$ be arbitrary morphisms
between $A$ and $i_*B$, for an object $B$ of $\FB$. It suffices to show that the
composition $fg$ belongs to the radical of $i_*B$. The category $\FB$ is
dense in $\FA$, hence it is pseudo-Abelian. Thus, Theorem~\ref{1C} 
can be applied. It yields decompositions
\[
B^r_1 \oplus B^s_1 = B = B^r_2 \oplus B^s_2 \; , 
\] 
which are respected by $fg$, such that 
\[
(fg)^r : B^r_1 \longto B^r_2
\]
belongs to the radical, and 
\[
(fg)^s : B^s_1 \isoto B^s_2 \; .
\]
In particular, $i_*B^s_1$ becomes a direct factor of $A$ \emph{via}
\[
i_*B^s_1 \longinto i_*B \stackrel{g}{\longto} A
\]
and
\[
A \stackrel{f}{\longto} i_*B \longonto i_*B^s_2 
\stackrel{((fg)^s)^{-1}}{\longto} i_*B^s_1 \; .
\]
Our assumption on the object $A$ implies that $B^s_1 = 0 = B^s_2$, hence 
\[
fg = (fg)^r \in \rad_\FB(B,B) = \rad_\FA(i_*B,i_*B) \; .
\]
Now let us turn to part~(b) of the claim. Let $f: A' \to A$ be
a morphism, such that $j^*f$ is an isomorphism. The functor
$j^*$ being full, we are reduced to assuming that $A' = A$, and
that $j^*f = \id_{j^*A}$.
Hence the difference $\id_A - f$ 
factors through an object of $\FB$.
By part~(a), this difference belongs to the radical. 
By definition of the latter, $f$ is an automorphism.
Next, assume that $g: A' \to A$ is a morphism, 
such that $j^* g$ belongs to the radical. 
For any $h: A \to A'$, the difference
$\id_{A'} - hg$ thus maps to an automorphism under $j^*$. By what we have just 
proved, it it therefore itself an automorphism. This shows that $g$
belongs to the radical. We have shown that
\[
(j^*)^{-1} \bigl( \rad_{\FA / \FB} (j^* A', j^* A) \bigr) \subset 
\rad_\FA (A',A) \; .
\]
In fact, this inclusion is an equality: $j^*$ is full, and therefore 
(\emph{cf.}\ \cite[Lemme~1.4.7]{AK}) it maps the radical of $\FA$
to the radical of $\FA / \FB$.

To prove (c), let $A'$ be in $\FA$, and assume $j^* A' \cong j^* A$.
The functor $j^*$ being full, there exist
\[
f : A \longto A' \quad \text{and} \quad g : A' \longto A
\]
such that $j^* (gf) = \id_{j^* A}$ and $j^* (fg) = \id_{j^* A'}$.
From (b), we conclude that $gf$ is an isomorphism, \emph{i.e.}, $A$ is indeed
a direct factor of $A'$:
\[
A' \isoto A \oplus B \; ,
\]
and this isomorphism maps to the identity under $j^*$. It follows that the
identity on $B$ maps to zero, meaning (since $\FB$ is dense in $\FA$)
that $B$ is an object of $\FB$.

Part (d) of the claim is a formal consequence of any of (b) or (c).
\end{Proof}

Corollary~\ref{1E}~(a) motivates the following.

\begin{Def} \label{1Ea}
(a)~Denote by $\Fg$ the two-sided ideal of $\FA$ generated by
\[
\Hom_\FA (A,i_*B) \quad \text{and} \quad 
\Hom_\FA (i_*B,A) \; ,
\]
for all objects $(A,B)$ of $\FA \times \FB$, 
such that $A$ admits no non-zero direct factor belonging to $\FB$. 
Denote by $\Fg_{\FB}$ the restriction of $\Fg$ to $\FB$. \\[0.1cm]
(b)~Denote by $\uA$ the quotient category $\FA/\Fg$, 
and by $\uB$ the quotient category $\FB/\Fg_{\FB}$.
\end{Def}

By definition, the inclusion $i_*:\FB \into \FA$ induces a
fully faithful functor $\uB \into \uA$. Slightly abusing notation,
we shall denote this functor by the same symbol $i_*$. 
The ideal $\Fg$ should be thought of as containing ``gluing data''
between $\FB$ and $\FA / \FB$. The categories $\uA$
and $\uB$ are the ``un-glued'' quotients of $\FA$ and $\FB$,
respectively. A morphism is in $\Fg$ if and only if
it is the sum of morphisms $g_n$, each of which factors through an element
of $\Hom_\FA (A_n,i_*B_n)$ or of $\Hom_\FA (i_*B_n,A_n)$, for $A_n$ and $B_n$ 
as in Definition~\ref{1Ea}~(a). In particular, $\Fg$ is contained in the
kernel of $j^*$. Under the hypotheses of Corollary~\ref{1E}, it is
also contained in $\rad_\FA$ (Corollary~\ref{1E}~(a)).
 
\begin{Cor} \label{1F}
Keep the situation of Corollary~\ref{1E}, \emph{i.e.},
assume $\FA$ to be pseudo-Abelian, $\FB$ to be semi-primary, 
and $A$ to be an object of $\FA$ admitting no non-zero direct factor belonging
to $\FB$. In addition, assume that there exists an object $B_0$ of $\FB$,
such that any endomorphism $f$ of $A$ with $j^* f = 0$
factors through $i_*B_0$. \\[0.1cm]
(a)~The kernel of 
\[
j^* : \End_\FA(A) \longonto \End_{\FA / \FB}(j^*A)
\]
equals $\Fg(A,A)$. It is 
is a nilpotent ideal contained in $\rad_\FA(A,A)$. 
In particular, any idempotent of $\End_{\FA / \FB}(j^*A)$
can be lifted to an idempotent of $\End_\FA(A)$. \\[0.1cm]
(b)~For any object $B$ of $\FB$, 
\[
\Fg (A \oplus i_*B,A \oplus i_*B)
\]
is a nilpotent ideal contained in $\rad_\FA(A \oplus i_*B,A \oplus i_*B)$. 
\end{Cor}

\begin{Proof}
Let $f$ be an endomorphism of $A$ such that $j^* f = 0$.
By our assumption, $f$
factors through $i_*B_0 \,$: there are morphisms 
$g : A \to i_*B_0$ and $h : i_*B_0 \to A$ such that $f = hg$.
Since $\FB$ is semi-primary, there is a positive integer $N$ 
such that $\rad_\FA(i_*B_0,i_*B_0)^N = \rad_\FB(B_0,B_0)^N = 0$.
Let $f_1,\ldots,f_{N+1}$ be endomorphisms of $A$ with trivial $j^* f_k$,
and choose factorizations 
\[
f_k: A \stackrel{g_k}{\longto} i_*B_0 \stackrel{h_k}{\longto} A
\]
as above, $k = 1,\ldots,N+1$. Then
\[
(g_1h_2)(g_2h_3)\cdots(g_Nh_{N+1}) = 0 : i_*B_0 \longto i_*B_0 \; .
\] 
Therefore, $f_1f_2\cdots f_{N+1} = 0: A \to A$.
This proves part~(a) of our claim. 

As for part~(b), in order to show that
\[
\Fg(A \oplus i_* B,A \oplus i_* B)
\]
is nilpotent, it suffices to show that the analogous statement 
holds for both $\Fg(A,A)$ and 
$\Fg(i_* B,i_* B)$ \cite[proof of Prop.~2.3.4~c)]{AK}.
Nilpotency of $\Fg(A,A)$ was proved in (a). For
\[
\Fg(i_* B,i_* B) = \Fg_\FB (B,B) \subset \rad_\FB (B,B) \; ,
\]
use semi-primality of $\FB$.
\end{Proof}

The hypotheses from Corollary~\ref{1F}
naturally lead to the following.

\begin{Def} \label{1G}
Let us denote the following hypothesis by $(\ast)$: 
\\[0.1cm]
Every object of $\FA / \FB$ admits a pre-image $A$
under $j^*$ without non-zero direct factors belonging to $\FB$, 
and an object $B_0$ of $\FB$,
such that any endomorphism $f$ of $A$ with $j^* f = 0$
factors through $i_*B_0$. 
\end{Def}

As we shall in Section~\ref{2}, hypothesis~$(\ast)$ is satisfied 
in the context of 
weight structures. 
Putting everything together, we get the following results.

\begin{Thm} \label{1H}
Let $i_* : \FB \into \FA$ be the inclusion of one $F$-linear 
pseudo-Abelian category
into another, with $\FB$ full and dense in $\FA$. 
Denote by $j^* : \FA \onto \FA / \FB$
the projection onto the quotient category. 
Assume hypothesis \ref{1G}~$(\ast)$, and that $\FB$ is semi-primary. \\[0.1cm]
(a)~Any object of $\FA$ is isomorphic to a direct sum $A \oplus i_* B$,
for some object $(A,B)$ of $\FA \times \FB$, 
such that $A$ admits no non-zero direct factor belonging to $\FB$. 
The object $(A,B)$ is unique
up to an isomorphism, which becomes unique in $\uA \times \uB$
(and hence in $\FA / \FB \times \uB$).
\\[0.1cm]
(b)~$\Fg$ is a nilpotent ideal contained in 
$\rad_\FA \cap \ker j^*$. \\[0.1cm]  
(c)~The $F$-linear categories $\FA / \FB$, 
$\uB$ and $\uA$ are pseudo-Abelian. \\[0.1cm]
(d)~The functors $j^*$ and $i_*$ induce a canonical
equivalence of $F$-linear categories
\[
\uA \isoto \FA / \FB \times \uB \; . 
\]
(e)~The equivalence from part~(d) descends further to induce a canonical
equivalence of $F$-linear categories
\[
\bA \isoto \overline{\FA / \FB} \times \bB \; . 
\]
(f)~If $\FB$ and $\FA / \FB$ are semi-primary, then so is $\FA$.  
\end{Thm}

\begin{Proof}
The first assertion of part~(a) follows from
hypothesis \ref{1G}~$(\ast)$ and Corollary~\ref{1E}~(c).
The second assertion of part~(a) is proved together with part~(d):
recall that by definition of $\Fg$,
\[
\Hom_{\uA}(A_1 \oplus i_* B_1,A_2 \oplus i_* B_2) 
= \Hom_{\uA}(A_1,A_2) \oplus \Hom_{\uA}(i_* B_1,i_* B_2)
\]
(for $(A_1,B_1)$ and $(A_2,B_2)$ in $\FA \times \FB$, 
such that the $A_k$ do not admit non-zero direct factors belonging to $\FB$). 
Again by definition, the right hand side equals
\[
\Hom_{\FA / \FB}(j^* A_1,j^* A_2) \oplus \Hom_{\uB}(B_1,B_2) \; .
\]

Part~(b) follows from (a), and from Corollary~\ref{1F}~(b).

The category $\FA / \FB$ is pseudo-Abelian
thanks to Corollary~\ref{1F}~(a) (take $B = 0$ in the decomposition from~(a)).
The categories $\uB$ and $\uA$ are pseudo-Abelian since by (b),
$\Fg$ is nilpotent. This proves part~(c) of the claim.

As for parts~(e) and (f), note that according to (b), 
the ideal $\Fg$ is contained in $\rad_\FA$.
It follows that $\Fg_\FB$ is contained in $\rad_\FB$.
Therefore, the category $\uB$ is semi-primary, and $i_*: \uB \into \uA$
satisfies the same hypotheses as $i_*: \FB \into \FA$.
Parts~(b) and (d) then allow us to assume that 
\[
\FA = \FA / \FB \times \FB \; ,
\]
in which case claims~(e) and (f) are trivial.
\forget{
\[
\End_{\uA}(A \oplus i_* B) = \End_{\uA}(A) \oplus \End_{\uA}(i_* B)
\]
(for $A$ and $B$ as above).
Since $i_*'$ is fully faithful, the second term is equal to $\End_{\FB'}(B)$. 
By Corollary~\ref{1E}~(b), the first term equals 
$\End_{\overline{\FA / \FB}} (j^* A)$. Thus, $\End_{\bA}(A \oplus i_* B)$
is semi-simple. In order to show that the radical
\[
\rad_\FA(A \oplus i_* B,A \oplus i_* B)
\]
is nilpotent, it suffices to show that the analogous statement 
holds for both $\rad_\FA(A,A)$ and 
$\rad_\FA(i_* B,i_* B)$ 
\cite[proof of Prop.~2.3.4~c)]{AK}. To obtain these, use 
Corollary~\ref{1F}~(c) and full faithfulness of $i_*$. 

This shows that $\FA$ is semi-primary. 
By Corollary~\ref{1F}~(a), the category $\FA / \FB$
is pseudo-Abelian.
The functors $i_*$ and $j^*$ are full, and therefore 
(\emph{cf.}\ \cite[Lemme~1.4.7]{AK}) they descend to give functors
\[
\bi_* : \bB \longinto \bA \quad \text{and} \quad 
\bj^* : \bA \longonto \overline{\FA / \FB} \; .
\] 
Since all three categories $\bA$, $\bB$ and $\overline{\FA / \FB}$,
are Abelian semi-simple \cite[Prop.~2.3.4~b)]{AK},
and $\bi_*$ identifies $\bB$ with the ``kernel'' of $\bj^*$, 
we get a cano\-ni\-cal
decomposition of the category $\bA$ (\emph{cf.}\ \cite[Lemme~2.1.3]{AK}): 
indeed, the complement to $\bB$ is the
full sub-category of objects without non-zero direct factors 
having a zero image under $\bj^*$. 
}
\end{Proof}

\begin{Rem} \label{1J}
Assume that $F$ is a finite direct product of perfect fields. 
According to \cite[D\'ef.~2.4.1]{AK}, the $F$-linear category 
$\FA$ is a \emph{Wedderburn category} 
if it is semi-primary, and if for all objects $A$ of $\FA$, 
the $F$-algebra $\End_{\bA} (A)$ is a finite $F$-module. In particular,
every object of $\bA$ is then of finite length.
Obviously, Theorem~\ref{1C}, Corollaries~\ref{1E}
and \ref{1F}, and Theorems~\ref{1H}~(a)--(e)
continue to hold when the respective hypotheses on semi-primality
are replaced by the (more restrictive) hypotheses of being Wedderburn.
It follows from the equivalence
\[
\bA \isoto \overline{\FA / \FB} \times \bB  
\]
that Theorem~\ref{1H}~(f) admits such a variant, too. More precisely,
keeping all the hypotheses of Theorem~\ref{1H} as they are, but supposing  
$\FB$ and $\FA / \FB$ to be Wedderburn categories, then 
$\FA$ is a Wedderburn category.
\end{Rem}


\bigskip
%
%

\section{The abstract theory of intermediate extensions}
\label{2}



We keep the coefficients $F$ fixed in the beginning,
and consider an $F$-linear triangulated category $\CC$.
Recall the following definitions.

\begin{Def}[{\cite[Def.~1.3.1]{Bo3}}] \label{2A}
(a)~A \emph{weight structure on $\CC$}
is a pair $w = (\CC_{w \le 0} , \CC_{w \ge 0})$ of full 
sub-categories of $\CC$, such that, putting
\[
\CC_{w \le n} := \CC_{w \le 0}[n] \quad , \quad
\CC_{w \ge n} := \CC_{w \ge 0}[n] \quad \forall \; n \in \BZ \; ,
\]
the following conditions are satisfied.
\begin{enumerate}
\item[(1)] The categories
$\CC_{w \le 0}$ and $\CC_{w \ge 0}$ are 
closed under formation of retracts in $\CC$.
\item[(2)] (Semi-invariance with respect to shifts.)
We have the inclusions
\[
\CC_{w \le 0} \subset \CC_{w \le 1} \quad , \quad
\CC_{w \ge 0} \supset \CC_{w \ge 1}
\]
of full sub-categories of $\CC$.
\item[(3)] (Orthogonality.)
For any pair of objects $X \in \CC_{w \le 0}$ and $Y \in \CC_{w \ge 1}$,
we have
\[
\Hom_{\CC}(X,Y) = 0 \; .
\]
\item[(4)] (Weight filtration.)
For any object $M \in \CC$, there exists an exact triangle
\[
X \longto M \longto Y \longto X[1]
\]
in $\CC$, such that $X \in \CC_{w \le 0}$ and $Y \in \CC_{w \ge 1}$.
\end{enumerate}
(b)~Let $w$ be a weight structure on $\CC$.
The \emph{heart of $w$} is the full $F$-linear sub-category $\CC_{w = 0}$
of $\CC$ whose objects lie 
both in $\CC_{w \le 0}$ and in $\CC_{w \ge 0}$.
\end{Def}

Slightly generalizing the terminology, for $n \in \BZ$,
we shall refer to any exact triangle
\[
X \longto M \longto Y \longto X[1]
\]
in $\CC$, with $X \in \CC_{w \le n}$ and $Y \in \CC_{w \ge n+1}$,
as a weight filtration of $M$. Note that when $\CC$ is pseudo-Abelian,
axiom~\ref{2A}~(1) says that the sub-categories
$\CC_{w \le 0}$ and $\CC_{w \ge 0}$ are dense in $\CC$, 
\emph{i.e.}, pseudo-Abelian themselves. \\

We shall be particularly interested in weight structures whose hearts
are semi-primary. In this context, Theorem~\ref{1C} will be used
as follows.

\begin{Thm} \label{2C}
Let $w$ be a weight structure on $\CC$, and assume $\CC_{w = 0}$
to be pseudo-Abelian and semi-primary. \\[0.1cm]
(a)~Any object $M$ of $\CC$
admits a \emph{minimal weight filtration, i.e.}, a weight filtration
\[
M_{\le 0} \longto M \longto M_{\ge 1} \stackrel{\delta}{\longto} M_{\le 0}[1]
\]
($M_{\le 0} \in \CC_{w \le 0}$, $M_{\ge 1} \in \CC_{w \ge 1}$) satisfying
\[
\delta \in \rad_\CC (M_{\ge 1},M_{\le 0}[1]) \; .
\]
(b)~Any two minimal
weight filtrations of the same object $M$ of $\CC$ are isomorphic,
the isomorphism being unique up to adding morphisms in $\rad_\CC$. 
The minimal weight filtration of $M$ is a direct factor
of any other weight filtration of $M$. \\[0.1cm]
(c)~Minimal weight filtrations are ``functorial up to $\rad_\CC$'':
let  
\[
M_{\le 0} \longto M \longto M_{\ge 1} \stackrel{\delta_M}{\longto} M_{\le 0}[1]
\]
and
\[
N_{\le 0} \longto N \longto N_{\ge 1} \stackrel{\delta_N}{\longto} N_{\le 0}[1]
\]
be minimal weight filtrations of objects $M$ and $N$ of $\CC$,
and $\alpha: M \to N$ a morphism. Then $\alpha$ extends to give a morphism
of exact triangles
\[
\vcenter{\xymatrix@R-10pt{
        M_{\le 0} \ar[d]_{\alpha_{\le 0}} \ar[r] &
        M \ar[d]_{\alpha} \ar[r] &
        M_{\ge 1} \ar[d]_{\alpha_{\ge 1}} \ar[r]^-{\delta_M} &
        M_{\le 0}[1] \ar[d]^{\alpha_{\le 0}[1]} \\
        N_{\le 0} \ar[r] &
        N \ar[r] &
        N_{\ge 1} \ar[r]^-{\delta_N} &
        N_{\le 0}[1]  
\\}}
\] 
and $\alpha_{\le 0}$ and $\alpha_{\ge 1}$ are unique up to 
adding morphisms in $\rad_\CC$.
\end{Thm}

\begin{Proof}
For part~(a) of the claim, start by choosing \emph{any} weight filtration
\[
X \longto M \longto Y \longto X[1]
\]
of $M$, with $X \in \CC_{w \le -1}$ and $Y \in \CC_{w \ge 0}$.
We leave it to the reader to construct from a minimal weight
filtration for $Y$ a minimal weight filtration for $M$.
Thus we may assume that $M \in \CC_{w \ge 0}$. Similarly,
starting from a weight filtration
\[
X \longto M \longto Y \longto X[1]
\]
of $M$, with $X \in \CC_{w \ge 0} \cap \CC_{w \le 1}$ and $Y \in \CC_{w \ge 2}$,
we may assume that $M \in \CC_{w \ge 0} \cap \CC_{w \le 1}$. 
This means that $M$ is actually a cone of a morphism in $\CC_{w = 0}$:
\[
A \stackrel{f}{\longto} B \longto M \longto A[1] 
\]
($A$ and $B$ in $\CC_{w = 0}$). Now apply Theorem~\ref{1C} to see that
$M$ is also the cone of a morphism $f^r$ in the radical.

Claim~(b) is a special case of claim~(c). Using orthogonality
(\emph{cf.}\ \cite[Lemma~1.4.1~1]{Bo}), it is easy to see that any morphism
$\alpha: M \to N$ can be extended to a morphism of minimal
(in fact, of arbitrary) weight filtrations. The difference
of two choices of extensions $(\alpha_{\le 0}^k,\alpha_{\ge 1}^k)$,
$k=1,2$, makes the following diagram commute.
\[
\vcenter{\xymatrix@R-10pt{
        M_{\le 0} \ar[d]_{\alpha_{\le 0}^2 - \alpha_{\le 0}^1} \ar[r] &
        M \ar[d]_{0} \ar[r] &
        M_{\ge 1} \ar[d]_{\alpha_{\ge 1}^2 - \alpha_{\ge 1}^1} 
                                                      \ar[r]^-{\delta_M} &
        M_{\le 0}[1] \ar[d]^{\alpha_{\le 0}^2[1] - \alpha_{\le 0}^1[1]} \\
        N_{\le 0} \ar[r] &
        N \ar[r] &
        N_{\ge 1} \ar[r]^-{\delta_N} &
        N_{\le 0}[1]  
\\}}
\] 
Recall that both $\delta_M$ and $\delta_N$ are in the radical, 
Thus, both differences $\alpha_{\le 0}^2 - \alpha_{\le 0}^1$ and
$\alpha_{\ge 1}^2 - \alpha_{\ge 1}^1$ factor through
a morphism in $\rad_\CC$. The radical being an ideal, 
they therefore belong to $\rad_\CC$, too.
\end{Proof}

\begin{Exs} \label{2D}
(a)~Let $w$ be a weight structure on $\CC$,
and assume in addition that $\CC$ carries a bounded $t$-structure 
$t = (\CC^{t \le 0},\CC^{t \ge 0})$, which is \emph{transversal}
to $w$ \cite[Def.~1.2.2]{Bo2}. According to 
\cite[Thm.~1.2.1~(i), Rem.~1.2.3~2.]{Bo2},
this implies in particular that the 
full sub-categories given by the intersections
\[
\CC_{w = 0}^{t = n} := \CC_{w = 0} \cap \CC^{t = n}
\]
are Abelian semi-simple, 
for all $n \in \BZ$, and that any object of $\CC_{w = 0}$
is isomorphic to a finite direct sum of objects in $\CC_{w = 0}^{t = n}$,
for varying $n$. \cite[proof of Prop.~2.3.4~c)]{AK} implies
(\emph{cf.}\ the proof of Corollary~\ref{1F}~(b))
that the heart $\CC_{w = 0}$ is semi-primary. 
\forget{
By orthogonality for the $t$-structure,
\[
\Hom_\CC (A,B) = \rad_\CC (A,B) \quad \text{and} \quad 
\Hom_\CC (B,A) = \rad_\CC (B,A) 
\]
whenever $A$ is an object of $\CC^{t = n}$,
and $B$ an object of $\CC^{t = m}$, for $n \ne m$.
(In fact, it is
a Wedderburn category 
as soon as the endomorphism ring of any object in
$\CC_{w = 0}^{t = n}$ is a finite $F$-module, for all $n \in \BZ$.) 
}
\\[0.1cm]
(b)~Let $X$ be a scheme over $\BC$, assume $F$ to be a field contained
in $\BR$, and put
$\CC := D^b ( \MHM_F X )$, the bounded derived category 
of \emph{algebraic mixed $F$-Hodge modules} on $X$ \cite[Sect.~4.2]{Sa}. 
According to \cite[Prop.~2.3.1~I]{Bo2}, $\CC$ carries a canonical
weight structure $w$: indeed, the categories $D^b ( \MHM_F X )_{w \le 0}$ and 
$D^b ( \MHM_F X )_{w \ge 0}$ are
equal to the sub-categories of complexes of
Hodge modules of weights $\le 0$ and $\ge 0$
(in the sense of \cite[Def.~4.5]{Sa}), respectively. 
Furthermore, the canonical $t$-structure
on $D^b ( \MHM_F X )$ is transversal to $w$. 
From (a), we conclude that $D^b ( \MHM_F X )_{w = 0}$ 
is semi-primary. Since $D^b ( \MHM_F X )$ is pseudo-Abelian,
Theorem~\ref{2C} can be applied, \emph{i.e.}, mi\-ni\-mal weight filtrations
exist. This seems to have non-trivial
consequences even for $X = \Spec \BC$ (in which case
$\MHM_F X = \MHS_F$, the category of 
mixed graded-polarizable $F$-Hodge structures). 
We intend to exploit 
these elsewhere. \\[0.1cm]
(c)~Minimal weight filtrations do not exist for arbitrary weight structures.
Take an Abelian $F$-linear
category $\FA$, and consider the bounded homotopy category
$\CC = K^b(\FA)$ with the ``stupid'' weight structure
\[
w = \bigl( K^{b,\ge 0}(\FA) , K^{b,\le 0}(\FA) \bigr)
\]
\cite[p.~400]{Bo}. We leave it to the reader to verify that the existence of
minimal weight filtrations in $\CC$ implies (and is in fact equivalent to) 
the validity of the conclusion of Theorem~\ref{1C}
for any morphism $f$ in $\FA$. 
Take $\FA :=$ the category of $F[X]$-modules
and $f :=$ multiplication by $X$ on $F[X]$ to see that this conclusion
is not valid for general Abelian $F$-linear categories. 
(According to \cite[Prop.~2.3.5]{AK}, it \emph{is} valid for Abelian $F$-linear
categories in which every object is of finite length.) 
\end{Exs}

\begin{Rem} \label{2Da}
In our previous work \cite{W3}, we followed a different approach to minimal weight
filtrations. Instead of deducing their existence from abstract properties of the
category $\CC$, we basically restricted our attention to objects $M$ admitting
a minimal weight filtration for an obvious reason. Indeed, for integers $m \le n$,
$M$ is said to 
admit a weight filtration \emph{avoiding weights $m,\ldots,n$} \cite[Def.~1.6]{W3}
if there exists an exact triangle
\[
M_{\le m-1} \longto M \longto M_{\ge n+1} \stackrel{\delta}{\longto} M_{\le m-1}[1]
\]
in $\CC$, with $M_{\le m-1} \in \CC_{w \le m-1}$ and $M_{\ge n+1} \in \CC_{w \ge n+1}$. 
The integer $m = (m-1)+1$ being strictly smaller than $n+1$, orthogonality
(condition~(3) in Definition~\ref{2A}~(a)) tells us that there 
are no non-zero morphisms
from $M_{\le m-1}[1]$ to $M_{\ge n+1}$. Therefore, the morphism
$\delta$ is necessarily in the radical of $\CC$. It follows that for objects $M$ of $\CC$
admitting a weight filtration avoiding weight $0$ or weight $1$, that weight filtration
is minimal in the sense of Theorem~\ref{2C}~(a). Note that in this case, there are proper
unicity and functoriality statements concerning the minimal weight filtrations 
\cite[Prop.~1.7]{W3}, \emph{i.e.}, the 
passage to the quotient by $\rad_\CC$ as in 
Theorem~\ref{2C}~(b) and (c) is not needed.
The draw-back of the method from \cite{W3} is that in the applications to motives
(see \emph{e.g.}\ \cite{W4}),
the condition on absence of weights is not necessarily easy to verify.
\end{Rem}

For the rest of this section, let us assume that our $F$-linear
triangulated category is obtained
by \emph{gluing}; modifying slightly our notation, let us fix
three $F$-linear triangulated 
categories $\CC(U)$, $\CC(X)$ and $\CC(Z)$,
together with six exact functors
\[
\vcenter{\xymatrix@R-10pt{
        \CC(U) \ar[r]^-{j_!} & \CC(X) \ar[r]^-{i^*} & \CC(Z) \\
        \CC(U) & \CC(X) \ar[l]_-{j^*} & \CC(Z) \ar[l]_-{i_*} \\ 
        \CC(U) \ar[r]^-{j_*} & \CC(X) \ar[r]^-{i^!} & \CC(Z)
\\}}
\] 
satisfying the axioms from \cite[Sect.~1.4.3]{BBD}. We assume that
$\CC(U)$, $\CC(X)$ and $\CC(Z)$ are equipped with weight structures
$w$ (the same letter for the three weight structures), and that
the one on $\CC(X)$ is actually obtained from the two others in a
way \emph{compatible with the gluing}, meaning 
that the left adjoints $j_!$, $j^*$, $i^*$
and $i_*$ respect the categories $\CC(\bullet)_{w \le 0}$, and the right
adjoints $j^*$, $j_*$, $i_*$ and $i^!$ respect the categories
$\CC(\bullet)_{w \ge 0}$. In particular, we have a fully faithful functor
\[
i_* : \CC(Z)_{w = 0} \longto \CC(X)_{w = 0} 
\]
compatible with the formation of finite direct sums, and
with dense image. We also have an exact functor 
\[
j^* : \CC(X)_{w = 0} \longto \CC(U)_{w = 0} \; .
\]

\begin{Prop} \label{2E}
The functor
\[
j^* : \CC(X)_{w = 0} \longto \CC(U)_{w = 0} 
\]
is full and essentially surjective.
\end{Prop}

\begin{Proof}
\forget{
We first establish essential surjectivity.
Let $M_U$ be an object of $\CC(U)_{w = 0}\, $, and consider the morphism
\[
m: j_! M_U \longto j_* M_U
\]
in $\CC(X)$ \cite[Axiom~(1.4.3.2)]{BBD}. 
Applying $j^*$ to $m$ yields an isomorphism \cite[Axiom~(1.4.3.5)]{BBD}.
Therefore, by localization \cite[Axiom~(1.4.3.4)]{BBD}, 
any cone of $m$ is of the form $i_* C$,
for an object $C$ of $\CC(Z)$. 
Choose and fix such a cone $i_* C$, as well as a weight filtration
\[
C_{\le 0} \stackrel{c_-}{\longto} C 
\stackrel{c_+}{\longto} C_{\ge 1} 
\stackrel{\delta}{\longto} C_{\le 0}[1]
\]
of $C \in \CC(Z)$. 
Thus, 
\[
C_{\le 0} \in \CC(Z)_{w \le 0} \quad \text{and} \quad 
C_{\ge 1} \in \CC(Z)_{w \ge 1} \; .
\] 
According to axiom TR4' of triangulated categories (see \cite[Sect.~1.1.6]{BBD}
for an equivalent formulation), the diagram of exact triangles
\[
\vcenter{\xymatrix@R-10pt{
        0 \ar[d] \ar[r] &
        i_* C_{\ge 1}[-1] \ar@{=}[r] &
        i_* C_{\ge 1}[-1] \ar[d]^{i_* \delta[-1]} \ar[r] &
        0 \ar[d] \\
        j_! M_U \ar@{=}[d] &  
         &
        i_* C_{\le 0} \ar[d]^{i_* c_-} \ar[r] &
        j_! M_U[1] \ar@{=}[d] \\
        j_! M_U \ar[d] \ar[r]^m &
        j_* M_U \ar[d] \ar[r] &
        i_* C \ar[d]^{i_* c_+} \ar[r] &
        j_! M_U[1] \ar[d] \\
        0 \ar[r] &
        i_* C_{\ge 1} \ar@{=}[r] &
        i_* C_{\ge 1} \ar[r] &
        0    
\\}}
\]
in $\CC(X)$ can be completed to give
\[
\vcenter{\xymatrix@R-10pt{
        0 \ar[d] \ar[r] &
        i_* C_{\ge 1}[-1] \ar[d] \ar@{=}[r] &
        i_* C_{\ge 1}[-1] \ar[d]^{i_* \delta[-1]} \ar[r] &
        0 \ar[d] \\
        j_! M_U \ar@{=}[d] \ar[r] &  
        M \ar[d] \ar[r] &
        i_* C_{\le 0} \ar[d]^{i_* c_-} \ar[r] &
        j_! M_U[1] \ar@{=}[d] \\
        j_! M_U \ar[d] \ar[r]^m &
        j_* M_U \ar[d] \ar[r] &
        i_* C \ar[d]^{i_* c_+} \ar[r] &
        j_! M_U[1] \ar[d] \\
        0 \ar[r] &
        i_* C_{\ge 1} \ar@{=}[r] &
        i_* C_{\ge 1} \ar[r] &
        0    
\\}}
\] 
with $M \in \CC(X)$. Since the composition of functors $j^* i_*$ 
is trivial \cite[Axiom~(1.4.3.3)]{BBD}, the inverse image $j^* M$ is
isomorphic to $M_U$. Now recall that the weight structures 
are supposed to be compatible with the gluing:
the functors $i_! = i_*$ 
and $j_!$ respect $\CC(\bullet)_{w \le 0}$, 
and $i_*$ and $j_*$ respect $\CC(\bullet)_{w \ge 0}$. Therefore, 
by the above diagram, 
the object $M$
is simultaneously an extension of objects of weights $\le 0$, 
and an extension of objects of weights $\ge 0$.
It follows easily (\emph{cf.}\ \cite[Prop.~1.3.3~3]{Bo})
that $M$ belongs to $\CC(X)_{w = 0}$.

In order to prove that $j^*$ is full, 
let $M$ and $N$ be objects of $\CC(X)_{w = 0}$, and assume that 
a morphism
\[
\beta_U : j^* M \longto j^* N
\]
between their restrictions to $U$ is given.
Consider the localization triangles for
$M$ and for $N$.
\[
\vcenter{\xymatrix@R-10pt{
        i_* i^* M [-1] \ar[r] &
        j_! j^* M \ar[d]^{j_! \beta_U} \ar[r] &
        M \ar[r] &
        i_* i^* M \\
        i_* i^* N [-1] \ar[r] &
        j_! j^* N \ar[r] &
        N \ar[r] &
        i_* i^* N  
\\}}
\] 
Thanks to compatibility of the weight structures with the gluing, 
these triangles are weight
filtrations of $j_! j^* M$ and of $j_! j^* N$, respectively. By orthogonality
(condition~(3) in Definition~\ref{2A}~(a)), any morphism
from $i_* i^* M [-1]$ to $N$ is zero. Therefore, the above 
diagram can be completed to give a morphism of exact triangles. 
}
Imitate the proof of \cite[Thm.~1.7]{W5}.
\end{Proof}

\begin{Cor} \label{2F}
The functor 
\[
j^* : \CC(X)_{w = 0} \longonto \CC(U)_{w = 0} 
\]
identifies $\CC(U)_{w = 0}$ with the categorical quotient
of $\CC(X)_{w = 0}$ by $\CC(Z)_{w = 0}$.
\end{Cor}

\begin{Proof}
By \cite[Prop.~1.4.5]{BBD}, the functor $j^*: \CC(X) \onto \CC(U)$
induces a canonical equivalence
\[
\CC(X) / \CC(Z) \isoto \CC(U) \; .
\] 
Since $\CC(X)_{w = 0} \cap i_* \CC(Z) = i_* \CC(Z)_{w = 0}$, this implies
that $\CC(X)_{w = 0} / \CC(Z)_{w = 0}$ becomes 
a full sub-category of $\CC(U)$ \emph{via} $j^*$. 
Now apply Proposition~\ref{2E}.
\end{Proof}

We are therefore in the context studied in Section~\ref{1}, with
\[
\FB = \CC(Z)_{w = 0} \; , \; \FA = \CC(X)_{w = 0} \; , \; \text{and} \;
\FA / \FB = \CC(U)_{w = 0} \; .
\]

\begin{Lem} \label{2G}
(a)~Let $M \in \CC(X)_{w = 0}$. Assume that 
the composition of the two adjunction morphisms
\[
i_* i^! M \longto M \longto i_* i^* M
\]
belongs to the radical 
\[
\rad_{\CC(X)}(i_* i^! M , i_* i^* M)
= \rad_{\CC(Z)}(i^! M , i^* M) \; .
\]
Then $M$ admits no non-zero direct factor belonging to 
$i_* \CC(Z)_{w = 0}$. \\[0.1cm]
(b)~Let $M \in \CC(X)_{w = 0}$. Then there exists 
an object $N_0$ of $\CC(Z)_{w = 0}$, 
such that any endomorphism $f$ of $M$ with $j^* f = 0$
factors through $i_*N_0$. \\[0.1cm]
(c)~Assume $\CC(Z)_{w = 0}$ to be pseudo-Abelian and semi-primary. Then
any object $M_U$ of $\CC(U)_{w = 0}$ admits a pre-image $M$
under $j^*$ such that the composition of the two adjunction morphisms
\[
i_* i^! M \longto M \longto i_* i^* M
\]
belongs to the radical. \\[0.1cm]
(d)~Assume $\CC(X)_{w = 0}$ and $\CC(Z)_{w = 0}$ to be pseudo-Abelian.
In addition, assume $\CC(Z)_{w = 0}$ to be semi-primary. Then 
hypothesis \ref{1G}~$(\ast)$ is satisfied for
\[
\FB = \CC(Z)_{w = 0} \; , \; \FA = \CC(X)_{w = 0} \; , \; \text{and} \;
\FA / \FB = \CC(U)_{w = 0} \; .
\] 
\end{Lem}

\begin{Proof}
(a): Assume that some composition of morphisms
\[
i_* N \longto M \longto i_* N 
\]
gives the identity on $i_* N$,
for some object $N$ of $\CC(Z)_{w = 0}$. Then
the adjunction properties of $i^*$, $i_*$ and $i^!$ 
show that $N$ is a direct factor of both
$i^! M$ and $i^* M$, 
and that the restriction of the composition $i^! M \to i^* M$
of the adjunction morphisms to this direct factor is the identity. 
But by assumption, this composition belongs
to $\rad_{\CC(Z)} (i^! M , i^* M)$, hence $N = 0$.

(b): Choose a weight filtration
\[
X \longto i^* M \longto N_0 \longto X[1]
\]
in $\CC(Z)$, with $X \in \CC(Z)_{w \le -1}$ and $N_0 \in \CC(Z)_{w = 0}$.
Any endomorphism $f$ of $M$ restricting trivially to $\CC(U)_{w = 0}$
factors through 
the adjunction morphism $M \to i_* i^* M$. Thanks to  
orthogonality, $f$ actually factors through $i_* N_0$. 

(c): Let $M_U$ be an object of $\CC(U)_{w = 0}$. 
Consider the morphism
\[
m: j_! M_U \longto j_* M_U
\]
in $\CC(X)$, and fix a cone $i_* C$ of $m$.
Any weight filtration
\[
C_{\le 0} \stackrel{c_-}{\longto} C 
\stackrel{c_+}{\longto} C_{\ge 1} 
\stackrel{\delta}{\longto} C_{\le 0}[1]
\]
of $C$ yields a diagram of exact triangles
\[
\vcenter{\xymatrix@R-10pt{
        0 \ar[d] \ar[r] &
        i_* C_{\ge 1}[-1] \ar[d] \ar@{=}[r] &
        i_* C_{\ge 1}[-1] \ar[d]^{i_* \delta[-1]} \ar[r] &
        0 \ar[d] \\
        j_! M_U \ar@{=}[d] \ar[r] &  
        M \ar[d] \ar[r] &
        i_* C_{\le 0} \ar[d]^{i_* c_-} \ar[r] &
        j_! M_U[1] \ar@{=}[d] \\
        j_! M_U \ar[d] \ar[r]^m &
        j_* M_U \ar[d] \ar[r] &
        i_* C \ar[d]^{i_* c_+} \ar[r] &
        j_! M_U[1] \ar[d] \\
        0 \ar[r] &
        i_* C_{\ge 1} \ar@{=}[r] &
        i_* C_{\ge 1} \ar[r] &
        0    
\\}}
\] 
with $M \in \CC(X)_{w = 0}$. 
This is true in particular when the weight filtration of $C$
is chosen as the minimal weight filtration of Theorem~\ref{2C}.
In this case, the morphism 
\[
\delta[-1]: C_{\ge 1}[-1] \longto C_{\le 0}
\]
belongs to the radical. But by the above diagram, this morphism
is isomorphic to the composition of the two adjunction morphisms
\[
i_* i^! M \longto M \longto i_* i^* M \; .
\]
(d): This follows from parts (c), (a) and (b).
\end{Proof}

\begin{Rem} \label{2Ga}
(a)~As the proof of Lemma~\ref{2G}~(c) shows, there is a quantitative
version of essential surjectivity of 
\[
j^* : \CC(X)_{w = 0} \longto \CC(U)_{w = 0} \; .
\]
Indeed, let $M_U$ be an object of $\CC(U)_{w = 0}$, and 
\[
m: j_! M_U \longto j_* M_U
\]
the morphism from \cite[Axiom~(1.4.3.2)]{BBD}. Choose and fix a cone
$i_* C$ of $m$. Then the map
\[
\big\{ ( C_{\le 0} , C_{\ge 1} ) \big\} / \cong \;
\longto 
\big\{ M \big\} / \cong
\]
is a bijection between \\[0.1cm]
(1)~the isomorphism classes of weight filtrations of $C$, \\[0.1cm]
(2)~the isomorphism classes of pre-images $M$ of $M_U$ under $j^*$. \\[0.1cm]
The inverse bijection maps the class of $M$ to the class of
$( i^* M , i^! M[1])$. \\[0.1cm]
(b)~In the situation of (a), assume $\CC(Z)_{w = 0}$ to be semi-primary.
As the proof of Lemma~\ref{2G} shows, the bijection 
\[
\big\{ ( C_{\le 0} , C_{\ge 1} ) \big\} / \cong \;
\longto 
\big\{ M \big\} / \cong
\]
from (a)
maps the class of the minimal weight filtration of Theorem~\ref{2C}
to the class of an object admitting 
no non-zero direct factor belonging to $i_* \CC(Z)_{w = 0} \,$;
note that according to Corollary~\ref{1E}~(d), this class is unique.
\end{Rem}

\begin{Thm} \label{2K}
Let the $F$-linear pseudo-Abelian
triangulated categories $\CC(U)$, $\CC(X)$ and $\CC(Z)$
be related by gluing, and equipped with weight structures $w$
compatible with the gluing. \\[0.1cm]
(a)~If $\CC(Z)_{w = 0}$ is semi-primary, then the conclusions
of Theorem~\ref{1H}~(a)--(e) hold for 
\[
\FB = \CC(Z)_{w = 0} \; , \; \FA = \CC(X)_{w = 0} \; , \; \text{and} \;
\FA / \FB = \CC(U)_{w = 0} \; .
\] 
In particular, the functors $j^*$ and $i_*$ then induce canonical
equivalences of $F$-linear categories
\[
\uCX \isoto \CC(U)_{w = 0} \times \uCZ 
\]
and
\[
\bCX \isoto \bCU \times \bCZ \; . 
\]
(b)~If both $\CC(Z)_{w = 0}$ and $\CC(U)_{w = 0}$ are semi-primary,
then so is $\CC(X)_{w = 0}$. If both $\CC(Z)_{w = 0}$ and $\CC(U)_{w = 0}$ 
are Wedderburn categories, then so is $\CC(X)_{w = 0}$.
\end{Thm}

\forget{
Proceeding as in Section~\ref{1}, denote by $\Fg$ 
the ideal generated by
\[
\Hom_{\CC(X)_{w = 0}} (M,i_*N) \quad \text{and} \quad 
\Hom_{\CC(X)_{w = 0}} (i_*N,M) \; ,
\]
for all objects $(M,N)$ of $\CC(X)_{w = 0} \times \CC(Z)_{w = 0}$, 
such that $M$ admits no non-zero direct factor belonging to $\CC(Z)_{w = 0}$. 
Put 
\[
\uCX := \CC(X)_{w = 0}/\Fg  \quad \text{and} \quad
\uCZ := \CC(Z)_{w = 0}/ \bigl( \Fg \cap \CC(Z)_{w = 0} \bigr) \; . 
\]
In the context of weight structures, Theorems~\ref{1H} and \ref{1I}
take on the following form.

\begin{Thm} \label{2H}
Let the $F$-linear pseudo-Abelian
triangulated categories $\CC(U)$, $\CC(X)$ and $\CC(Z)$
be related by gluing, and equipped with weight structures $w$
compatible with the gluing. Assume 
$\CC(Z)_{w = 0}$ to be semi-primary. \\[0.1cm]
(a)~Any object of $\CC(X)_{w = 0}$ is isomorphic to a direct sum 
$M \oplus i_* N$,
for some object $(M,N)$ of $\CC(X)_{w = 0} \times \CC(Z)_{w = 0}$, 
such that $M$ admits no non-zero direct factor belonging to $\CC(Z)_{w = 0}$. 
The object $(M,N)$ is unique
up to an isomorphism, which becomes unique in $\uCX \times \uCZ$
(and hence in $\CC(U)_{w = 0} \times \uCZ$). \\[0.1cm]
(b)~$\Fg$ is a nilpotent ideal contained in 
$\rad_{\CC(X)_{w = 0}} \cap \ker j^*$. \\[0.1cm]  
(c)~The $F$-linear categories $\uCX$ and $\uCZ$ are pseudo-Abelian. \\[0.1cm]
(d)~The functors $j^*$ and $i_*$ induce a canonical
equivalence of $F$-linear categories
\[
\uCX \isoto \CC(U)_{w = 0} \times \uCZ \; . 
\]
\end{Thm}

\begin{Thm} \label{2I}
Keep the situation of Theorem~\ref{2H}. In addition, assume
$\CC(Z)_{w = 0}$ and $\CC(U)_{w = 0}$ to be semi-primary. \\[0.1cm]
(a)~The category $\CC(X)_{w = 0}$ is semi-primary. \\[0.1cm] 
(b)~The equivalence from Theorem~\ref{2H}~(d) 
descends further to induce a canonical
equivalence of $F$-linear Abelian semi-simple categories
\[
\bCX \isoto \bCU \times \bCZ \; . 
\]
\end{Thm}
}
\begin{Proof}
By Corollary~\ref{2F} and
Lemma~\ref{2G}~(d), the hypotheses of Theorem~\ref{1H}
are fulfilled.
\end{Proof}

\begin{Def} \label{2J} 
Let the $F$-linear pseudo-Abelian
triangulated catego\-ries $\CC(U)$, $\CC(X)$ and $\CC(Z)$
be related by gluing, and equipped with weight structures $w$
compatible with the gluing. 
Assume $\CC(Z)_{w = 0}$ to be semi-primary. 
Define the \emph{intermediate extension}
\[
\ujast : \CC(U)_{w = 0} \longinto \uCX
\] 
as the fully faithful functor corresponding to $(\id_{\CC(U)_{w = 0}},0)$
under the equivalence
$\uCX \isoto \CC(U)_{w = 0} \times \uCZ$
of Theorem~\ref{2K}~(a). 
\end{Def}

\begin{Rem}
According to the one of the main results of \cite{AK}, the categorical
version of Wedderburn's theorem holds \cite[Thm.~12.1.1]{AK}.
Thus, assuming that $\CC(Z)_{w = 0}$ and $\CC(U)_{w = 0}$ 
are Wedderburn categories (hence so is $\CC(X)_{w = 0}$), there exist sections
\[
\iota : \bCX \longinto \CC(X)_{w = 0} 
\]
of the projection $pr_{\CC(X)_{w = 0}}$ from $\CC(X)_{w = 0}$ to $\bCX$. 
In an earlier approach, we defined  
an intermediate extension to be a functor
\[
j_{!*} : \CC(U)_{w = 0} \longto \CC(X)_{w = 0}
\]
of the form $\iota \circ \bj_{!*} \circ pr_{\CC(U)_{w = 0}}$
\[
\vcenter{\xymatrix@R-10pt{
    \CC(U)_{w = 0} \ar@{->>}[d]_{pr_{\CC(U)_{w = 0}}} \ar@{-->}[r]^-{j_{!*}} &
    \CC(X)_{w = 0}  \\
    \bCU \ar@{^{ (}->}[r]^-{\bj_{!*}} &
    \bCX \ar[u]_{\iota} 
\\}}
\] 
where $\iota$ is a choice of such a section, and $\bj_{!*}$ is defined
as in Definition~\ref{2J}. The price to pay for 
a functor with target $\CC(X)_{w = 0}$ is that the relation
$j^* \circ j_{!*} \cong \id_{\CC(U)_{w = 0}}$ does not necessarily hold;
\emph{a priori}, it is true only up to elements in the 
radical of $\CC(U)_{w = 0}$.
\end{Rem} 

Let us conclude the section by listing the properties of $j_{!*}$.

\begin{Sum} \label{2M}
Let the $F$-linear pseudo-Abelian
triangulated catego\-ries $\CC(U)$, $\CC(X)$ and $\CC(Z)$
be related by gluing, and equipped with weight structures $w$
compatible with the gluing. 
Assume $\CC(Z)_{w = 0}$ to be semi-primary. \\[0.1cm]
(a)~Let $M_U$ be an object of $\CC(U)_{w = 0}$. Then
$j_{!*} M_U$ extends $M_U$, \emph{i.e.}, 
$j^* j_{!*} M_U \cong M_U$. 
Furthermore, $j_{!*} M_U$ satisfies the following properties,
any of which characterizes $j_{!*} M_U$ 
among the objects $M$ of $\CC(X)_{w = 0}$ extending $M_U$,
up to an iso\-morphism, which becomes unique in $\uCX$ (and hence in $\bCX$).
\begin{enumerate}
\item[(1)] $j_{!*} M_U$ admits no non-zero direct
factor belonging to $i_* \CC(Z)_{w = 0}$. 
\item[(2a)] For all objects $N$ of $\CC(Z)_{w = 0}$,
\[
\Hom_{\CC(X)_{w = 0}} (j_{!*} M_U,i_*N) 
= \rad_{\CC(X)_{w = 0}} (j_{!*} M_U,i_*N) \; .
\] 
\item[(2b)] For all objects $N$ of $\CC(Z)_{w = 0}$, 
\[
\Hom_{\CC(X)_{w = 0}} (i_*N,j_{!*} M_U) 
= \rad_{\CC(X)_{w = 0}} (i_*N,j_{!*} M_U) \; .
\]
\item[(3a)] The composition of the two adjunction morphisms
\[
i_* i^! j_{!*} M_U \longto j_{!*} M_U \longto i_* i^* j_{!*} M_U
\]
belongs to the radical 
\[
\rad_{\CC(X)}(i_* i^! j_{!*} M_U , i_* i^* j_{!*} M_U)
= \rad_{\CC(Z)}(i^! j_{!*} M_U , i^* j_{!*} M_U) \; .
\]
\item[(3b)] There exist weight filtrations
\[
N_{\le -1} \longto i^* j_{!*} M_U \stackrel{m_+}{\longto} N_0 
\longto N_{\le -1} [1]
\]
of $i^* j_{!*} M_U$ (with $N_{\le -1} \in \CC(Z)_{w \le -1}$ and 
$N_0 \in \CC(Z)_{w = 0}$), and
\[
L_0 \stackrel{m_-}{\longto} i^! j_{!*} M_U \longto L_{\ge 1} 
\longto L_0 [1]
\]
of $i^! j_{!*} M_U$ (with $L_0 \in \CC(Z)_{w = 0}$ and 
$L_{\ge 1} \in \CC(Z)_{w \ge 1}$), such that
the composition
\[
f: i_* L_0 \stackrel{i_*m_-}{\longto} i_*i^! j_{!*} M_U 
\longto i_*i^* j_{!*} M_U \stackrel{i_*m_+}{\longto} i_* N_0
\] 
of $i_* m_-$, of the composition of
the adjunction morphisms, and of $i_* m_+$ belongs to the radical:
\[
f \in \rad_{\CC(X)_{w = 0}} (i_*L_0,i_*N_0) 
= \rad_{\CC(Z)_{w = 0}}(L_0 , N_0) \; .
\] 
\item[(3c)] For all weight filtrations
\[
N_{\le -1} \longto i^* j_{!*} M_U \stackrel{m_+}{\longto} N_0 
\longto N_{\le -1} [1]
\]
of $i^* j_{!*} M_U$, and
\[
L_0 \stackrel{m_-}{\longto} i^! j_{!*} M_U \longto L_{\ge 1} 
\longto L_0 [1]
\]
of $i^! j_{!*} M_U$, 
the composition
\[
f: i_* L_0 \stackrel{i_*m_-}{\longto} i_*i^! j_{!*} M_U 
\longto i_*i^* j_{!*} M_U \stackrel{i_*m_+}{\longto} i_* N_0
\] 
belongs to the radical:
\[
f \in \rad_{\CC(X)_{w = 0}} (i_*L_0,i_*N_0) 
= \rad_{\CC(Z)_{w = 0}}(L_0 , N_0) \; .
\] 
\item[(4a)] The kernel of 
\[
j^* : \End_{\CC(X)_{w = 0}}(j_{!*} M_U) \longonto \End_{\CC(U)_{w = 0}}(M_U)
\]
is contained in $\rad_{\CC(X)_{w = 0}} (j_{!*} M_U,j_{!*} M_U)$. 
\item[(4b)] Any element of the kernel of
\[
j^* : \End_{\CC(X)_{w = 0}}(j_{!*} M_U) \longonto \End_{\CC(U)_{w = 0}}(M_U)
\]
is nilpotent.
\end{enumerate}
(b)~Any object $M$ of $\CC(X)_{w = 0}$ is isomorphic to a direct sum
$j_{!*} M_U \oplus i_* N$,
for an object $M_U$ of $\CC(U)_{w = 0}$ and an object $N$ of $\CC(Z)_{w = 0}$.
$M_U$ is unique up to unique isomorphism; in fact, $M_U \cong j^* M$.
$N$ is unique up to an isomorphism, which becomes unique in $\uCZ$
(and hence in $\bCZ$).  \\[0.1cm] 
(c)~The indecomposable objects of $\CC(X)_{w = 0}$ are precisely those 
of the form 
\[
j_{!*} M_U \; ,
\]
for an indecomposable object $M_U$ of $\CC(U)_{w = 0}$, or of the form
\[
i_* N \; ,
\]
for an indecomposable object $N$ of $\CC(Z)_{w = 0}$. \\[0.1cm]
(d)~Let $f_U : M_U' \to M_U$ be a morphism in $\CC(U)_{w = 0}$. Then any extension of $f_U$ to a morphism 
\[
f: j_{!*} M_U' \longto j_{!*} M_U
\]
in $\CC(X)_{w = 0}$,
\emph{i.e.}, satisfying $j^* f = f_U$, 
maps to $j_{!*} f_U$ in the quotient category $\CC(X)_{w = 0}^u$. \\[0.1cm]
(e)~Let $f_U$ be an idempotent endomorphism of an object of $\CC(U)_{w = 0}$.
Then $j_{!*} f_U \in \CC(X)_{w = 0}^u$ can be lifted idempotently
to $\CC(X)_{w = 0}$.
\end{Sum}

\begin{Proof}
By definition, property (1) of part~(a) 
characterizes and is satisfied by $j_{!*} M_U$.

Part~(b) follows from Theorem~\ref{1H}~(a), part~(d) 
from the definitions, and part~(e) from Theorem~\ref{1H}~(b).

Part~(c) is implied by (b) and property~(1) of (a).

The object $j_{!*} M_U$ satisfies properties~(2a)--(4b) of (a): use
Theorem~\ref{1H}~(b) and Lemma~\ref{2G} (see Remark~\ref{2Ga}).  
Now let $M$ be an object of $\CC(X)_{w = 0}$ such that $j^* M \cong M_U$.
According to (b), $M$ is isomorphic to a direct sum
$j_{!*} M_U \oplus i_* N$,
for an object $N$ of $\CC(Z)_{w = 0}$.
It is then easy to see that none of properties~(2a)--(4b) holds when
$j_{!*} M_U$ is replaced by $M$, unless $N = 0$.
\end{Proof}

\begin{Ex} \label{2N}
Let $F$ be field contained
in $\BR$, and $X$ a scheme over $\BC$. Put
$\CC(X) := D^b ( \MHM_F X )$ as in Example~\ref{2D}~(b).
As explained there, the additional presence of the $t$-structure
ensures that $D^b ( \MHM_F X )_{w = 0}$ 
is a semi-primary category. Since $D^b ( \MHM_F X )$
is pseudo-Abelian, the intermediate extension
\[
j_{!*} : D^b ( \MHM_F U )_{w = 0} \longto D^b ( \MHM_F X )_{w = 0}^u
\]
can be defined for any
open sub-scheme $U \subset X$. 
Let us clarify its relation with the 
intermediate extension defined as the image of the canonical
transformation from $H^0 j_!$ to $H^0 j_*$ from \cite[(4.5.10)]{Sa}, 
which we shall denote by 
\[
j_{!*}^{t=0} : \MHM_F U \longto \MHM_F X \; . 
\]
Let $M_U \in \MHM_F U$ be pure of weight $n \in \BZ$.
Then $M_U[-n]$ belongs to $D^b ( \MHM_F U )_{w = 0}$.
According to \cite[(4.5.2)]{Sa}, 
\[
(j_{!*}^{t=0} M_U)[-n] \in D^b ( \MHM_F X )_{w = 0} \; .
\]
Now the composition of the adjunction morphisms
\[
i_* i^! (j_{!*}^{t=0} M_U) \longto j_{!*}^{t=0} M_U 
\longto i_* i^* (j_{!*}^{t=0} M_U)
\]
belongs to the radical (since the source is in $D^b ( \MHM_F X )^{t \ge 1}$, 
and the
target in $D^b ( \MHM_F X )^{t \le -1}$). According to Summary~\ref{2M}~(a)~(3),
\[
(j_{!*}^{t=0} M_U)[-n] \cong j_{!*} (M_U[-n]) \; .
\]
In fact, the isomorphism is unique once we require it to give the identity
on $U$: recall that $D^b ( \MHM_F X )_{w = 0}^{t = n}$ is Abelian semi-simple,
therefore, the radical 
$\rad_{D^b ( \MHM_F X )} (j_{!*} (M_U[-n],j_{!*} (M_U[-n])$
is trivial.
We get
\[
[-n] \circ j_{!*}^{t=0} \circ [n] = j_{!*} \quad \text{on} \quad
D^b ( \MHM_F X )_{w = 0}^{t = n} \; ,
\]
canonically for any integer $n$. 
\end{Ex}

\begin{Rem} \label{2O}
Example~\ref{2N} illustrates that the intermediate
extension behaves best on objects which are concentrated simultaneously 
in a single $t$-degree and in a single $w$-degree.
The absence of control of the $t$-degree leads to 
choices in $\CC(X)_{w = 0}$, which in general become unique only  
in the quotient category $\uCX$. 
The absence of control of the $w$-degree leads to non-exactness
of $j_{!*}^{t=0}$. Historically, the second phenomenon seems to be well
tolerated (actually, the author knows of no successful application of 
$j_{!*}^{t=0}$ to objects which are not pure...). 
We think of the first as its weight-structural counterpart. 
\end{Rem}

\begin{Rem} \label{2P}
(a)~The intermediate extension behaves well under passage 
to opposite categories. 
\forget{
This is true first for its underlying notions.
Examples: an $F$-linear category $\FB$ is semi-primary if and only if
the opposite category $\FB^{opp}$ is; the triangulated catego\-ries 
$\CC(U)$, $\CC(X)$ and $\CC(Z)$ are related by gluing if and only if
the opposite categories $\CC(U)^{opp}$, $\CC(X)^{opp}$ and $\CC(Z)^{opp}$
are --- here of course, we have to make the obvious modifications on the level
of the six exact functors used to glue (\emph{e.g.}, the functor 
$j_*^{opp}: \CC(U)^{opp} \to \CC(X)^{opp}$ is opposite to $j_!$). 
Similarly, these ca\-te\-gories carry weight structures $w$ compatible with the
gluing if and only if the opposite categories do --- here, in order to
obtain the weight structure $w^{opp}$ opposite to $w$, we have to
invert the sign of the weight (\emph{e.g.}, the category 
$\CC(X)_{w^{opp} \le 0}^{opp}$ is opposite to $\CC(X)_{w \ge 0}$). 
In particular,
\[
\CC(X)^{opp}_{w^{opp} = 0} 
= \bigl( \CC(X)_{w = 0} \bigr)^{opp} \; ,
\]
and similarly for $U$ and $Z$ instead of $X$.
Thus, the functor 
\[
j_{!*}^{opp} : \bigl( \CC(U)_{w = 0} \bigr)^{opp} 
\longinto \bigl( \uCX \bigr)^{opp} 
\] 
is defined as in Definition~\ref{2J}.
One notes then that $j_{!*}^{opp}$ is in fact 
the functor opposite to $\ujast$.
} \\[0.1cm]
(b)~The intermediate extension is compatible with
duality. 
\forget{
More precisely, assume that equivalences
\[
\BD_\argdot : \CC(\argdot) \isoto \CC(\argdot)^{opp} \; ,
\argdot \in \{ U,X, Z \}
\]
are given, that they are compatible with the gluing 
(\emph{e.g.}, $\BD_X \circ j_* \cong j_*^{opp} \circ \BD_U$)
and with the weight structures (\emph{e.g.}, 
$\BD_X \CC(X)_{w \ge 0} = \CC(X)^{opp}_{w^{opp} \ge 0})$. 
It then follows formally from (a) that 
\[
\BD_X \circ j_{!*} \cong j_{!*}^{opp} \circ \BD_U \; .
\]}
\end{Rem}

One can study compatibility of intermediate extensions with other types 
of operations, for example direct and inverse images (see Section~\ref{6}).
In that context, the following notion turns out to be essential.

\begin{Def}[{\cite[D\'ef.~1.4.6]{AK}}] \label{2Q}
An $F$-linear functor $r: \FA_1 \to \FA_2$ between $F$-linear categories
is called \emph{radicial} if it maps the radical of $\FA_1$ to the radical
of $\FA_2$:
\[
r(\rad_{\FA_1}) \subset \rad_{\FA_2} \; .
\]
\end{Def}

\begin{Prop} \label{6E}
Assume that two triples of $F$-linear pseudo-Abelian
triangulated catego\-ries, related by gluing, and equipped 
with weight structures $w$ compatible with the gluing
are given: $\CC_1(U)$, $\CC_1(X)$, $\CC_1(Z)$, 
and $\CC_2(U)$, $\CC_2(X)$, $\CC_2(Z)$ (we use the same symbols
$j^*$, $i_*$ \emph{etc.} for the two sets of gluing functors).
Assume that furthermore, $F$-linear exact functors 
\[
r_U: \CC_1(U) \longto \CC_2(U) \; , \;
r_X: \CC_1(X) \longto \CC_2(X) \; , \;
r_Z: \CC_1(Z) \longto \CC_2(Z) 
\]
are given, which commute with the gluing functors, and which 
are \emph{$w$-exact}: 
\[
r_\bullet \bigl( \CC_1(\argdot)_{w \le 0} \bigr) 
\subset \CC_2(\argdot)_{w \le 0} \quad \text{and} \quad
r_\bullet \bigl( \CC_1(\argdot)_{w \ge 0} \bigr) 
\subset \CC_2(\argdot)_{w \ge 0} \; . 
\]
In particular, they induce functors 
$r_\bullet: \CC_1(\argdot)_{w = 0} \to \CC_2(\argdot)_{w = 0}$.
Assume that $\CC_1(Z)_{w = 0}$ and $\CC_2(Z)_{w = 0}$ 
are semi-primary. \\[0.1cm]
(a)~If $r_Z: \CC_1(Z)_{w = 0} \to \CC_2(Z)_{w = 0}$ is radicial,  
then $r_X$ maps any object
$M$ of $\CC_1(X)_{w = 0}$ without non-zero direct factors 
in $i_* \CC_1(Z)_{w = 0}$ to 
an object $r_XM$ of $\CC_2(X)_{w = 0}$ 
without non-zero direct factors in $i_* \CC_2(Z)_{w = 0}$. 
The functor $r_X$ descends to induce a functor 
\[
r_X^u : \CC_1(X)_{w = 0}^u \longto \CC_2(X)_{w = 0}^u \; ,
\]
and $r_X^u$ respects the decompositions
\[
\CC_m(X)_{w = 0}^u \isoto \CC_m(U)_{w = 0} \times \CC_m(Z)_{w = 0}^u \; ,
\]
$m = 1,2$, from Theorem~\ref{2K}~(a). In particular, 
the diagram
\[
\vcenter{\xymatrix@R-10pt{
    \CC_1(U)_{w = 0} \ar[d]_{r_U} \ar@{^{ (}->}[r]^-{\ujast} &
    \CC_1(X)_{w = 0}^u \ar[d]_{r_X^u} \\
    \CC_2(U)_{w = 0} \ar@{^{ (}->}[r]^-{\ujast} &
    \CC_2(X)_{w = 0}^u
\\}}
\]
commutes. \\[0.1cm]
(b)~If both $r_Z: \CC_1(Z)_{w = 0} \to \CC_2(Z)_{w = 0}$ and 
$r_U: \CC_1(U)_{w = 0} \to \CC_2(U)_{w = 0}$ are radicial, then
so is $r_X: \CC_1(X)_{w = 0} \to \CC_2(X)_{w = 0}$. 
The functor $r_X^u$ descends further
to induce a functor 
\[
\overline{r_X} : \overline{\CC_1(X)_{w = 0}} 
\longto \overline{\CC_2(X)_{w = 0}} \; ,
\]
and $\overline{r_X}$ respects the decompositions
\[
\overline{\CC_m(X)_{w = 0}} \isoto 
\overline{\CC_m(U)_{w = 0}} \times \overline{\CC_m(Z)_{w = 0}} \; ,
\]
$m = 1,2$, from Theorem~\ref{2K}~(a). 
\end{Prop}

\begin{Proof}
(a): Since $M$ is supposed to satisfy condition~\ref{2M}~(a)~(1),
it satisfies the equivalent condition~\ref{2M}~(a)~(3b).
The latter concerns a morphism belonging to the radical of $\CC_1(Z)_{w = 0}$.
Thanks to our assumption on the restriction of $r_Z$ to $\CC_1(Z)_{w = 0}$, 
condition~\ref{2M}~(a)~(3b)
is then equally satisfied by $r_XM$. Therefore, condition~\ref{2M}~(a)~(1)
is also fulfilled, which proves the claim concerning $r_XM$. 
The claims concerning $r_X^u$ then follow from
Definition~\ref{1Ea}.

(b): Theorem~\ref{2K}~(a) allows us to assume that the radicals 
of $\CC_2(\argdot)_{w = 0}$ are all trivial, and that we are in a trivial 
gluing situation:
\[
\CC_2(X)_{w = 0} = \CC_2(U)_{w = 0} \times \CC_2(Z)_{w = 0} \; .
\]
The assumption is that $r_Z$ and $r_U$ are trivial on the
radicals of $\CC_1(Z)_{w = 0}$ and $\CC_1(U)_{w = 0}$,
respectively. Let $M$ and $N$ be objects of $\CC_1(X)_{w = 0}$, and
\[
f : M \longto N
\]
a morphism belonging to the radical. We need to show that $r_Xf = 0$.
According to Theorem~\ref{2K}~(a), we may assume that $M$ belongs to
the image of $j_{!*}$ (meaning that it does not admit non-zero
direct factors in the image of $i_*$) or to the image of $i_*$,
and likewise for $N$. 
If both $M$ and $N$ are in the image of $j_{!*}$, then according to (a),
so are $r_XM$ and $r_XN$. Furthermore, according to our assumption,
$j^*r_Xf = 0$. But then $r_Xf = 0$. One argues similarly if $M$ and $N$
are both in the image of $i_*$. If one of them is on the image of $j_{!*}$
and the other in the image of $i_*$, then there are no non-trivial
morphisms in $\CC_2(X)_{w = 0}$ between their images under $r_X$. 
\end{Proof}


\bigskip
%
%

\section{The motivic picture}
\label{3}



From now on, $F$ is assumed to be a finite direct product of fields
of characteristic zero.
In addition, let us fix a base scheme
$\BB$, which is of finite type over an excellent scheme 
of dimension at most two. The conventions on schemes and morphisms
are those fixed in our Introduction.
We use the triangulated, $\BQ$-linear categories
$\DBcXM$ of \emph{constructible Beilinson motives} over $X$ 
\cite[Def.~15.1.1]{CD},  
indexed by schemes $X$ (always in the sense of the above conventions). 
In order to have an $F$-linear theory at one's disposal, 
one re-does the construction, but using $F$
instead of $\BQ$ as coefficients \cite[Sect.~15.2.5]{CD}. This yields 
triangulated, $F$-linear categories
$\DBcFXM$ satisfying the $F$-linear analogues of the properties
of $\DBcXM$. In particular, these categories 
are pseudo-Abelian (see \cite[Sect.~2.10]{H}).
Furthermore, the canonical functor $\DBcXM \otimes_\BQ F \to \DBcFXM$
is fully faithful \cite[Sect.~14.2.20]{CD}.
As in \cite{CD}, the symbol $\one_X$
is used to denote the unit for the tensor product in $\DBcFXM$.
We shall employ the full formalism of six operations developed in
\loccit . The reader may choose to consult \cite[Sect.~2]{H} or
\cite[Sect.~1]{W5} for concise presentations of this formalism. \\

Note that since Beilinson
motives satisfy the localization property \cite[Def.~2.3.2]{CD}, 
nilregular schemes are to us ``as good as regular schemes''. Indeed,
the pull-back
from a scheme $X$ to its underlying reduced scheme $X_{red}$ induces
an equivalence $\DBcFXM \isoto \DBcM(X_{red})_F$ \cite[Prop.~2.3.6]{CD}. \\

Beilinson motives can be endowed with weight structures,
thanks to the main results from \cite{H}
(see \cite[Prop.~6.5.3]{Bo} for the case
$X = \Spec k$, for a perfect field $k$). More precisely, the following holds.

\begin{Thm}[{\cite[Thm.~3.3~(ii), Thm.~3.8~(i)--(ii), Rem.~3.4]{H}}] \label{3A}
There are 
ca\-no\-ni\-cal weight structures $w$ on the categories 
$\DBcM(\argdot)$. They have the following properties. 
\begin{enumerate}
\item[(1)] The objects $\one_X(p)[2p]$ belong
to the heart $\DBcM(X)_{w = 0}$, for all integers $p$, whenever $X$ is nilregular.
\item[(2)] For a morphism of schemes $f$, left adjoint functors 
$f^*$, $f_!$ respect the full sub-categories
$\DBcM(\argdot)_{w \le 0}$, and right adjoint
functors $f_*$, $f^!$ and $f^*$  (the latter for smooth $f$) respect
$\DBcM(\argdot)_{w \ge 0}$.
\item[(3)] For a fixed scheme $X$, the heart of $w$
is the pseudo-Abelian completion of the additive envelope of the category
of motives over $X$ of the form
\[
f_* \one_S (p)[2p] \; ,
\]
for proper morphisms $f: S \to X$ with nilregular source $S$, and integers $p$. 
The same statement holds, with ``proper'' replaced by ``projective''.
\end{enumerate}
Furthermore, properties~(1) and (2) 
together characterize the weight structures $w$ uniquely.
\end{Thm}

The proof of \cite[Thm.~3.3, Thm.~3.8]{H} can be imitated
to show that the analogue of Theorem~\ref{3A}
holds for the $F$-linear versions $\DBcM(\argdot)_F$ of the categories
$\DBcM(\argdot)$.
Let us refer to the weight structure $w$ on $\DBcM(\argdot)_F$ as the
\emph{motivic weight structure}. 

\begin{Def}[{cmp.\ \cite[Def.~1.5]{W5}}] \label{3B}
The category 
$\CHFXM$ of \emph{Chow motives} over $X$
is defined as the heart $\DBcM(X)_{F,w = 0}$ of the motivic weight structure
on $\DBcFXM$.
\end{Def}

Thus, the category $\CHFXM$ equals the pseudo-Abelian completion
of $\CHXM \otimes_\BQ F$. \\

If $i: Z \into X$
and $j: U \into X$ are complementary closed, resp.\ open immersions
of schemes, then the category $\DBcFXM$ is obtained by gluing 
$\DBcFUM$ and $\DBcFZM$ \emph{via} the functors 
\[
\vcenter{\xymatrix@R-10pt{
        \DBcFUM \ar[r]^-{j_!} & \DBcFXM \ar[r]^-{i^*} & \DBcFZM \\
        \DBcFUM & \DBcFXM \ar[l]_-{j^*} & \DBcFZM \ar[l]_-{i_*} \\ 
        \DBcFUM \ar[r]^-{j_*} & \DBcFXM \ar[r]^-{i^!} & \DBcFZM
\\}}
\] 
\cite[Prop.~2.3.3~(2), (3), Thm.~2.2.14~(3), Sect.~2.3.1]{CD}. 
We are therefore in the situation studied in Section~\ref{2}. \\

Even though the results of this paper are unconditional, it may be
useful to indicate that our global vision is guided by the following.

\begin{Conj} \label{3C}
Let $X$ be a scheme, and $M$ an indecomposable object
of $\CHFXM$. Then the radical
\[
\rad_{\CHFXM}(M,M) \subset \End_{\CHFXM}(M)
\]
is trivial.
\end{Conj}

\begin{Conj} \label{3D}
Let $X$ be a scheme.
Then the category $\CHFXM$ is semi-primary. 
\end{Conj}

When $X$ is the spectrum of a field $k$, then Conjecture~\ref{3D}
is a consequence of a conjecture proposed independently by Kimura and O'Sullivan
(\emph{cf.} \cite[Sect.~12.1.2]{A}). It predicts that all objects
of $CHM (\Spec k)_F$ are
\emph{finite dimensional}. This notion will turn out
to be very important for us; its precise definition will 
be recalled later (Definition~\ref{5C}). 
For the moment, let us note that according to 
\cite[Thm.~9.2.2]{AK} (see also \cite[Lemma~4.1]{O'S}), 
any pseudo-Abelian
$F$-linear symmetric rigid tensor category is 
(Wedderburn, hence) semi-primary, as soon as all of its objects are 
finite dimensional. \\

Conjectures~\ref{3C} and \ref{3D} 
would be immediate consequences of the existence of a
$t$-structure on $\DBcFXM$ that is transversal to the motivic
weight structure (see Example~\ref{2D}~(a)). 
Conjecture~\ref{3D} guarantees the existence of the intermediate extension
\[
j_{!*} : \CHFUM \longto CHM(X)_F^u \; .
\]
If $M_U$ is an indecomposable object of $\CHFUM$,
then according to Summary~\ref{2M}~(c), $j_{!*} M_U$
is indecomposable in $\CHFXM$. Conjecture~\ref{3C} then predicts
the triviality of its radical, hence the projection
\[
\End_{\CHFXM}(j_{!*} M_U) \longonto \End_{CHM(X)_F^u}(j_{!*} M_U)
\]
would be an isomorphism. In particular, any endomorphism $h_U$
of $M_U$ would then admit a unique extension
to an endomorphism of $j_{!*} M_U$ in $\CHFXM$.

\begin{Ex} \label{3E}
Assume $U$ to be nilregular, and dense in $X$. According to Theorem~\ref{3A}~(1),
the object $\one_U$ then belongs to $\CHFUM$.
On the one hand, Conjecture~\ref{3D} would allow to define 
the value of intermediate extension on $\one_U$,
\[
j_{!*} \one_U \in CHM(X)_F^u \; .
\]
On the other hand, the \emph{motivic intersection
complex}, denoted by the same symbol $j_{!*} \one_U$,
is defined as an extension of $\one_U$ to a Chow motive over $X$,
satisfying 
\[
j^* : \End_{\CHFXM}(j_{!*} \one_U) \isoto \End_{\CHFUM}(\one_U) 
\]
\cite[Def.~2.1]{W5}. If it exists, then Conjecture~\ref{3C}
holds for $M = j_{!*} \one_U$, since $\End_{\CHFXM}(j_{!*} \one_U)$
is then equal to a finite product of copies of $F$.
By the above discussion, the two notions are compatible: 
the motivic intersection complex
would thus equal the value of the functor $j_{!*}$ on $\one_U$. 
\end{Ex}

The main issue to be treated in Sections~\ref{4} and \ref{5} is the 
construction of full sub-categories of Chow motives that 
can actually be shown to satisfy Conjecture~\ref{3D}, \emph{i.e.},
to be semi-primary.
 

\bigskip
%
%

\section{Gluing of motives}
\label{4}



The construction we aim at relies on the process of gluing
certain ``nice'' sub-categories of motives over the members of a stratification.
The following abstract criterion tells us when it is actually possible to glue
given sub-categories.

\begin{Prop} \label{4A}
Let the
three $F$-linear triangulated 
ca\-te\-gories $\CC(U)$, $\CC(X)$ and $\CC(Z)$ be related by gluing:
\[
\vcenter{\xymatrix@R-10pt{
        \CC(U) \ar[r]^-{j_!} & \CC(X) \ar[r]^-{i^*} & \CC(Z) \\
        \CC(U) & \CC(X) \ar[l]_-{j^*} & \CC(Z) \ar[l]_-{i_*} \\ 
        \CC(U) \ar[r]^-{j_*} & \CC(X) \ar[r]^-{i^!} & \CC(Z)
\\}}
\] 
Let $\CD(U) \subset \CC(U)$ and $\CD(Z) \subset \CC(Z)$ be
full, triangulated sub-catego\-ries. \\[0.1cm]
(a)~$\CD(U)$ and $\CD(Z)$ can be glued 
\emph{via} $(j_!,j^*,j_*,i^*,i_*,i^!)$ to give a full,
triangulated sub-category $\CD(X)$ of $\CC(X)$ if and only if
for all objects $M_U$ of $\CD(U)$, both $i^*j_* M_U$ and $i^! j_! M_U$
belong to $\CD(Z)$ (equivalently: $i^*j_* M_U$ \emph{or} $i^! j_! M_U$
belongs to $\CD(Z)$). \\[0.1cm]
(b)~Assume that the condition from (a) is satisfied. Let $M \in \CC(X)$,
and assume that $j^* M \in \CD(U)$. Then the following conditions are
equivalent.
\begin{enumerate}
\item[(1)] $M \in \CD(X)$.
\item[(2)] $i^* M \in \CD(Z)$.
\item[(3)] $i^! M \in \CD(Z)$.
\end{enumerate}
\end{Prop}
 
\begin{Proof}
Note that by \cite[Formula~(1.4.6.4)]{BBD}, the functors $i^*j_*$
and $i^! j_![1]$ are isomorphic. The category $\CD(Z)$ is triangulated,
hence the condition $i^*j_* M_U \in \CD(Z)$ is equivalent to
$i^! j_! M_U \in \CD(Z)$. Given this observation, the proof of part~(a)
is straightforward, as is the description of $\CD(X)$:
an object $M$ of $\CC(X)$ belongs to $\CD(X)$ if and only if
$j^* M \in \CD(U)$, $i^* M \in \CD(Z)$ and $i^! M \in \CD(Z)$. 

In order to prove part~(b), it thus suffices to establish that claim~(1)
is implied by both claim~(2) and claim~(3). In order to show
``$(2) \Rightarrow (1)$'', let $M \in \CC(X)$, and
assume that $j^* M \in \CD(U)$ and $i^* M \in \CD(Z)$. We have to
show that $i^! M \in \CD(Z)$. In order to do so, consider the 
localization triangle
\[
j_! j^* M \longto M \longto i_* i^* M \longto j_! j^* M[1]
\] 
\cite[Axiom~(1.4.3.4)]{BBD}, and apply $i^!$, to get
\[
i^! j_! j^* M \longto i^! M \longto i^* M \longto i^! j_! j^* M[1] \; .
\] 
Since $j^* M \in \CD(U)$ and the condition from (a) is satisfied,
the term $i^! j_! j^* M$ belongs to $\CD(Z)$. The same is true for $i^* M$.
Therefore, $i^! M \in \CD(Z)$. The proof of ``$(3) \Rightarrow (1)$''
is dual.
\end{Proof}

\begin{Rem} \label{4B}
Assume that $\CD(U) \subset \CC(U)$ and $\CD(Z) \subset \CC(Z)$ can be glued 
to give $\CD(X) \subset \CC(X)$. \\[0.1cm]
(a)~It follows from the localization triangle
that $\CD(X)$ can be described as the full, triangulated sub-category of $\CC(X)$
generated by $j_! \CD(U)$ and $i_* \CD(Z)$. Dually, it equals the
full, triangulated sub-category of $\CC(X)$
generated by $j_* \CD(U)$ and $i_* \CD(Z)$. \\[0.1cm]
(b)~The sub-category $\CD(X)$ is dense in $\CC(X)$ if and only
$\CD(U)$ and $\CD(Z)$ are dense in $\CC(U)$ and $\CC(Z)$, respectively.
\end{Rem}

\forget{
Iterated applications
of Proposition~\ref{4A} yield the following.  

\begin{Cor} \label{4C}
Let $S(\FS) = \coprod_{\sigma \in \FS} \Ss$ be a good stratification
of a scheme $S(\FS)$. Assume that a full, triangulated sub-category
$\CD(\Ss)$ of $\DBcFYpM$ is given for any $\sigma \in \FS$. \\[0.1cm]
(a)~The $\CD(\Ss)$ can be glued to give a full,
triangulated sub-category $\CD(S(\FS))$ of $\DBcFSSM$ if and only if
the following condition is satisfied for any $\sigma \in \FS$, and for any
$\psi$ such that $Y_\psi \subset \bSs$: letting $\jp : \Ss \into \bSs$
and $i_\psi : Y_\psi \into \bSs$ denote the immersions, the composition
\[
i_\psi^* j_{\sigma,*} : \DBcFYpM \longto \DBcM(Y_\psi)_F
\]
maps $\CD(\Ss)$ to $\CD(Y_\psi)$. Equivalently, the composition 
$i_\psi^! j_{\sigma,!}$ maps $\CD(\Ss)$ to $\CD(Y_\psi)$, for any $\sigma$
and $\psi$ as above. \\[0.1cm]
(b)~Assume that the condition from (a) is satisfied. Let $M \in \DBcFSSM$.
Then the following conditions are equivalent.
\begin{enumerate}
\item[(1)] $M \in \CD(S(\FS))$.
\item[(2)] $j^* M \in \CD(\Ss)$, for all $\sigma \in \FS$, where $j$
denotes the immersion $\Ss \into S(\FS)$.
\item[(3)] $j^! M \in \CD(\Ss)$, for all $\sigma \in \FS$. 
\end{enumerate}
(c)~Assume that the condition from (a) is satisfied. Then $\CD(S(\FS))$
is the full, triangulated sub-category of $\DBcFSSM$
generated by the $j_! \CD(\Ss)$, for $\sigma \in \FS$. It also equals the
full, triangulated sub-category of $\DBcFSSM$
generated by the $j_* \CD(\Ss)$, for $\sigma \in \FS$. 
\end{Cor}
}

Let us now turn to motives. 

\begin{Def}[{cmp.\ \cite[Def.~3.1]{L}, \cite[Sect.~3.3]{EL}}] \label{4D}
Let $X$ be a sche\-me. Define the $F$-linear \emph{category of Tate motives over $X$}
as the strict, full, $F$-linear
triangulated sub-category $\DFTXM$ of $\DBcFXM$ generated by
the $\one_{X^0}(p)$, for $p \in \BZ$, and $X^0$ running through the connected components
of $X$. 
\end{Def}

By definition, a strict sub-category is closed under
isomorphisms in the ambient category. Recall that
$F$ is a finite product of fields $F_m$ of characteristic zero; thus,  
\[
\DFTXM = \prod_m DMT(X)_{F_m} \; .
\]
In the sequel, Tate motives will be used to construct certain sub-categories
of the category $\DBcM(\argdot)_F$ of Beilinson motives. It will be important
to know that for appropriate choices, these sub-categories inherit a weight structure
from the motivic weight structure on $\DBcM(\argdot)_F$. \\

Let us illustrate the point
for $\DFTXM$ itself, assuming that $X$ is nilregular. 
The triangulated category
$\DFTXM$ is ge\-ne\-rated by the motives isomorphic
to direct factors of $\one_{X^0}(p)[2p]$, for $p \in \BZ$,
and $X^0$ running through the connected components of $X$. These motives are
Chow motives since the scheme $X$ is nilregular (Theorem~\ref{3A}~(1)).
Define $\CK$ as the strict, full, $F$-linear sub-category of $\DFTXM$ 
of objects, which are finite direct sums of objects isomorphic to
$\one_{X^0}(p)[2p]$, for $p \in \BZ$, and $X^0$ a connected component of $X$.
Thus, $\CK$ generates the triangulated category
$\DFTXM$, and all objects of $\CK$
are Chow motives. In particular, by orthogonality~\ref{2A}~(3)
for the motivic weight structure, 
$\CK$ is \emph{negative}, meaning that
\[
\Hom_{\DFTXM}(M_1,M_2[i]) = \Hom_{\DBcM(X)_F}(M_1,M_2[i]) = 0
\]
for any two objects $M_1, M_2$ of $\CK$, and any integer $i > 0$. \\

Therefore,
\cite[Thm.~4.3.2~II~1]{Bo} can be applied to ensure the existence
of a bounded weight structure $v$ on $\DFTXM \,$, uniquely characterized by the
property of containing $\CK$ in its heart. Furthermore
\cite[Thm.~4.3.2~II~2]{Bo}, the heart of $v$
is equal to the category $\CK'$ of retracts of $\CK$ in $\DFTXM$. 
In particular, it is contained in the heart
of the motivic weight structure. The existence of weight 
filtrations~\ref{2A}~(4) for the weight structure $v$ then formally implies that
\[
DMT (X)_{F,v \le 0} \subset \DBcM(X)_{F,w \le 0} \; ,
\]
and that
\[
DMT (X)_{F,v \ge 0} \subset \DBcM(X)_{F,w \ge 0} \; .
\]
We leave it to the reader to prove from this (cmp.~\cite[Lemma~1.3.8]{Bo})
that in fact
\[
DMT (X)_{F,v \le 0}
= \DFTXM \cap \DBcM(X)_{F,w \le 0} 
\]
and
\[
DMT (X)_{F,v \ge 0} 
= \DFTXM \cap \DBcM(X)_{F,w \ge 0} \; . 
\]
In other words, $v$ is induced by $w$. 

\begin{Rem} \label{4L}
We used the standard pattern 
for the induction of (boun\-ded) weight structures. Its abstract ingredients
are as follows: $\CD \subset \CC$ is a full, dense, triangulated
sub-category of a triangulated ca\-te\-gory $\CC$ underlying a weight structure 
$w$, and as a triangulated category, $\CD$ is generated by a family $\CF$
of objects lying in the heart
$\CC_{w=0}$. Then \cite[Thm.~4.3.2~II]{Bo} 
$w$ induces a bounded weight structure on $\CD$, whose heart
equals the category of retracts of the full, additive sub-category
generated by $\CF$.  
\end{Rem}

Note that even without the assumption on generators of $\CD$,
we could define full sub-categories
\[
\CD_{w \le n} := \CD \cap \CC_{w \le n}
\]
and
\[
\CD_{w \ge n} := \CD \cap \CC_{w \ge n}
\]
of $\CD$. While axioms~\ref{2A}~(1)--(3) are trivially satisfied, the existence of
weight filtrations~\ref{2A}~(4) \emph{within the category $\CD$} is not in general. 
This reflects the fact that in general, the heart $\CD_{w=0}$ may not be big enough,
\emph{i.e.}, may not generate the whole of $\CD$. The reader may find it useful to
consider the case of $\DFTXM$ when $X$ is not nilregular. \\

Let $S$ be a scheme.
By definition, a \emph{good stratification} of $S$ 
indexed by a finite set $\FS$
is a collection of locally closed sub-schemes $\Ss$ indexed by
$\sigma \in \FS$, such that 
\[
S = \coprod_{\sigma \in \FS} \Ss
\]
on the set-theoretic level, and such that the closure $\bSs$
of any stratum $\Ss$ is a union of strata $S_\tau$. In that situation,
we shall often write $S(\FS)$ instead of $S$. 

\begin{Thm} \label{4E}
Let $S(\FS) = \coprod_{\sigma \in \FS} \Ss$ be a good stratification
of a scheme $S(\FS)$. Assume the following for all $\sigma \in \FS$:
($\alpha$)~the stratum $\Ss$ is nilregular, ($\beta$)~for any stratum 
$\iat: \St \into \bSs$ contained in the closure $\bSs$ of $\Ss$, 
the functor
\[
i_\tau^! : \DBcM(\bSs)_F \longto \DBcM(S_\tau)_F
\]
maps $\one_{\bSs}$ to a Tate motive over $S_\tau$. \\[0.1cm]
(a)~The categories of Tate motives $\DFTSsM$, $\sigma \in \FS$, 
can be glued to give a full,
triangulated sub-category $\DFTSSM$ of $\DBcFSSM$. \\[0.1cm]
(b)~Let $M \in \DBcFSSM$.
Then the following conditions are equivalent.
\begin{enumerate}
\item[(1)] $M \in \DFTSSM$.
\item[(2)] $j^* M \in \DFTSsM$, for all $\sigma \in \FS$, where $j$
denotes the immersion $\Ss \into S(\FS)$.
\item[(3)] $j^! M \in \DFTSsM$, for all $\sigma \in \FS$. 
\end{enumerate}
In particular,
the triangulated category $DMT(S(\FS))_F$ of Tate motives over $S(\FS)$
is contained in $\DFTSSM$. \\[0.1cm]
(c)~The category 
$\DFTSSM$ is the strict, full, dense, $F$-linear triangulated sub-category of $\DBcFSSM$
generated by the $\bj_* \one_{\bSs}(p)$, for $\sigma \in \FS$ and 
$p \in \BZ$, where $\bj$ denotes the closed immersion $\bSs \into S(\FS)$. \\
(d)~The category $\DFTSSM$ is pseudo-Abelian.
\end{Thm}

From the equivalence ``$(1) \Leftrightarrow (2)$'' in Theorem~\ref{4E}~(b), we see that
the ca\-te\-gory $\DFTSSM$ is stable under tensor product. Furthermore, if $S(\FS)$ is nilregular,
then duality \cite[Thm.~15.2.4]{CD} implies that $\DFTSSM$ is rigid. 

\begin{Def} \label{4F}
Let $S(\FS) = \coprod_{\sigma \in \FS} \Ss$ be a good stratification
of a scheme $S(\FS)$. Assume the following for all $\sigma \in \FS$:
($\alpha$)~$\Ss$ is nilregular, ($\beta$)~for any 
$\iat: \St \into \bSs$, the functor $i_\tau^!$
maps $\one_{\bSs}$ to a Tate motive over $S_\tau$. \\[0.1cm]
(a)~The category $\DFTSSM$ of Theorem~\ref{4E}
is called the \emph{category of $\FS$-constructible
Tate motives over $S(\FS)$}. \\[0.1cm]
(b)~The full sub-category $\CHTFSSM$ of $\CHFSSM$
of objects lying in both $\DFTSSM$ and $\CHFSSM$ is called the
\emph{category of $\FS$-constructible Chow--Tate motives over $S(\FS)$}. 
\end{Def}

\begin{Rem} \label{4G}
Let $S(\FS) = \coprod_{\sigma \in \FS} \Ss$ be a good stratification
of $S(\FS)$. Assume that the closures $\bSs$ of all strata
$\Ss$, $\sigma \in \FS$ are nilregular.
\emph{Absolute purity} \cite[Thm.~14.4.1]{CD} then implies that 
$i_\tau^! \one_{\bSs} (p) \in DMT(S_\tau)_F$.
In other words, the hypotheses of Theorem~\ref{4E} and Definition~\ref{4F}
are satisfied.
\end{Rem}

\begin{Ex} \label{4Fa}
Let $K$ be a number field, and denote by $S$ the spectrum of its ring
of integers $\Fo_K$. Let $\FS$ be a set consisting
of a finite number of closed points of $S$
(\emph{i.e.}, of maximal ideals of $\Fo_K$), and of the open complement
of their union. Then $S = \coprod_{\sigma \in \FS} \Ss$  
satisfies the hypotheses of Remark~\ref{4G}, hence of Theorem~\ref{4E} and Definition~\ref{4F}.
Thus, the category $\DFToKM$ 
of $\FS$-constructible Tate motives over $\Spec \Fo_K$ is defined,
and so is the category $\CHTFoKM$ of 
$\FS$-constructible Chow--Tate motives over $\Spec \Fo_K$.
By passing to the limit of $\DFToKM$ over all $\FS$, for $F = \BQ$, one recovers the category denoted
${\bf DTM}(\Spec \Fo_K)$ in \cite[Def.~2.2]{Sc}. Theorem~\ref{4E}~(b) and (c) continues to apply
(with $\FS =$ the set of all points of $\Spec \Fo_K$). In particular,
we recover \cite[Thm.~2.4]{Sc}.  
\end{Ex}

\begin{Proofof}{Theorem~\ref{4E}}
Let us start by considering the case of a trivial stratification,
\emph{i.e.}, of a nilregular scheme $S(\FS) = \Ss$.
All claims are tri\-vial, except for the density of $\DFTSsM$ in $\DBcFSsM$ (part~(c)),
and of part~(d). The latter implying the former, we need to show that
$\DFTSsM$ is pseudo--Abelian. 

We may assume that $F$ is a field.
According to hypothesis~($\alpha$), the scheme $\Ss$ is nilregular. Therefore,
the motivic weight structure $w$ induces a weight structure, still denoted $w$,
on $\DFTSsM$. Its heart equals the category $\CK'$ of retracts of 
the category $\CK$ defined as the strict, full, $F$-linear sub-category of $\DFTSsM$ 
of objects, which are finite direct sums of objects isomorphic to
$\one_{S_\sigma^0}(p)[2p]$, for $p \in \BZ$, and $S_\sigma^0$ a connected component of $\Ss$.

According to \cite[Lemma~5.2.1]{Bo}, since the weight structure
on the ca\-te\-gory $\DFTSsM$ is bounded, the latter is pseudo-Abelian if $\CK'$ is.
We are thus reduced to showing that $\CK$ is pseudo-Abelian (hence $\CK' = \CK$).

We may assume that $\Ss$ is connected.
For $p \in \BZ$, write $\CK_p$ for the strict, full, $F$-linear sub-category
of $\CK$ of objects isomorphic to a finite direct sum of copies of $\one_{\Ss}(p)[2p]$.
Thus, any object of $\CK$ is a finite direct sum of objects of $\CK_p$, $p \in \BZ$.
We claim that 
\begin{enumerate}
\item[(i)] for every integer $p$, the category $\CK_p$ is Abelian semi-simple
(hence in particular, pseudo-Abelian),
\item[(ii)] for every pair of integers $p_1, p_2$,
\[
\Hom_\CK (M_1 , M_2) = 0 \; , \; \forall \; M_1 \in \CK_{p_1}, M_2 \in \CK_{p_2} \; ,
\]
as soon as $p_1$ is strictly greater than $p_2$.
\end {enumerate}
Indeed, according to \cite[Cor.~14.2.14]{CD}, 
\[
(\ast) \quad \Hom_\CK (\one_{\Ss}(p_1)[2p_1] , \one_{\Ss}(p_2)[2p_2]) = 
\Gr_\gamma^{p_2-p_1} K_0 \bigl( S_{\sigma,red} \bigr)_\BQ \otimes_\BQ F \; ,
\]
where $S_{\sigma,red}$ denotes the (regular) reduced scheme underlying 
the nilregular scheme $\Ss$, and $\Gr_\gamma^n K_0 ( \bullet )_\BQ$ 
the $n$-th graded object with respect to the \emph{$\gamma$-filtration}
on $K_0 ( \bullet )_\BQ := K_0 ( \bullet ) \otimes_\BZ \BQ$ \cite[Ex.~3.9.1, 
Sect.~3.10]{SGA}. Claim~(1) concerns the case $p_1 = p = p_2$,
and follows from the canonical isomorphism given by \emph{augmentation} 
\[
\Gr_\gamma^0 K_0 \bigl( S_{\sigma,red} \bigr)_\BQ \cong \BQ
\]
\cite[Ex.~3.9.1, Sect.~3.10]{SGA}, and
the fact that its composition with $(\ast)$ is inverse to the structure morphism
of the $F$-algebra
$\End_\CK (\one_{\Ss}(p)[2p])$; in fact, the category $\CK_p$ is equivalent
to the category of finite dimensional $F$-vector spaces. As for claim~(2), observe that the
$\gamma$-filtration on $K_0 ( \bullet )_\BQ$ is concentrated in non-negative degrees \cite[Sect.~3.10]{SGA}.

Now let $M = \oplus_{i=1}^n M_i$ be an object of $\CK$, with $M_i \in \CK_{p_i}$, such that
$p_i > p_j$ whenever $i < j$. Let $e$ be an idempotent endomorphism of $M$.
In order to show that $e$ admits a kernel, we apply induction on the number of components
$n$, the case $n = 1$ resulting from (i). For the induction step, rewrite $M$ as a direct
sum 
\[
M = M_1 \oplus N \; ;
\]  
according to (ii), with respect to this direct sum, 
\[
e =
\left( \begin{array}{cc}
A & B \\
0 & D
\end{array} \right) \; ,
\]
with $A \in \End_\CK (M_1)$, $B \in \Hom_\CK (N,M_1)$, and $D \in \End_\CK (N)$.
The relation $e^2 = e$ is equivalent to the system of relations
\[
(\ast \ast) \quad \quad A^2 = A \; , \; D^2 = D \; , \; AB + BD = B \; .
\] 
By our induction hypothesis, $\ker(A) \subset M_1$ and $\ker(D) \subset N$ exist
(and so do $\ker(\id_{M_1} - A)$ and $\ker(\id_N - D)$).
We leave it as an exercice to the reader to prove, using $(\ast \ast)$, that
the morphism
\[
\left( \begin{array}{cc}
\id_{\ker(A)} & -B \\
0 & \id_{\ker(D)}
\end{array} \right) 
\]
from $\ker(A) \oplus \ker(D)$ to $M$ is a kernel of $e$.

We now turn to the general case.
For parts (a)--(c), we proceed by induction on the number of strata.
We already treated the case where this number equals one. 

As for the induction step, 
let $\Ss$ be an open stratum. Our induction hypothesis implies that 
Theorem~\ref{4E} is true for any locally closed union of strata of 
the complement $S(\FS) - \Ss$, with its induced stra\-ti\-fi\-cation.
Write $j$ for the open immersion $\Ss \into S(\FS)$, and $i$ for the
closed immersion $S(\FS) - \Ss \into S(\FS)$. 

As far as part~(a) is concerned, our
aim is to glue $\DFTSsM$ and $DMT_\FS(S(\FS) - \Ss)_F$.
More precisely, we intend to verify the criterion from Proposition~\ref{4A}~(a).
Given the shape of the latter, we may assume that $S(\FS) = \bSs$. 

Since Theorem~\ref{4E} holds for $\bSs - \Ss$,
the motive $i^* \one_{\bSs} (p) = \one_{\bSs - \Ss}(p)$ belongs
to $DMT_\FS(\bSs - \Ss)_F$, for all $p \in \BZ$. 
We claim that the same is true for the motive
$i^! \one_{\bSs} (p)$.
Indeed, according to criterion~(b)~(3), applied to $\bSs - \Ss$, 
it is sufficient
to show that $i_\tau^! \one_{\bSs} (p)$ is a Tate motive over $S_\tau$,
for the immersion $i_\tau: S_\tau \into \bSs$ 
of any stratum $S_\tau$ of $\bSs - \Ss$. But
this is hypothesis~($\beta$). 

Now $i^* j_* \one_{\Ss} (p)$ is a cone of the canonical morphism 
$i^! \one_{\bSs} (p) \to 
i^* \one_{\bSs} (p)$ 
(see the proof of Proposition~\ref{4A}); 
from what precedes, we see that $i^* j_* \one_{\Ss} (p)$
belongs indeed to $DMT_\FS(\bSs - \Ss)_F$.
By induction hypothesis, the category $DMT_\FS(\bSs - \Ss)_F$ is dense in $\DBcM(\bSs - \Ss)_F$. 
The
object $i^* j_* \one_{\Ss} (p)$ being the direct sum of the $i^* j_* \one_{S_\sigma^0} (p)$,
where $S_\sigma^0$ runs through the connected components of $\Ss$, the functor
$i^* j_*$ maps all $\one_{S_\sigma^0} (p)$ to objects of $DMT_\FS(\bSs - \Ss)_F$.

Since $i^* j_*$ maps the generators of $\DFTSsM$
to $DMT_\FS(\bSs - \Ss)_F$, the whole of $\DFTSsM$ is in fact mapped to
$DMT_\FS(\bSs - \Ss)_F$ under $i^* j_*$. This shows that the criterion
from Proposition~\ref{4A}~(a) is fulfilled, thereby proving 
part~(a) of our claim.

Part~(b) is Proposition~\ref{4A}~(b), applied to our situation.

By Remark~\ref{4B}~(a),
$DMT_\FS(S(\FS))_F$ is contained in the strict, full, 
dense, $F$-linear triangulated sub-category 
of $\DBcM(S(\FS))_F$
generated by $j_! \one_{\Ss}(p)$, for all $p \in \BZ$, 
and by the category $i_* DMT_\FS(S(\FS) - \Ss)_F$. 
But modulo
$i_* DMT_\FS(S(\FS) - \Ss)_F$, the $j_! \one_{\Ss}(p)$
and the $\bj_* \one_{\bSs}(p)$
ge\-ne\-rate the same full, triangulated sub-category, as can be seen 
from the localization triangle
\[
j_! \one_{\Ss}(p) \longto \bj_* \one_{\bSs}(p) \longto i_* \one_{\bSs - \Ss}(p) 
\longto j_! \one_{\Ss}(p)[1] \; .
\]
Part~(c) is therefore implied by part~(d). 

As for the latter, note that all ambiant categories $\DBcM(\bullet)_F$
are pseudo-Abelian. Given (a), and Remark~\ref{4B}~(b), we are thus
reduced to the case of a single stratum, which is already established.
\end{Proofof}

As we shall see (Main Theorem~\ref{5A}), the category $\CHTFSSM$
of $\FS$-constructible Chow--Tate motives
over $S(\FS)$ is (Wedderburn, hence) semi-primary,
thus providing us with a non-trivial example of a sub-category of Chow motives
for which the analogue of Conjecture~\ref{3D} is true. Actually,
we aim at a more general result; to obtain it, 
let us generalize the geometric setting. \\

Fix $S(\FS) = \coprod_{\sigma \in \FS} \Ss$,
add a second good stratification 
$Y(\Phi) = \coprod_{\varphi \in \Phi} \Yp$ 
of a scheme $Y(\Phi)$, and a morphism $\pi: S(\FS) \to Y(\Phi)$.
We assume that the pre-image $\pi^{-1}(\Yp)$ of 
any stratum $\Yp$ of $Y(\Phi)$ is a union of strata $\Ss$.
To abbreviate, let us refer to $\pi$ as a \emph{morphism
of good stratifications}.

\begin{Def} \label{4H}
Let $\pi: S(\FS) \to Y(\Phi)$ be a morphism of good stratifications.
Assume the following for all $\sigma \in \FS$:
($\alpha$)~$\Ss$ is nilregular, ($\beta$)~for any 
$\iat: \St \into \bSs$, the functor $i_\tau^!$
maps $\one_{\bSs}$ to a Tate motive over $S_\tau$. 
Define the category $\pi_* \! \DFTSSM$ 
as the strict, full, $F$-linear triangulated sub-category of $\DBcFYPM$ generated by
the images under $\pi_*$ of the objects of $\DFTSSM$. 
\end{Def}

Recall that according to the conventions fixed in our Introduction, 
the pseudo-Abelian completion
of an additive category $\FA$  
is denoted by $\FA^\natural$. Note that by \cite[Thm.~1.5]{BS}, the category
$\FA^\natural$ is triangulated if $\FA$ is. \\

Using proper base change \cite[Thm.~2.4.50~(4)]{CD}, 
one deduces the following from Theorem~\ref{4E} and Remark~\ref{4B}~(b). 

\begin{Cor} \label{4I}
Let $\pi: S(\FS) \to Y(\Phi)$ be a proper morphism of good stratifications.
Assume the following for all $\sigma \in \FS$:
($\alpha$)~$\Ss$ is nilregular, ($\beta$)~for any 
$\iat: \St \into \bSs$, the functor $i_\tau^!$
maps $\one_{\bSs}$ to a Tate motive over $S_\tau$. \\[0.1cm]
(a)~$\pi_* \! \DFTSSM$ is obtained by gluing the 
$\pi_* DMT_\FS(\pi^{-1}(\Yp))_F$, $\varphi \in \Phi$. \\[0.1cm]
(b)~$\pi_* \! \DFTSnatSM$ is obtained by gluing the 
$\pi_* DMT_\FS(\pi^{-1}(\Yp))_F^\natural$, $\varphi \in \Phi$. \\[0.1cm]
(c)~The category $\pi_* \! \DFTSnatSM$ is the strict, full, dense, $F$-linear triangulated 
sub-catego\-ry of $\DBcFYPM$
generated by the $\pi_{\bar{\sigma},*} \one_{\bSs}(p)$, 
for $\sigma \in \FS$ and 
$p \in \BZ$, where $\pi_{\bar{\sigma}}$ denotes the 
restriction of $\pi$ to $\bSs$.
\end{Cor}

\begin{Cor} \label{4K}
Let $\pi: S(\FS) \to Y(\Phi)$ be a proper morphism of good stratifications.
Assume the following for all $\sigma \in \FS$:
($\alpha$)~$\Ss$ is nilregular, ($\beta$)~for any 
$\iat: \St \into \bSs$, the functor $i_\tau^!$
maps $\one_{\bSs}$ to a Tate motive over $S_\tau$. \\[0.1cm]
(a)~The motivic weight structure $w$ on $\DBcFYPM$ induces a weight structure,
again denoted $w$, on
$\pi_* \! \DFTSSM \subset \DBcFYPM$. \\[0.1cm]
(b)~The weight structure $w$ on $\pi_* \! \DFTSSM$ is bounded. \\[0.1cm]
(c)~The motivic weight structure $w$ on $\DBcFYPM$ induces a weight structure,
still denoted $w$, on
$\pi_* \! \DFTSnatSM \subset \DBcFYPM$. \\[0.1cm]
(d)~The weight structure $w$ on $\pi_* \! \DFTSnatSM$ is bounded. \\[0.1cm]
(e)~The heart of the weight structure on $\pi_* \! \DFTSnatSM$
equals 
\[
\pi_* DMT_\FS(S(\FS))_{F,w=0}^\natural \; , 
\]
the pseudo-Abelian completion of $\pi_* DMT_\FS(S(\FS))_{F,w=0}$.
\end{Cor}

The special case $\pi = \id_{S(\FS)}$ will actually enter the proof of
Corollary~\ref{4K}; let us spell it out.

\begin{Cor} \label{4M}
Assume the following for all $\sigma \in \FS$:
($\alpha$)~$\Ss$ is nilregular, ($\beta$)~for any 
$\iat: \St \into \bSs$, the functor $i_\tau^!$
maps $\one_{\bSs}$ to a Tate motive over $S_\tau$. \\[0.1cm]
(a)~The motivic weight structure on $\DBcFSSM$ induces a weight structure on
$\DFTSSM \subset \DBcFSSM$. \\[0.1cm]
(b)~The weight structure on $\DFTSSM$ is bounded. In other words,
the triangulated category $\DFTSSM$ is generated 
by $\CHTFSSM$, the category of $\FS$-constructible Chow--Tate motives over 
$S(\FS)$. 
\end{Cor}

\begin{Proof}
According to Theorem~\ref{4E}~(a), the category $\DFTSSM$ is obtained by gluing.
The functors $j^* = j^!$, $i^*$, $i^!$ associated to open immersions $j$ and closed 
immersions $i$ respect weights as specified in Theorem~~\ref{3A}~(2);
therefore, the motivic weight structure on $\DFTSSM$ is obtained 
by gluing in the sense of \cite[Sect.~8.2]{Bo}. 
In particular, the formulae from \cite[bottom of page~493]{Bo} hold.

We may thus assume
that the stratification of $S(\FS)$ is trivial, \emph{i.e.}, that $S(\FS) = \Ss$,
and $\DFTSSM = DMT (\Ss)_F$.
In that situation, our claims were established after Definition~\ref{4D}.
\end{Proof}

\medskip

\begin{Proofof}{Corollary~\ref{4K}}
Let $\CK$ be the strict, full, $F$-linear sub-category of $\DBcFYPM$ 
of objects isomorphic to objects in the image  under $\pi_*$ of $\CHTFSSM$.
According to our definition and Corollary~\ref{4M}~(b), 
$\CK$ ge\-ne\-rates the triangulated category $\pi_* \! \DFTSSM$.
Following Remark~\ref{4L}, 
the motivic weight structure induces a weight structure on $\pi_* \! \DFTSSM$,
whose heart contains $\CK$. 
This shows part~(a). 

As we already established, the category $\CK$ generates
the triangulated ca\-te\-gory $\pi_* \! \DFTSSM$; \emph{a fortiori}, the same 
holds for $\pi_* DMT_\FS(S(\FS))_{F,w = 0} \,$. This shows part~(b).

Repeat the same argument with  
$\pi_* DMT_\FS(S(\FS))_{F,w = 0}^\natural$ instead of $\CK$. We get a bounded weight structure,
induced by the motivic weight structure,
on the sub-category $\CD$ of $\DBcFYPM$ generated by 
$\pi_* DMT_\FS(S(\FS))_{F,w = 0}^\natural \,$.
By Remark~\ref{4L},
its heart equals
$\pi_* DMT_\FS(S(\FS))_{F,w = 0}^\natural$ (which contains all of its retracts).
Parts~(c)--(e) are implied by \cite[Prop.~5.2.2]{Bo}, which states that 
$\CD = \pi_* \! \DFTSnatSM \,$.
\end{Proofof}

\forget{
\begin{Rem} \label{4N}
In practice, Corollaries~\ref{4I} and \ref{4K}
will be applied in slightly greater generality: fix a good stratification
$Y(\Phi) = \coprod_{\varphi \in \Phi} \Yp$, but let $S(\FS)$ vary.
More precisely, let $\pi_\alpha: S(\FS_\alpha) \to Y(\Phi)$, $\alpha \in I$ 
be a family of morphisms of good stratifications.
Assume that the closures $\bSs$ of all strata
$\Ss$ of $S(\FS_\alpha)$, $\alpha \in I$ are nilregular, 
and that all $\pi_\alpha$ are proper. Write 
\[
S(\FS) := \coprod_{\alpha \in I} S(\FS_\alpha) \quad \text{and} \quad 
\pi := \coprod_{\alpha \in I} \pi_\alpha : S(\FS) \longto Y(\Phi) \; .
\]
The category $\pi_* \! \DFTSSM$ is then defined  
as the full, triangulated sub-category of $\DBcFYPM$ generated by
the individual $\pi_{\alpha,*} \! \DFTSnSM$, $\alpha \in I$. 
By passing to the direct limit
over finite sub-families of $I$, one sees that
the statements of Corollaries~\ref{4I} and \ref{4K} continue to hold.
\end{Rem}
}


\bigskip
%
%

\section{Construction of semi-primary categories of Chow motives}
\label{5}



This section contains our main structural result. 
We first need to recall some termino\-lo\-gy.

\begin{Def}[{cmp.\ \cite[Def.~5.16]{L2}}] \label{5B}
Let $X$ be a nilregular scheme. Define the category $\CHFXMs$ 
of \emph{smooth Chow motives} over $X$ as the strict, full, dense 
sub-category of $\CHFXM$ generated by the motives of the form
\[
f_* \one_S (p)[2p] \; ,
\]
for proper, smooth morphisms $f: S \to X$, and integers $p$.
\end{Def}

It follows from \emph{duality} \cite[Thm.~15.2.4]{CD}
that the pseudo-Abelian
$F$-linear symmetric tensor category $\CHFXMs$ is \emph{rigid}.
The relation of morphisms in $\DBcM(\argdot)_F$ to $K$-theory
\cite[Cor.~14.2.14]{CD}
shows that morphisms in $\CHFXMs$ are
given by the classical relation to Chow groups: for proper, smooth morphisms
$f: S \to X$, supposed to be of pure relative dimension $d_S$,
and $g: T \to X$, we have 
\[
\Hom_{\CHFXMs} \bigl( f_* \one_S (p)[2p] , g_* \one_T (q)[2q] \bigr)
= \ch^{d_S+q-p} (S \times_X T) \otimes_\BZ F \; .
\] 
Note that according to  \cite[Prop.~5.19, Cor.~6.14]{L2}, this identification 
is compatible with composition (of morphisms on the left hand,
of correspondences in the classical sense on the right hand).

\begin{Rem} 
Assume that $j: U \into X$ is the immersion of a dense open subscheme
of $X$. Then
\cite[Prop.~5.1.1]{O'S2} implies that the functor
\[
j^* : \CHFXMs \longto CHM^s(U)_F
\]
induces an isomorphism
\[
CHM^s(X)_F^u \isoto CHM^s(U)_{X,F} \; ,
\]
where $CHM^s(U)_{X,F} \subset CHM^s(U)_F$ denotes the full subcategory
of Chow motives admitting a smooth extension to $X$. This shows that
the theory of intermediate extensions on $CHM^s(U)_{X,F}$ is trivial. 
\end{Rem}

\begin{Def}[{\cite[D\'ef.~9.1.1]{AK}, cmp.\ \cite[Def.~3.7]{K}}] \label{5C}
Let $X$ be a nilregular scheme, and $M \in \CHFXMs$ a smooth Chow motive
over $X$. \\[0.1cm]
(a)~The motive $M$ is \emph{evenly finite dimensional} if there exists 
a positive integer $m$ such that $\Lambda^m M = 0$. \\[0.1cm]
(b)~The motive $M$ is \emph{oddly finite dimensional} if there exists 
a positive integer $m$ such that $\Sym^m M = 0$. \\[0.1cm]
(c)~The motive $M$ is \emph{finite dimensional} if there exists 
a direct sum decomposition $M = M^+ \oplus M^-$, where 
$M^+$ is evenly finite dimensional and $M^-$ is oddly finite
dimensional.
\end{Def}

Let us get back to the situation considered at the end of Section~\ref{4}.
The proof of the following uses everything said so far.

\begin{MainThm} \label{5A}
Let $\pi: S(\FS) \to Y(\Phi)$ be a proper
morphism of good stratifications.
Assume that all strata $\Yp$ and $\Ss$, $\varphi \in \Phi$, $\sigma \in \FS$,
are nilregular, and that for any 
$\iat: \St \into \bSs$, the functor $i_\tau^!$
maps $\one_{\bSs}$ to a Tate motive over $S_\tau$. Furthermore, assume 
that for all $\sigma \in \FS$ such that $\Ss$ is a stratum
of $\pi^{-1}(\Yp)$, $\varphi \in \Phi$,
the following holds: the morphism induced 
by $\pi$ on $\Ss$,
\[
\pis : \Ss \longto \Yp
\]
factors over a scheme $\Bs$,
\[
\pis = \pios \circ \pits : \Ss \stackrel{\pits}{\longto} \Bs 
\stackrel{\pios}{\longto} \Yp \; ,
\]
and the morphisms $\pios$ and $\pits$ satisfy the following conditions.
\begin{enumerate}
\item[$(1)_\sigma$] $\pios$ is proper and smooth, 
and the smooth Chow motive 
\[
\pi'_{\sigma,*} \one_{\Bs} \in \CHFYpMs \subset \CHFYpM
\]
is finite dimensional.
\item[$(2)_\sigma$] The motive 
\[
\pi''_{\sigma,*} \one_{\Ss} \in \DBcFBsM
\]
belongs to the category $\DFTBsM$ of Tate motives over $\Bs$.
\end{enumerate}
Then the heart $\pi_* DMT_\FS(S(\FS))_{F,w=0}^\natural$
of the motivic weight structure $w$ on  
$\pi_* \! \DFTSnatSM \subset \DBcFYPM$ (Corollary~\ref{4K}~(d))
is a Wedderburn category. In particular, it is semi-primary.
\end{MainThm} 

\begin{Proof}
Condition~$(2)_\sigma$ implies that the functor $\pi''_{\sigma,*}$
respects the categories of Tate motives:
\[
\pi''_{\sigma,*} \DFTSsM \subset \DFTBsM \; .
\]
It follows that
\[
\pi_{\sigma,*} \DFTSnatsM \subset \pi'_{\sigma,*} \DFTBnatsM 
\]
for all $\sigma \in \FS$. Thus, the
category $\pi_* DMT_{\FS}(\pi^{-1}(\Yp))_F^\natural$,
for fixed $\varphi \in \Phi$, 
is contained in the strict, full, dense, triangulated sub-category $\CD_\varphi$
of $\DBcFYpM$ generated by the 
\[
\pi'_{\sigma,*} \one_{\Bs} (p) \; ,
\]
for all $p \in \BZ$ and $\sigma \in \FS$,
such that $\Ss$ is a stratum of $\pi^{-1}(\Yp)$. 
According to Remark~\ref{4L},
the motivic weight structure on $\DBcFYpM$ induces
a weight structure on $\CD_\varphi$, whose heart $\CD_{\varphi, w=0}$ 
is the pseudo-Abelian
completion of the strict, full, additive sub-category of $\DBcFYpM$ 
of objects, which are finite direct sums of objects isomorphic to
$\pi'_{\sigma,*} \one_{\Bs}(p)[2p]$, for $p \in \BZ$ and $\sigma \in \FS$ as
above. Condition~$(1)_\sigma$ says that $\pi'_{\sigma,*} \one_{\Bs}$
is a finite dimensional smooth Chow motive. But then so are all 
$\pi'_{\sigma,*} \one_{\Bs}(p)[2p]$, $p \in \BZ$: the canonical
action of the symmetric group on $\one_{\Bs}(p)[2p] \otimes \one_{\Bs}(p)[2p]$
is trivial \cite[first part of Prop.~16.2.10]{CD}. Hence any decomposition
\[
\pi'_{\sigma,*} \one_{\Bs} = M^+ \oplus M^- \; ,
\]
where $M^+$ is evenly finite dimensional, and $M^-$
oddly finite dimensional. will indeed induce a decomposition
\[
\pi'_{\sigma,*} \one_{\Bs}(p)[2p] = M^+(p)[2p] \oplus M^-(p)[2p] \; ,
\]
where $M^+(p)[2p]$ is evenly finite dimensional, and $M^-(p)[2p]$
oddly finite dimensional. By \cite[Prop.~9.1.12]{AK} 
(cmp.\ \cite[Prop.~6.9]{K}), all objects of the category $\CD_{\varphi, w=0}$
are finite dimensional. The same is therefore true for the objects
of the category $\pi_* DMT_{\FS}(\pi^{-1}(\Yp))_{F,w=0}^\natural$. By
\cite[Thm.~9.2.2]{AK}, the latter
is a Wedderburn category, for all $\varphi \in \Phi$.

Now apply Corollary~\ref{4I}~(b) and Theorem~\ref{2K}~(b).
\end{Proof}

For future reference,
let us extract the following from the proof of Main Theorem~\ref{5A}.

\begin{Prop} \label{5AL}
Keep the notation and the hypotheses of Main Theorem~\ref{5A},
and fix $\varphi \in \Phi$. Then every object of
\[
\pi_* DMT_{\FS}(\pi^{-1}(\Yp))_{F,w=0}^\natural
\]
is a direct factor of a finite direct sum of objects isomorphic to
$\pi'_{\sigma,*} \one_{\Bs}(p)[2p]$, for $p \in \BZ$ and $\sigma \in \FS$,
such that $\Ss$ is a stratum of $\pi^{-1}(\Yp)$. 
\end{Prop} 

Main Theorem~\ref{5A} yields categories of Chow motives
to which the theory of intermediate extensions developed in Section~\ref{2}
applies (see in particular Summary~\ref{2M}). For better legibility, let us
fix the hypotheses.

\begin{Ass} \label{As1}
Let $\pi: S(\FS) \to Y(\Phi)$ be a 
proper morphism of good stratifications. 
Let $\Phi_U$ be an open subset of $\Phi$, meaning that it gives rise to
a good stratification $Y(\Phi_U)$, and that $Y(\Phi_U) \into Y(\Phi)$
is an open immersion. Denote
by $\Phi_Z$ its complement. Write $\FS_U := \pi^{-1} \Phi_U$ and
$\FS_Z := \pi^{-1} \Phi_Z$. \\[0.1cm]
(a)~For all $\sigma \in \FS$, the strata
$\Ss$ are nilregular, and for any 
$\iat: \St \into \bSs$, the functor $i_\tau^!$
maps $\one_{\bSs}$ to a Tate motive over $S_\tau$. \\[0.1cm]
(b)~For all $\varphi \in \Phi_Z$, 
the strata $\Yp$ are nilregular. \\[0.1cm]
(c)~For all $\varphi \in \Phi_Z$ and $\sigma \in \FS_Z$
such that $\Ss$ is a stratum
of $\pi^{-1}(\Yp)$, the morphism
$\pis : \Ss \to \Yp$ can be factorized,
\[
\pis = \pios \circ \pits : \Ss \stackrel{\pits}{\longto} \Bs 
\stackrel{\pios}{\longto} \Yp \; ,
\]
such that conditions~\ref{5A}~$(1)_\sigma$ and $(2)_\sigma$ are satisfied. 
\end{Ass}

Note that only Assumption~\ref{As1}~(a) is needed for the results on gluing
(Section~\ref{4}). In particular (Corollary~\ref{4I}~(b)),
$\pi_* \! \DFTSnatSM$ is obtained by gluing
$\pi_* \! \DFTSUnatSM$ and 
$\pi_* \! \DFTSZnatSM \,$.

\begin{Cor} \label{5Aa}
Let $\pi: S(\FS) \to Y(\Phi)$ be a 
proper morphism of good stratifications. 
Let $\Phi_U \subset \Phi$ be open, 
and denote by $\FS_U \subset \FS$ the pre-image of $\Phi_U$ under $\pi$.
Suppose that Assumption~\ref{As1} holds. Then the hypotheses
of Theorem~\ref{2K}~(a) are fulfilled for $\CC(X) = \pi_* \! \DFTSnatSM$, 
$\CC(U) = \pi_* \! \DFTSUnatSM$ and $\CC(Z) = \pi_* \! \DFTSZnatSM \,$. 
In par\-ti\-cu\-lar, the functor
\[
\ujast : \pi_* DMT_{\FS_U}(S(\FS_U))_{F,w = 0}^\natural
\longinto \pi_* DMT_\FS(S(\FS))_{F,w=0}^{\natural,u} 
\]
is defined, and it satisfies the properties listed in Summary~\ref{2M}. 
\end{Cor}

\begin{Proof}
Main Theorem~\ref{5A}, applied to 
$\pi_Z: S(\FS_Z) \to Y(\Phi_Z)$, tells us that 
$\pi_* DMT_{\FS_Z}(S(\FS_Z))_{F,w = 0}^\natural$ is semi-primary.
\end{Proof}

In the applications we have in mind, the morphisms $\pits$
will turn each $\Ss$ into a torsor under a split torus over $\Bs$.
We leave it to the reader to prove, using the \emph{$\BA^1$-homotopy
property}, that in this case, condition~\ref{5A}~$(2)_\sigma$ holds.  
Let us now identify concrete situations where
condition~\ref{5A}~$(1)_\sigma$ is satisfied. 
The following principles turn out to be very useful.

\begin{Prop} [{\cite[pp.~54--55]{O'S2}}] \label{5D}
Let $g: X' \to X$ be a morphism of nilregular schemes, $M \in \CHFXMs$
and $M' \in \CHFXoMs$. \\[0.1cm]
(a)~If $M$ is finite dimensional (resp.\ evenly finitely dimensional,
resp.\ oddly finitely dimensional), then so is $g^* M \in \CHFXoMs$. \\[0.1cm] 
(b)~If $g$ is dominant, and $g^* M$ is finite dimensional (resp.\ 
evenly finitely dimensional,
resp.\ oddly finitely dimensional), then so is $M$. \\[0.1cm]
(c)~If $g$ is finite, \'etale, and $M'$ is finite dimensional
(resp.\ evenly finitely dimensional, resp.\ oddly finitely dimensional), 
then so is $g_* M' \in \CHFXMs$. 
\end{Prop} 

\forget{
\begin{Proof}
$g^*$ is an additive tensor functor \cite[Ex.~1.1.22]{CD},
and therefore, part~(a) of our claim is obvious. 

As for part~(b), let us show the claim first for 
evenly, resp.\ oddly finite dimensional $g^* M$.
It suffices for this to show that $g^*$ is conservative on $\CHFXMs$. 
If $g$ is surjective, then
this follows from \cite[Thm.~14.3.3]{CD},
which says that the categories $\DBcM (\argdot)_F$ are 
\emph{separated} in the sense of \cite[Def.~2.1.7]{CD}.
We are therefore reduced to the case of an immersion $g$
of an open, dense sub-scheme $X'$. In this case, conservativity
of $g^*$ on $\CHFXMs$ follows from \cite[Prop.~5.1.1]{O'S2}. 
\forget{
the kernel of $g^*$ is then \emph{strongly tensor nilpotent}:
there exists a positive integer $N$ such that for any $N$-tuple
of elements $\beta_i \in \Hom_{\CHFXMs} (M_i,L_i)$ with trivial
images in $\Hom_{\CHFXoMs} (g^*M_i,g^*L_i)$, the tensor product
\[
\beta_1 \otimes \ldots \otimes \beta_N 
\in \Hom_{\CHFXMs} \bigl( \otimes_{i=1}^N M_i,\otimes_{i=1}^N L_i)  
\]
is trivial.

In particular, it is nilpotent, \emph{whatever the composition law}.
}
The case of finite dimensional $g^* M$ will be treated in a moment.

In order to prove part~(c), we may assume $M'$ to be evenly (resp.\ oddly)
finitely dimensional. Furthermore, thanks to the cases of part~(b)
we already proved,
and finite, \'etale base change, we may assume that $X'$
is a finite disjoint union of copies of $X$, and $g$ an isomorphism on each
component. In that case, $g_*M'$ is a direct sum of smooth Chow
motives over $X$, each
of which is evenly (resp.\ oddly) finite dimensional.
Our claim thus follows from \cite[Prop.~9.1.4~a)]{AK}.

It remains to prove part~(b), assuming that $g^* M$ is only finite dimensional. 
If $g$ is surjective, then \cite[Prop.~5.1.2]{O'S2} allows to reduce 
to the case where $g$ is finite, \'etale (see the discussion in 
\cite[top of p.~55]{O'S2}). According to (c), 
the motive $g_*g^* M \in \CHFXMs$ is then
finite dimensional.
Now since $g$ is finite and \'etale, the functors $g_*$ and $g^*$
are adjoint to each other in each direction. Furthermore, the composition
\[
\id_X \longto g_*g^* \longto \id_X
\]
of the adjunction morphisms is multiplication by a positive integer (the degree)
on each connected component of $X$ (this can be checked on the level
of $K$-theory, thanks to \cite[Cor.~14.2.14]{CD}). The category
$\CHFXMs$ being $\BQ$-linear, the motive $M$ is therefore a direct factor 
of the finite dimensional motive $g_*g^* M$. 
Thus \cite[Prop.~9.1.12]{AK} 
(cmp.\ \cite[Prop.~6.9]{K}), it is itself finite dimensional.
Thus, we are again reduced to the case of an immersion
$g$ of an open, dense sub-scheme $X'$. Our claim follows then from
the fact that $g^* : \CHFXMs \to \CHFXoMs$ is full, and reflects
split monomorphisms and split epimorphisms  \cite[Prop.~5.1.1]{O'S2}. 
\end{Proof}
}
\begin{Exs} \label{5E}
Let $f: S \to X$ be a proper, smooth morphism between nilregular schemes.
We give a list of cases where $M:= f_* \one_S$ is
finite dimensional. \\[0.1cm]
(a)~The morphism $f$ is finite and \'etale (in which case $M$ 
is an \emph{Artin motive} over $X$). Indeed,
$M$ is then evenly finite dimensional according to 
Proposition~\ref{5D}~(c). \\[0.1cm]
(b)~$S$ is an Abelian scheme over $X$ \emph{via} $f$.
Indeed, as observed in \cite[lower half of p.~61]{O'S2},
Proposition~\ref{5D}~(b) allows to reduce to the case where
$X$ is the spectrum of an algebraically closed field. 
The claim then follows from a
classical result of Shermenev \cite[Theorem]{Sh}.
It generalizes \cite[Thm.~(3.3.1)]{Kue},
which was proved using the Fourier--Mukai transform for Abelian schemes,
and assuming that $X$ is quasi-projective and smooth over a field. 
\\[0.1cm]
(c)~More generally (cmp.~\cite[Claim~11 on p.~48]{Sh}), 
each fibre of $f$ over a geometric generic point of $X$ is  
the target of a proper, surjective morphism, whose source is a product
of smooth, projective curves. Indeed, apply Proposition~\ref{5D}~(b)
and \cite[Prop.~6.9, Rem.~6.6, Cor.~5.11 and Cor.~4.4]{K}.\\[0.1cm]
(d)~Still more generally, and still thanks to Proposition~\ref{5D}~(b),
$M$ is finite dimensional as soon as it is isomorphic to one of the earlier
cases after a dominant base change. For example, $M = f_* \one_S$ is
finite dimensional if $S$ is a torsor under an Abelian scheme over $X$
(for the $fpqc$-topology, say).
\end{Exs}

\forget{
Here is one case where the hypotheses of Main Theorem~\ref{5A}
are obviously satisfied.

\begin{Cor} \label{5Ea}
Assume that the closures $\bSs$ of all strata
$\Ss$, $\sigma \in \FS$, are nilregular. Then the category $\CHTFSSM$ of 
$\FS$-constructible Chow--Tate motives over $S(\FS)$ is semi-primary.
\end{Cor}

\begin{Proof}
$\CHTFSSM$ is the heart of the motivic weight structure $w$ on $\DFTSSM$.
Now apply Main Theorem~\ref{5A} (with $\pi = \id_{S(\FS)}$).
\end{Proof}

\begin{Ex} \label{5Eb}
The conclusion of Corollary~\ref{5Ea} holds in particular
in the situation considered in Example~\ref{4Fa}, where $S = \Spec \Fo_K$
for a number field $K$. Passage to the limit shows that it holds equally
for the category \? of \cite{Sc}. Note that in this particular geometrical
context, an alternative proof of Corollary~\ref{5Ea} can be given.
Indeed, the category \? then carries a bounded
$t$-structure \cite[\?]{Sc}, which is orthogonal to the motivic weight
structure. Now apply Example~\ref{2D}~(a). 
\forget{
\\[0.1cm]
(b)~The formal reasoning employed in Example~\ref{2N} shows that
the theory of intermediate extensions from Section~\ref{2}, when
applied to the category \? of \cite{Sc}, is compatible, in the sense
analogous to Example~ref{2D}, with the intermediate extension
defined on the heart of the $t$-structure \cite[\?]{Sc}.   
}
\end{Ex}
}
The following consequence of Proposition~\ref{5D}~(c) will be needed
in the sequel.

\begin{Cor} \label{5F}
Let $\pi: S(\FS) \to Y(\Phi)$ be a 
morphism of good stra\-ti\-fi\-cations.
Assume that all strata $\Yp$ of $Y(\Phi)$, $\varphi \in \Phi$ are nilregular. 
Let $a: Y(\Phi) \to W(\Phi)$
be a surjective finite morphism of good stratifications, with the same 
index set $\Phi$: $a^{-1} (\Wp) = \Yp$ for all $\varphi \in \Phi$. 
Assume that the morphisms induced by $a$ on $\Yp$,
\[
a_\varphi : \Yp \longto \Wp \; ,  
\]
are all finite and \'etale (hence all $\Wp$, $\varphi \in \Phi$
are nilregular).
Furthermore, assume 
that for all $\sigma \in \FS$
such that $\Ss$ is a stratum
of $\pi^{-1}(\Yp)$, $\varphi \in \Phi$, a factorization of 
$\pis : \Ss \to \Yp$ is given,
\[
\pis = \pios \circ \pits : \Ss \stackrel{\pits}{\longto} \Bs 
\stackrel{\pios}{\longto} \Yp \; .
\]
Then if all $\pios$ satisfy condition~\ref{5A}~$(1)_\sigma$,
then so do all $a \circ \pios$. In particular, 
if the morphism $\pi$ satisfies the hypotheses of 
Main Theorem~\ref{5A}, then so does the morphism 
$a \circ \pi: S(\FS) \to W(\Phi)$.
\end{Cor}


\bigskip
%
%

\section{Compatibility with certain direct and inverse images}
\label{6}



In this section, we aim at compatibility statements of type 
``$a_* \circ j_{!*} = j_{!*} \circ a_*$'' and 
``$b^* \circ j_{!*} = j_{!*} \circ b^*$'', for 
finite morphisms $a$ satisfying certain conditions
(Theorem~\ref{6G}), and certain smooth morphisms $b$ (Theorem~\ref{6H}). 
Actually, Theorem~\ref{6G} should hold for arbitrary finite $a$, as
is suggested by \cite[Cor.~4.1.3]{BBD}:
in the context of triangulated categories of sheaves, $a_*$ is exact
for the (perverse) $t$-structure. 
Similar remarks apply to Theorem~\ref{6H}.
Our starting point is the following result.

\begin{Prop} [{\cite[p.~55]{O'S2}}] \label{6A}
Assume that $F$ is a field of characteristic zero. 
Let $g: X' \to X$ be a dominant morphism of nilregular connected schemes.
Then the inverse image $g^* : \CHFXMs \to \CHFXoMs$
maps the maximal tensor ideal $\CN_X$ of $\CHFXMs$ to the maximal tensor ideal
$\CN_{X'}$ of $\CHFXoMs$:
\[
g^*(\CN_X) \subset \CN_{X'} \; .
\]
\end{Prop} 

A word of explanation is in order: due to the restrictions on $F$, $X$
and $X'$, the categories $\CHFXMs$ and $\CHFXoMs$ are not only $F$-linear
rigid symmetric tensor categories, but they also satisfy the additional
hypothesis from \cite[Sect.~7]{AK}: $F$ is a field, 
\[
F = \End_{\CHFXMs}(\one_X) \quad 
\text{and} \quad F = \End_{\CHFXoMs}(\one_{X'}) 
\]
\cite[Cor.~14.2.14]{CD}.
Therefore, \cite[Prop.~7.1.4~b)]{AK} can be applied, ensuring 
existence and unicity of a \emph{maximal tensor ideal, i.e.},
an ideal $\CN_X$, resp.\ $\CN_{X'}$, 
which is maximal among the tensor ideals unequal to $\CHFXMs$,
resp.\ $\CHFXoMs$. \\

With an appropriate definition of the tensor ideal $\CN_\argdot$,
Proposition~\ref{6A} admits an obvious generalization to the setting
we have considered so far ($F$ a finite product of fields of
characteristic zero, schemes which are not necessarily irreducible).
Denote by $CHM^s(\argdot)_F^{fd}$ the full sub-category
of $CHM^s(\argdot)_F$ formed by finite dimensional objects.  
By \cite[Thm.~9.1.12, Prop.~9.1.4~b)]{AK} 
(cmp.\ \cite[Prop.~6.9, Prop.~5.10]{K}), 
$CHM^s(\argdot)_F^{fd}$ is a dense tensor sub-category of 
$CHM^s(\argdot)_F$. It remains rigid according to \cite[Prop.~9.1.4~c)]{AK}.

\begin{Cor} \label{6C}
Let $g: X' \to X$ be a morphism of nilregular schemes,
whose restriction to any connected component of $X'$ 
is dominant over some connected component of $X$.
Then $g^* : \CHFfXMs \to \CHFfXoMs$ is radicial (see Definition~\ref{2Q}).
\end{Cor}

\begin{Proof}
According to \cite[Thm.~9.2.2]{AK},
the ideal $\CN_\argdot$ of $CHM^s(\argdot)_F^{fd}$
coincides with the radical. 
The claim then follows from Proposition~\ref{6A}. 
\end{Proof}

\begin{Cor} \label{6D}
Let $g: X' \to X$ be a finite, \'etale morphism of nilregular schemes.
Then both the inverse image $g^* : \CHFfXMs \to \CHFfXoMs$
and the direct image $g_*: \CHFfXoMs \to \CHFfXMs$ (Proposition~\ref{5D}~(c))
are radicial.
\end{Cor}

\begin{Proof}
For $g^*$, we have the more general statement of Corollary~\ref{6C}.

In order to prove the claim for $g_*$, note first that it suffices to show it
after a surjective base change $f: Y \to X$ with nilregular source $Y$: indeed,
according to \cite[Thm.~14.3.3]{CD},
the categories $\DBcM (\argdot)_F$ are 
\emph{separated} in the sense of \cite[Def.~2.1.7]{CD}.
Therefore, the inverse image $f^*$ is conservative. 
According to Corollary~\ref{6C}, it is radicial on $CHM^s(\argdot)_F^{fd}$.
By \cite[Lemme~1.4.7]{AK},
\[
f^* : \CHFfXMs \longto CHM^s(Y)_F^{fd}
\]
detects elements in the radical. 

Then, taking $f$ to be an appropriate
finite, \'etale morphism,
we may assume that $X'$
is a finite disjoint union of copies of $X$, and $g$ an isomorphism on each
component. In that case, our claim follows from 
\cite[Cor.~1.4.5]{AK}.
\end{Proof}

Let us now show radiciality of $a_*$, for morphisms $a$ of the type
considered at the end of Section~\ref{5}. 

\begin{Cor} \label{6F}
Let $\pi: S(\FS) \to Y(\Phi)$ be a
proper morphism of good stratifications.
Assume that all strata $\Yp$ and $\Ss$, $\varphi \in \Phi$, $\sigma \in \FS$,
are nilregular, and that for any 
$\iat: \St \into \bSs$, the functor $i_\tau^!$
maps $\one_{\bSs}$ to a Tate motive over $S_\tau$. 
Furthermore, assume 
that for all $\sigma \in \FS$
such that $\Ss$ is a stratum
of $\pi^{-1}(\Yp)$, the morphism
$\pis : \Ss \to \Yp$ can be factorized,
\[
\pis = \pios \circ \pits : \Ss \stackrel{\pits}{\longto} \Bs 
\stackrel{\pios}{\longto} \Yp \; ,
\]
such that conditions~\ref{5A}~$(1)_\sigma$ and $(2)_\sigma$ are satisfied. 
Let $a: Y(\Phi) \to W(\Phi)$
be a surjective finite morphism of good stratifications, with the same 
index set $\Phi$. 
Assume that the morphisms 
$a_\varphi : \Yp \to \Wp$
are all finite and \'etale. Then the direct image
\[
a_* : \pi_* DMT_\FS(S(\FS))_{F,w=0}^\natural \longto 
(a \circ \pi)_* DMT_\FS(S(\FS))_{F,w=0}^\natural
\]
is radicial.
\end{Cor}

By definition, $a_*$ maps 
$\pi_* \! \DFTSnatSM$ to
$(a \circ \pi)_* \! \DFTSnatSM$. 
Since $a$ is finite, hence proper, it is
$w$-exact; therefore, it respects the hearts.

\medskip

\begin{Proofof}{Corollary~\ref{6F}}
According to Corollary~\ref{4I}~(b), 
the categories $\pi_* \! \DFTSnatSM$
and $(a \circ \pi)_* \! \DFTSnatSM$ are obtained by gluing. 
Thanks to Proposition~\ref{6E}~(b), we are reduced to showing that each 
\[
a_{\varphi,*} : \pi_* DMT_{\FS}(\pi^{-1}(\Yp))_{F,w=0}^\natural
\longto 
(a \circ \pi)_* DMT_{\FS}(\pi^{-1}(\Yp))_{F,w=0}^\natural
\]
is radicial. By Proposition~\ref{5AL},
the source of $a_{\varphi,*}$ is contained in $CHM^s(\Yp)_F^{fd}$, 
and the target
in $CHM^s(\Wp)_F^{fd}$. 
Our claim thus follows from our assumption on $a_\varphi$,
and from Corollary~\ref{6D}.
\end{Proofof}

We are ready to state the main results of this section. Again, let us fix the 
hypotheses.

\begin{Ass} \label{As}
Let $\pi: S(\FS) \to Y(\Phi)$ and $a: Y(\Phi) \to W(\Phi)$ be 
proper, resp.\ 
surjective finite morphisms of good stratifications. 
Let $\Phi_U$ be an open subset of $\Phi$, and denote
by $\Phi_Z$ its complement. Write $\FS_U := \pi^{-1} \Phi_U$ and
$\FS_Z := \pi^{-1} \Phi_Z$. \\[0.1cm]
(a)~For all $\sigma \in \FS$, the strata
$\Ss$ are nilregular, and for any 
$\iat: \St \into \bSs$, the functor $i_\tau^!$
maps $\one_{\bSs}$ to a Tate motive over $S_\tau$. \\[0.1cm]
(b)~For all $\varphi \in \Phi_Z$, 
the strata $\Wp$ and $\Yp$ are nilregular. \\[0.1cm]
(c)~For all $\varphi \in \Phi_Z$, the morphism 
$a_\varphi : \Yp \to \Wp$
is finite and \'etale. \\[0.1cm]
(d)~For all $\varphi \in \Phi_Z$ and $\sigma \in \FS_Z$
such that $\Ss$ is a stratum
of $\pi^{-1}(\Yp)$, the morphism
$\pis : \Ss \to \Yp$ can be factorized,
\[
\pis = \pios \circ \pits : \Ss \stackrel{\pits}{\longto} \Bs 
\stackrel{\pios}{\longto} \Yp \; ,
\]
such that conditions~\ref{5A}~$(1)_\sigma$ and $(2)_\sigma$ are satisfied. 
\end{Ass}

\begin{Thm} \label{6G}
Let $\pi: S(\FS) \to Y(\Phi)$ and $a: Y(\Phi) \to W(\Phi)$ be 
proper, resp.\ 
surjective finite morphisms of good stratifications. 
Let $\Phi_U$ be an open subset of $\Phi$.
Denote by $\FS_U \subset \FS$ the pre-image of $\Phi_U$ under $\pi$,
and by $\FS_Z$ its complement.
Suppose that Assumption~\ref{As} holds.
Then the functor
\[
a_{Z,*} : \pi_* DMT_{\FS_Z}(S(\FS_Z))_{F,w = 0}^\natural \longto 
(a \circ \pi)_* DMT_{\FS_Z}(S(\FS_Z))_{F,w = 0}^\natural
\]
is radicial, and the diagram
\[
\vcenter{\xymatrix@R-10pt{
    \pi_* DMT_{\FS_U}(S(\FS_U))_{F,w = 0}^\natural 
                    \ar[d]_{a_{U,*}} \ar@{^{ (}->}[r]^-{\ujast} &
    \pi_* DMT_{\FS}(S(\FS))_{F,w=0}^{\natural,u} \ar[d]_{a_*^u} \\
    (a \circ \pi)_* DMT_{\FS_U}(S(\FS_U))_{F,w = 0}^\natural 
                    \ar@{^{ (}->}[r]^-{\ujast} &
    (a \circ \pi)_* DMT_{\FS}(S(\FS))_{F,w=0}^{\natural,u}
\\}}
\]
commutes. 
\end{Thm}

\begin{Proof}
When applied to 
$\pi_Z: S(\FS_Z) \to Y(\Phi_Z)$
and $a_Z: Y(\Phi_Z) \to W(\Phi_Z)$, Corollary~\ref{6F}
tells us that $a_{Z,*}$ is indeed radicial . 
Note that according to Main Theorem~\ref{5A}
and Corollary~\ref{5F},
source and target of $a_{Z,*}$ are semi-primary.

Now apply Proposition~\ref{6E}~(a).
\end{Proof}

The proof of the following result is along the same lines;
we leave the details to the reader.

\begin{Thm} \label{6H}
Let $\pi: S(\FS) \to Y(\Phi)$ and $b: Y'(\Phi) \to Y(\Phi)$ be 
proper, resp.\ smooth morphisms of good stratifications. 
Denote by $\pi': S'(\FS) \to Y'(\Phi)$ the base change of $\pi$
to $Y'(\Phi)$. Suppose that Assumption~\ref{As1} holds.
Then the functor
\[
b_Z^*: \pi_* DMT_{\FS_Z}(S(\FS_Z))_{F,w = 0}^\natural \longto 
\pi'_* DMT_{\FS_Z}(S'(\FS_Z))_{F,w = 0}^\natural 
\]
is radicial, and the diagram
\[
\vcenter{\xymatrix@R-10pt{
    \pi_* DMT_{\FS_U}(S(\FS_U))_{F,w = 0}^\natural 
                    \ar[d]_{b_U^*} \ar@{^{ (}->}[r]^-{\ujast} &
    \pi_* DMT_{\FS}(S(\FS))_{F,w=0}^{\natural,u} \ar[d]_{b^{*,u}} \\
    \pi'_* DMT_{\FS_U}(S'(\FS_U))_{F,w = 0}^\natural 
                    \ar@{^{ (}->}[r]^-{\ujast} &
    \pi'_* DMT_{\FS}(S'(\FS))_{F,w=0}^{\natural,u}
\\}}
\]
commutes. 
\end{Thm}

\begin{Rem} \label{6I}
Theorem~\ref{6H} applies in particular to morphisms
$b$ induced by smooth changes
$\BB' \to \BB$ of finite type 
of our base scheme (see the conventions fixed in our Introduction). 
Actually, thanks to \emph{continuity} \cite[Prop.~14.3.1]{CD}, 
the statement remains valid
without the finiteness condition on the base change. This is true
for example when $\BB' \to \BB$ is the inclusion of a generic point of $\BB$. 
\end{Rem}


\bigskip
%
%

\section{Compatibility with realizations}
\label{7}



Throughout this section, we fix a generic point $\Spec k$ of our base scheme
$\BB$. For any scheme $X$ in the sense of the conventions fixed in the Introduction, 
the base change
$X \times_\BB \Spec k$ is separated and of finite type over $k$.
We assume that we are in one of the following situations. 
\begin{enumerate}
\item[(i)] $k$ is embedded into $\BC$ \emph{via} a morphism 
$\eta: k \into \BC$, yielding a geometric point
of $\BB$, denoted by the same symbol
\[
\eta: \Spec \BC \longto \Spec k 
\longinto \BB \; . 
\]
The \emph{Betti realization} is defined in \cite[D\'ef.~2.1]{Ay2}. 
It is a family of covariant exact functors
\[
R_{\eta,Y}: \SH(Y) \longto D(Y) \; ,
\] 
indexed by quasi-projective $k$-schemes $Y$. 
The source of $R_{\eta,Y}$ is
the \emph{stable homotopy category of $Y$-schemes} \cite[Sect.~4.5]{Ay}.
Its target is 
the (unbounded) derived category of the Abelian category 
of sheaves with values in Abelian groups
on the topological space $Y(\BC)$ of points 
of $Y$ with values in $\BC$ with respect to $\eta$. The functors $R_{\eta,Y}$ 
are symmetric monoidal \cite[Lemme~2.2]{Ay2}.
Accor\-ding to \cite[Prop.~2.4, Thm.~3.4, Thm.~3.7]{Ay2}, they
commute with the functors $f^*$, $f_*$, $f_!$, $f^!$, provided the latter
are applied to constructible objects (note that commutation holds without
this restriction for the two functors $f^*$ and $f_!$). 
In parti\-cu\-lar, they commute with Tate twists.
In \cite[Ex.~17.1.7]{CD},
it is shown how to obtain from the $R_{\eta,Y}$ a family of exact functors
with analogous properties, and which we denote by the same symbols
\[
R_{\eta,Y} : \DBcYM \longto D^b_c(Y) \; ,
\]     
where the right hand denotes the full triangulated sub-category
of $D(Y)$ of classes of bounded complexes with constructible
cohomology objects. The construction can be imitated to obtain  
$F$-linear versions of the Betti realization. Composing with the 
base change \emph{via} $\Spec k \into \BB$, we finally obtain a family
of exact tensor functors
\[
R_{\eta, X} : \DBcFXM \longto D^b_c(X_k)_F \; ,
\]  
still referred to as the Betti realization,   
and indexed by schemes $X$, whose base change
$X_k := X \times_\BB \Spec k$ is quasi-projective over $k$. The $R_{\eta, X}$
are symmetric monoidal; in particular, they respect the unit objects. 
They commute with the functors 
$f^*$, $f_*$, $f_!$, $f^!$ since $\Spec k \into \BB$ is a projective limit
of open immersions (use \cite[Prop.~14.3.1]{CD}),
\item[(ii)] $k$ is of characteristic zero, and $\ell$ is a prime. 
The \emph{$\ell$-adic realization} is defined in 
\cite[Sect.~7.2, see in part.\ Rem.~7.2.25]{CD2}. 
It is a family of covariant exact functors
\[
R_{\ell,Y}: \DBcYM \longto D^b_c(Y) \; ,
\] 
indexed by $k$-schemes $Y$ of finite dimension. Its target is 
the bounded ``derived category'' of
constructible $\BQ_\ell$-sheaves on $Y$ \cite[Sect.~6]{E}.
The functors $R_{\ell,Y}$ 
are symmetric monoidal, and they
commute with the functors $f^*$, $f_*$, $f_!$, $f^!$ \cite[Thm.~7.2.24]{CD2}. 
In parti\-cular, they commute with Tate twists.
The construction can be imitated to obtain  
$F$-linear versions of the $\ell$-adic realization. Composing with the 
base change \emph{via} $\Spec k \into \BB$, we finally obtain a family
of exact tensor functors
\[
R_{\ell, X} : \DBcFXM \longto D^b_c(X_k)_F \; ,
\]  
still referred to as the $\ell$-adic realization. The $R_{\ell, X}$
are symmetric monoidal, and they
commute with the functors $f^*$, $f_*$, $f_!$, $f^!$.
\end{enumerate}

\begin{Ass} \label{As2}
Let $\pi: S(\FS) \to Y(\Phi)$ be a 
proper morphism of good stratifications of schemes $S(\FS)$ and $Y(\Phi)$.
In the context of the Betti rea\-li\-zation,
assume that the base change $Y(\Phi)_k$ is quasi-projective
over $\Spec k$.
Let $\Phi_U$ be an open subset of $\Phi$, and denote
by $\Phi_Z$ its complement. Write $\FS_U := \pi^{-1} \Phi_U$ and
$\FS_Z := \pi^{-1} \Phi_Z$. \\[0.1cm]
(a)~For all $\sigma \in \FS$, the strata
$\Ss$ are nilregular, and for any 
$\iat: \St \into \bSs$, the functor $i_\tau^!$
maps $\one_{\bSs}$ to a Tate motive over $S_\tau$. \\[0.1cm]
(b)~For all $\varphi \in \Phi_Z$, 
the strata $\Yp$ are nilregular. \\[0.1cm]
(c)~For all $\varphi \in \Phi_Z$ and $\sigma \in \FS_Z$
such that $\Ss$ is a stratum
of $\pi^{-1}(\Yp)$, the morphism
$\pis : \Ss \to \Yp$ can be factorized,
\[
\pis = \pios \circ \pits : \Ss \stackrel{\pits}{\longto} \Bs 
\stackrel{\pios}{\longto} \Yp \; ,
\]
such that condition~\ref{5A}~$(2)_\sigma$ is satisfied, 
such that $\pios$ is proper and smooth, 
and such that the pull-back of the base change 
\[
\pikos := \pios \times_\BB \Spec k: \Bks \longto \Ykp  
\] 
to any geometric
point of $\Ykp$ lying over a generic point
is isomorphic to a finite disjoint union of Abelian varieties.
\end{Ass}

Assumption~\ref{As2} is a variant of Assumption~\ref{As1}; it is more
restrictive in order to respect the hypothesis from \cite{Ay2}
on quasi-projectivity,
and also to guarantee the validity of the main result of this section, 
which we state now. \\

Let 
\[
R_{\bullet} : \DBcM(\bullet)_F \longto D^b_c(\bullet_k)_F 
\]  
be one of the two realizations considered above (Betti or $\ell$-adic).
The categories $D^b_c(X_k)_F$ 
are equipped with a perverse $t$-structure;
write $\Perv_c (X_k,F)$ for its heart, and 
\[
H^m : D^b_c(X_k)_F \longto \Perv_c (X_k,F) \; , \; m \in \BZ \; ,
\] 
for the perverse cohomology functors. 
Denote by 
\[
j_{!*}: \Perv_c(Y(\Phi_U)_k,F) \into \Perv_c(Y(\Phi)_k,F)
\]
the intermediate extension of perverse sheaves \cite[D\'ef.~1.4.22]{BBD}.

\begin{Thm} \label{7C}
Let $\pi: S(\FS) \to Y(\Phi)$ be a 
proper morphism of good stratifications of schemes $S(\FS)$ and $Y(\Phi)$.
Let $\Phi_U$ be an open subset of $\Phi$.
Denote by $\FS_U \subset \FS$ the pre-image of $\Phi_U$ under $\pi$,
and by $\FS_Z$ and $\Phi_Z$ the respective complements.
Suppose that Assumption~\ref{As2} holds. \\[0.1cm]
(a)~For any integer $m$, the restriction of the composition
\[
H^m \circ R_{Y(\Phi)}: \DBcFYPM \longto \Perv_c (Y(\Phi)_k,F) 
\]
to $\pi_* DMT_{\FS}(S(\FS))_{F,w =0}^\natural$
factors over $\pi_* DMT_{\FS}(S(\FS))_{F,w=0}^{\natural,u} \,$. \\[0.1cm]
(b)~For any integer $m$, the diagram
\[
\vcenter{\xymatrix@R-10pt{
    \pi_* DMT_{\FS_U}(S(\FS_U))_{F,w = 0}^\natural 
              \ar[d]_{H^m \circ R_{Y(\Phi_U)}} \ar@{^{ (}->}[r]^-{\ujast} &
    \pi_* DMT_{\FS}(S(\FS))_{F,w=0}^{\natural,u} \ar[d]^{H^m \circ R_{Y(\Phi)}} \\
    \Perv_c(Y(\Phi_U)_k,F)
                    \ar@{^{ (}->}[r]^-{j_{!*}} &
    \Perv_c(Y(\Phi)_k,F)
\\}}
\]
commutes. 
\end{Thm}

The rest of this section is devoted to the proof of Theorem~\ref{7C}.
As was the case for the proof of the main results of Section~\ref{6},
it necessitates the study of the statement under gluing. Let us return again to 
the abstract categorical setting. We fix two triples of $F$-linear
triangulated cate\-go\-ries, which 
are related by gluing: $\CC_1(U)$, $\CC_1(X)$, $\CC_1(Z)$, 
and $\CC_2(U)$, $\CC_2(X)$, $\CC_2(Z)$ (we use the same symbols
$j^*$, $i_*$ \emph{etc.} for the two sets of gluing functors).
The $\CC_1(\argdot)$ are assumed pseudo-Abelian, and equipped 
with a weight structure $w$, but this time
the $\CC_2(\argdot)$ are supposed to carry a $t$-structure $t$. Both $w$ and
$t$ are assumed to be compatible with the gluing. Denote by $\CC_2(\argdot)^{t = 0}$ 
the heart of the latter, and by $H^m: \CC_2(\argdot) \to \CC_2(\argdot)^{t = 0}$, $m \in \BZ$,
the cohomology functors associated to $t$. Furthermore, $F$-linear exact functors 
\[
r_U: \CC_1(U) \longto \CC_2(U) \; , \;
r_X: \CC_1(X) \longto \CC_2(X) \; , \;
r_Z: \CC_1(Z) \longto \CC_2(Z) 
\]
are supposed to be given, and to commute with the gluing functors. 
The following should be seen as an analogue of Proposition~\ref{6E}~(b)
in the absence of a weight structure on $\CC_2(\argdot)$.

\begin{Prop} \label{7D}
Assume that $\CC_1(Z)_{w = 0}$ is semi-primary, and
that for any object $M \in \CC_1(X)_{w = 0}$ without non-zero direct factors
in $i_* \CC_1(Z)_{w = 0}$, and any integer $m$, both 
\[
H^0 i^! H^m (r_X M) \quad \text{and} \quad H^0 i^* H^m (r_X M)
\]
are zero.
Denote the restriction of $r_\bullet$ to $\CC_1(\argdot)_{w = 0}$ by
$r_{\bullet,w=0}$.
If both $r_{Z,w=0}$ and $r_{U,w=0}$ map the radical to the 
kernel of $\oplus_{m \in \BZ} H^m$, 
then so does $r_{X,w=0}$. 
\end{Prop}

\begin{Proof}
Let $M$ and $N$ be objects of $\CC_1(X)_{w = 0}$, and
\[
f : M \longto N
\]
a morphism belonging to the radical.
According to Theorem~\ref{2K}~(a), we may assume that $M$ belongs to
the image of $j_{!*}$ (meaning that it does not admit non-zero
direct factors in the image of $i_*$) or to the image of $i_*$,
and likewise for $N$. 

If both $M$ and $N$ are in the image of $i_*$, then the claim is identical to
our assumption on $r_{Z,w=0}$. Now let us treat the case
where $M$ has no non-zero direct factors in the image of $i_*$. 
Since $j^* f$ is in the radical of $\CC_1(U)_{w = 0}$, 
\[
j^* H^m (r_X f) = H^m (r_U (j^* f)) = 0 
\]
for all $m \in \BZ$,
by our assumption on $r_{U,w=0}$. It follows that for all $m$, the morphism $H^m (r_X f)$
in $\CC_2(X)^{t = 0}$ factors through $i_* H^0 i^* H^m (r_X M) = 0$.
The case where $N$ is without non-zero direct factors
in $i_* \CC_1(Z)_{w = 0}$, is treated dually.
\end{Proof}

Now recall \cite[Prop.~2.1.2~1]{Bo} that any covariant additive functor
$\CH$ from a triangulated category $\CC$ carrying a weight structure $w$,
to an Abelian category $\FA$ admits a canonical \emph{weight filtration}
by sub-functors
\[
\ldots \subset W_n \CH \subset W_{n+1} \CH \subset \ldots \subset \CH \; .
\]
According to \cite[Def.~2.1.1]{Bo} (use the normalization of \cite[Def.~1.3.1]{Bo3} for the
signs of the weights), for an object $M$ of $\CC$, and $n \in \BZ$,
the sub-object $W_n \CH(M) \subset \CH(M)$ is defined as the image of the morphism
$\CH(\iota_{w \le n})$, for \emph{any} weight filtration
\[
M_{w \le n} \stackrel{\iota_{w \le n}}{\longto} M \longto M_{w \ge n+1} \longto M_{w \le n}[1]
\]
(with $M_{w \le n} \in \CC_{w \le n}$ and $M_{w \ge n+1} \in \CC_{w \ge n+1}$).  
For any $m \in \BZ$, one defines 
\[
\CH^m : \CC \longto \FA \; , \; M \longto \CH(X[m]) \; ;
\]
according to the usual convention, the weight filtration of $\CH^m(M)$ \emph{equals} the
weight filtration of $\CH(X[m])$, \emph{i.e.}, it differs by \emph{d\'ecalage}
from the intrinsic weight filtration of the covariant additive functor $\CH^m$. The following is a direct consequence of the definitions.

\begin{Lem} \label{7E}
Let $r: \CC_1 \to \CC_2$
be an $F$-linear exact functor between $F$-linear
triangulated cate\-go\-ries. Assume that $\CC_1$ is equipped 
with a weight structure $w$, and that $\CC_2$ carries a $t$-structure $t$.
Let 
\[
\alpha \in \Hom_{\CC_1} (B,A) \; ,
\]
for objects $A \in \CC_{1,w \le 0}$ and $B \in \CC_{1,w \ge 0}$.
Assume the following for all integers $m$: (i)~the morphism $H^m \circ r (\alpha)$ in $\CC_2^{t = 0}$ 
is strict with respect to the weight filtrations of $H^m \circ r (B)$ and $H^m \circ r (A)$, 
(ii)~$W_m H^m \circ r (\alpha) = 0$. \\[0.1cm]
Then $r (\alpha)$ is zero on cohomology: $H^m \circ r (\alpha) = 0$, $\forall \, m \in \BZ$.
\end{Lem}

\begin{Lem} \label{7F}
Assume that $\CC_1(Z)_{w = 0}$ is semi-primary, 
and that $r_{Z,w=0}$ maps the radical of $\CC_1(Z)_{w = 0}$ to the 
kernel of $\oplus_{m \in \BZ} H^m$.
Let $M \in \CC_1(X)_{w = 0}$ be without non-zero direct factors
in $i_* \CC_1(Z)_{w = 0}$. Assume that for all $m \in \BZ$, the functor
$H^m \circ r_Z$ maps the composition 
\[
i^! M \longto i^* M
\]
of the two adjunction morphisms
to a morphism which is strict with respect to the weight filtrations.
Then the composition
\[
i^! (r_X M) \longto i^* (r_X M)
\]
is zero on cohomology.
\end{Lem}

\begin{Proof}
Apply Lemma~\ref{7E} to $\CC_1 := \CC_1(Z)$, $\CC_2 := \CC_2(Z)$, and 
$\alpha : i^! M \to i^* M$ equal to the composition of the adjunction morphisms. 
This is indeed possible:
first, hypothesis~\ref{7E}~(i) is satisfied by assumption. 
Second, note that according to Summary~\ref{2M}~(a)~(3c), the composition 
\[
L_0 \longto i^! M \stackrel{\alpha}{\longto} i^* M \longto N_0
\]
belongs to the radical of $\CC_1(Z)_{w = 0}$, for any pair of objects $L_0$, $N_0$
occurring in weight filtrations of $i^! M$ and $i^* M$, respectively.  
Our assumption on $r_{Z,w=0}$ thus ensures
the validity of hypothesis~\ref{7E}~(ii).

Lemma~\ref{7E} then tells us that the morphism $r_Z (\alpha)$ is zero on cohomology. But
since the functors $r_{\argdot}$ commute with the gluing,
the morphism $r_Z (\alpha)$ equals the composition 
\[
i^! (r_X M) \longto i^* (r_X M)
\]
of the adjunction morphisms. 
\end{Proof}

In the context of the realization
\[
R_{\bullet} : \DBcM(\bullet)_F \longto D^b_c(\bullet_k)_F 
\]  
we now have to check the hypotheses of Proposition~\ref{7D} and Lemma~\ref{7F}.
For the Betti realization $R_{\eta,\bullet}$, it will be necessary to prove a number of
technical results, culminating in Corollary~\ref{7H}, which states that
the morphism $i^! M \to i^* M$ from Lemma~\ref{7F} is indeed strict.
We start with the following partial compatibility result.

\begin{Prop} \label{7Fa}
Let $f: S \to X$ be a proper morphism of schemes.
Assume that the base change $X_k$ is quasi-projective. 
Then for all integers $m$, the weight filtration step $W_{m-1}$ of \cite[Def.~2.1.1]{Bo} on
the perverse sheaf
\[
H^m \circ R_{\eta,X} (f_* \one_S) 
\]
on $X(\BC)$ coincides with the weight filtration step $W_{m-1}$ underlying the 
algebraic mixed Hodge module
\[
H^m (f_* \BQ_{S_\BC}) \; .
\]
\end{Prop}

Here, we make use of the functor
\[
D^b \bigl( \MHM_\BQ (\argdot \times_\BB \Spec \BC) \bigr) \longto
D^b_c(\bullet_k)_F 
\]
from the bounded derived category $D^b ( \MHM_\BQ (\argdot \times_\BB \Spec \BC) )$ 
of algebraic mixed Hodge modules on $\argdot \times_\BB \Spec \BC$ \cite[Sect.~4.2]{Sa}. By definition, it is the composition of the forgetful functor
\[
D^b \bigl( \MHM_\BQ (\argdot \times_\BB \Spec \BC) \bigr) \longto
D^b_c(\bullet_k)_\BQ 
\]
and the functor 
\[
D^b_c(\bullet_k)_\BQ \longto
D^b_c(\bullet_k)_F
\]
given by extension of coefficients.
It is $t$-exact with respect to the obvious $t$-structure on 
$D^b ( \MHM_\BQ (\argdot \times_\BB \Spec \BC) )$ and the perverse
$t$-structure on $D^b_c(\bullet_k)_F$, and it commutes with
the Grothendieck operations $f^*$, $f_*$, $f_!$, $f^!$ on both sides \cite[Thm.~4.3, Sect.~4.4]{Sa}.

\medskip

\begin{Proofof}{Proposition~\ref{7Fa}}
Without loss of generality, we may assume that $\BB = \Spec k$,
and that $k = \BC$. 

Suppose first that $S$ is regular. Then $f_* \one_S$
is a Chow motive (Theorem~\ref{3A}~(1), (2)), hence the weight filtration
on $H^m \circ R_{\eta,X} (f_* \one_S)$ is concentrated in weight $m$.
But the same is then true for the Hodge theoretical weight filtration
\cite[Thm.~5.3.1]{Sa1}.

In general case, let $s: \tilde{S} \to S$ be a surjective, proper morphism
with regular source. Consider the adjunction morphism
\[
adj: \one_S \longto s_* \one_{\tilde{S}} \; .
\]
Choose $C \in \DBcM(S)_F$ fitting into an exact triangle
\[
C \longto \one_S \stackrel{adj}{\longto} s_* \one_{\tilde{S}} \longto C[1] \; .
\]
Our claim follows from the next Lemma.
\end{Proofof}

\begin{Lem} \label{7lem}
(a)~The motive $C$ lies in $\DBcM(S)_{F,w \le -1}$. The exact triangle
\[
f_* C \longto f_* \one_S \stackrel{f_* adj}{\longto} (fs)_* \one_{\tilde{S}} \longto f_* C[1] 
\]
is a weight filtration of $f_* \one_S$. \\[0.1cm]
(b)~For all integers $m$, the Hodge theoretical weight filtration step $W_{m-1}$ on
the perverse sheaf $H^m \circ R_{\eta,X} (f_* \one_S)$ coincides with the kernel of the morphism
induced by $f_* adj$, 
\[
H^m \circ R_{\eta,X} (f_* \one_S) \longto H^m \circ R_{\eta,X} ((fs)_* \one_{\tilde{S}}) \; .
\]
\end{Lem}

\begin{Proof}
The second statement of (a) follows from the first, since $f$ is proper (Theorem~\ref{3A}~(2)).

For the remaining claims, let us start by assuming $S$ to be nilregular. In that case, 
we imitate the proof of \cite[Prop.~4.2.6.1]{A} in order to show that 
the morphism $adj$ is split monomorphic. This, together with stability of
purity under direct images of proper morphisms
\cite[Thm.~5.3.1]{Sa1} establishes both (a) and (b).

In the general case, we use induction on the dimension $d$ of $S$.
When $d=0$, then $S$ is nilregular.
Thus assume $d>0$, and choose a dense, nilregular open subscheme $j: U \into S$.
Denote by $i: Z \into S$ the complementary closed immersion ($Z$ carries the
reduced scheme structure, say). 

In order to prove claim~(a), consider the localization triangles
for $\one_S$ and for $s_* \one_{\tilde{S}}$.
\[
\vcenter{\xymatrix@R-10pt{
        j_! \one_U \ar[d]_{j_! adj} \ar[r] &
        \one_S \ar[d]_{j_! adj} \ar[r] &
        i_* \one_Z \ar[d]_{i_* adj} \ar[r]^-{[1]} & \\
        j_! s_* \one_{\tilde{U}} \ar[r] &
        s_* \one_{\tilde{S}} \ar[r] &
        i_* s_* \one_{\tilde{Z}} \ar[r]^-{[1]} &
\\}}
\]
Here, we denote by $\tilde{U}$ and $\tilde{Z}$ the pre-images under $s$ of $U$ and $Z$,
respectively. Complete this morphism of exact triangles to a diagram with exact rows
and columns. 
\[
\vcenter{\xymatrix@R-10pt{
        j_! C_U  \ar[d] \ar[r] &
        C \ar[d] \ar[r] &
        i_* C_Z \ar[d] \ar[r]^-{[1]} & \\
        j_! \one_U \ar[d]_{j_! adj} \ar[r] &
        \one_S \ar[d]_{j_! adj} \ar[r] &
        i_* \one_Z \ar[d]_{i_* adj} \ar[r]^-{[1]} & \\
        j_! s_* \one_{\tilde{U}} \ar[d]_-{[1]} \ar[r] &
        s_* \one_{\tilde{S}} \ar[d]_-{[1]} \ar[r] &
        i_* s_* \one_{\tilde{Z}} \ar[d]_-{[1]} \ar[r]^-{[1]} & \\
        &&&
\\}}
\] 
According to the nilregular case, $C_U$ is of weights $\le -1$, hence so is 
$j_! C_U$ (Theorem~\ref{3A}~(2)). In order to limit the weights occurring 
in $C_Z$, hence in $i_* C_Z$, let $z: \tilde{\tilde{Z}} \to \tilde{Z}$ be
surjective and proper, with nilregular source. The composition $sz: \tilde{\tilde{Z}} \to Z$
is also surjective and proper. We complete the diagram
\[
\vcenter{\xymatrix@R-10pt{
        \one_Z \ar[d]_{adj} \ar[r]^-{adj} &
        (sz)_* \one_{\tilde{\tilde{Z}}} \ar@{=}[d]  \\
        s_* \one_{\tilde{Z}} \ar[r]^-{s_* adj} &
        (sz)_* \one_{\tilde{\tilde{Z}}}
\\}}
\]
to a diagram of exact triangles,
and apply the induction hypothesis to both $adj: \one_Z \to (sz)_* \one_{\tilde{\tilde{Z}}}$ and 
$adj: \one_{\tilde{Z}} \to z_* \one_{\tilde{\tilde{Z}}}$. Thus, we see
that indeed $C_Z \in \DBcM(Z)_{F,w \le -1}$.

In order to prove (b), note that the images under $R_{\eta,X}$ of the above diagrams
are isomorphic to the images under the forgetful functor of the Hodge theoretical
analogues of the diagrams. We may thus repeat the reasoning, observing that 
the analogues of the inequalities for weights from Theorem~\ref{3A}~(2)
hold in the Hodge theoretical context thanks to \cite[(4.5.2)]{Sa}.
\end{Proof}

\begin{Rem}
(a)~Lemma~\ref{7lem}~(b) (is most probably known to the experts and)
generalizes the first statement of \cite[Prop.~(8.2.5)]{D1},
which concerns the case $X = \Spec \BC$. \\[0.1cm]
(b)~Using simplicial resolutions as in \cite[Sect.~8.1]{D1}, it can be shown that the 
whole weight filtration on
\[
H^m \circ R_{\eta,X} (f_* \one_S) 
\]
coincides with the Hodge theoretical weight filtration; this holds more generally for
separated morphisms $f$. This result
should be compared to \cite[Prop.~2.2.2]{LW}, 
which concerns the case when $X = \BB = \Spec k$ and $f$ is smooth. 
In the sequel, we shall only use Proposition~\ref{7Fa} as stated.
\end{Rem}

\begin{Cor} \label{7H}
Let $X$ be a scheme, whose base change
$X_k$ is quasi-projective over $k$. 
Let $i: Z \into X$ be the immersion of a closed subscheme
of $X$. Let $M \in \CHFXM$. Then for all $m \in \BZ$, the functor
$H^m \circ R_{\eta,Z}$ maps the composition 
\[
i^! M \longto i^* M
\]
of the two adjunction morphisms
to a morphism which is strict with respect to the weight filtrations.
\end{Cor}

\begin{Proof}
Our claim is clearly stable under
passage to direct factors and finite direct sums. 
Given the description of objects 
of $\CHFXM$ from Theo\-rem~\ref{3A}~(3), it therefore suffices to check
the claim for motives $M$ of the form
\[
f_* \one_S (p)[2p] \; ,
\]
for a projective morphism $f: S \to X$ with nilregular source $S$, 
and an integer $p$. The claim is invariant under twists and shifts,
meaning that we may assume $p=0$. We may also assume $F = \BQ$.
Furthermore, the base change
$S_k$ is quasi-projective over $k$. 
The functors $R_{\eta,\bullet}$ commute with $f_*$, $i^!$ and $i^*$;
therefore, $R_{\eta,Z}$ maps the morphism $i^! M \to i^* M$ to
the composition 
\[
i^! f_* (R_{\eta,S} \one_S) \longto i^* f_* (R_{\eta,S} \one_S)
\]
of the adjunction maps in $D^b_c(Z(\BC),\BQ)$.
But this composition lies in the image of the forgetful functor
\[
D^b ( \MHM_\BQ (Z \times_\BB \Spec \BC) ) \longto
D^b_c(Z_k)_\BQ \; .
\]
Altogether, this establishes that for all $m \in \BZ$, the functor
$H^m \circ R_{\eta,Z}$ maps the composition 
\[
i^! f_* \one_S \longto i^* f_* \one_S
\]
of the two adjunction morphisms
to a morphism $\alpha_m$ which is strict with respect to the weight filtrations
of Hodge theory. 

It remains to establish the analogous statement for the weight filtrations
of \cite[Def.~2.1.1]{Bo}. 

Observe first (Theorem~\ref{3A}~(2), \cite[(4.5.2)]{Sa})
that the perverse sheaf
$H^m \circ R_{\eta,Z} (i^! f_* \one_S)$ is of weights $\ge m$, whichever of the two
notions of weight is chosen. Likewise, $H^m \circ R_{\eta,Z} (i^* f_* \one_S)$
is of weights $\le m$. Strictness of
\[
\alpha_m: H^m \circ R_{\eta,Z} (i^! f_* \one_S) \longto H^m \circ R_{\eta,Z} (i^* f_* \one_S)
\]
is thus equivalent to the following statement: the image of $\alpha_m$ coincides
with the image of the restriction of $\alpha_m$ to $W_m H^m \circ R_{\eta,Z} (i^! f_* \one_S)$.
This statement being already established for the Hodge theoretical weight filtration,
we need to show that the filtration steps $W_m$ coincide on 
$H^m \circ R_{\eta,Z} (i^! f_* \one_S)$. 

Second, apply duality for motives \cite[Thm.~15.2.4~(d)]{CD}
and for Hodge modules \cite[(4.2.3)]{Sa}. It shows that the desired claim
is equivalent to the following: the filtration steps $W_{m-1}$ coincide on 
$H^m \circ R_{\eta,Z} (i^* f_* \one_S)$ (recall that $f$ is projective,
hence proper). By proper base change \cite[Thm.~2.4.50~(4)]{CD},
\[
i^* f_* \one_S = f'_* \one_{S_Z} \; ,
\]
where $f': S_Z \to Z$ is the base change of $f$ \emph{via} $i$. 
Now apply Proposition~\ref{7Fa}.
\end{Proof}

Let us return to the simultaneous treatment of the Betti and $\ell$-adic
rea\-li\-zation
\[
R_{\bullet} : \DBcM(\bullet)_F \longto D^b_c(\bullet_k)_F \; .
\]  
The following is a direct consequence of the Decomposition
Theorem \cite[Thm.~6.2.5]{BBD}.

\begin{Thm} \label{7G}
Let $X$ be a scheme. 
In the context of the Betti rea\-li\-zation,
assume that the base change $X_k$ is quasi-projective over $\Spec k$.
Let $M \in \CHFXM$. Then $R_X M$ has the following properties. \\[0.1cm]
(a)~Let $\alpha$ be one of the four operations $f^*$, $f_*$, $f_!$, $f^!$
associated to a morphism $f$ with source $X$ (for $\alpha = f_*$ or $f_!$) or
target $X$ (for $\alpha = f^*$ or $f^!$). Then the spectral sequence
\[
E_2^{p,q} = H^p \alpha H^q (R_X M) \Longrightarrow H^{p+q} \alpha (R_X M) 
\]
of perverse sheaves degenerates at $E_2$. \\[0.1cm]
(b)~Let $i: Z \into X$ a closed immersion. Then the composition 
\[
H^0 i^! H^m (R_X M) \longto H^0 i^* H^m (R_X M)
\]
of the two adjunction morphisms is an isomorphism of perverse sheaves on $Z_k$.
\end{Thm} 

\begin{Proof}
The desired properties of $R_X M$ are stable under
passage to direct factors. 
According to Theorem~\ref{3A}~(3), we may therefore assume $M$ to equal
$f_* \one_S$, 
for a projective morphism $f: S \to X$ with nilre\-gu\-lar source $S$. 

In the context of the Betti realization, 
part~(a) follows from the commutation of $R_\bullet$ with $f_*$ and   
from the Decomposition Theorem \cite[Thm.~6.2.5]{BBD}: indeed,
there exists an isomorphism
\[
\beta: R_X M \isoto \bigoplus_{m \in \BZ} H^m(R_X M)[-m]
\]
in $D^b_c(X_k)_F$ such that $H^m \beta = \id_{H^m(R_X M)}$, 
for all $m \in \BZ$. As for part~(b), note that by
the Decomposition Theorem \cite[Thm.~6.2.5]{BBD},
the perverse sheaves $H^m(R_X M)$ are direct sums
of simple perverse sheaves, for all $m \in \BZ$.
Now use
the description of simple perverse sheaves on $X_k(\BC)$ \cite[Prop.~1.4.26]{BBD},
and \cite[Cor.~1.4.25]{BBD}.

The same argument works for the $\ell$-adic realization, provided that $k$ is algebraically closed
(see the remark following \cite[Thm.~6.2.5]{BBD}). For the descent to 
an arbitrary field $k$ (of characteristic zero), use the fact that base change
from $\Perv_c (\bullet_k,F)$ to $\Perv_c (\bullet_{\bar{k}},F)$ 
is faithful (hence (a)) and conservative (hence (b)).
\end{Proof}

\forget{
The following is a direct consequence of the Decomposition
Theorem \cite[Thm.~6.2.5]{BBD}.

\begin{Thm} \label{7G}
Let $X$ be a scheme, whose base change
$X \times_\BB \Spec k$ is quasi-projective over $k$. Then the essential image
of the restriction of 
\[
R_X : \DBcFXM \longto D^b_c(X(\BC),F) 
\]
to $\CHFXM$ is contained in the full sub-category of $D^b_c(X(\BC),F)$
of \emph{semi-simple} objects: for all $M \in \CHFXM$, 
the following holds. \\[0.1cm]
(a)~The object $R_X M$ is \emph{split}, \emph{i.e.},
there exists an isomorphism
\[
\beta: R_X M \isoto \bigoplus_{n \in \BZ} H^n(R_X M)[-n]
\]
in $D^b_c(X(\BC),F)$ such that $H^n \beta = \id_{H^n(R_X M)}$, 
for all $n \in \BZ$. \\[0.1cm]
(b)~For all $n \in \BZ$, the perverse sheaf $H^n(R_X M)$ is a direct sum
of simple perverse sheaves.
\end{Thm} 

\begin{Proof}
An object $K$ of $D^b_c(X(\BC),F)$ is split if and only if for all integers
$m$, the boundary morphism $\delta$ in the exact truncation triangle
\[
\tau_{\le n-1} K \longto \tau_{\le n} K \longto H^n(K)[-n] 
\stackrel{\delta}{\longto} (\tau_{\le n-1} K)[1] 
\]
is trivial. This shows that the property of being semi-simple is stable under
passage to direct factors. 
Given the description of objects 
of $\CHFXM$ from Theorem~\ref{3A}~(3), it therefore suffices check
the claim for motives $M$ of the form
\[
f_* \one_S (p)[2p] \; ,
\]
for a projective morphism $f: S \to X$ with nilre\-gu\-lar source $S$, 
and an integer $p$. The claim is invariant under twists and shifts,
meaning that we may assume $p=0$. 
Now use the commutation of $R_\bullet$ with $f_*$,    
and the Decomposition Theorem \cite[Thm.~6.2.5]{BBD}.
\end{Proof}

\begin{Cor} \label{7H}
Let $X$ be a scheme, whose base change
$X \times_\BB \Spec k$ is quasi-projective over $k$. 
Let $i: Z \into X$ be the immersion of a closed subscheme
of $X$. Let $M \in \CHFXM$. Then for all $m \in \BZ$, the functor
$H^m \circ R_Z$ maps the composition 
\[
i^! M \longto i^* M
\]
of the two adjunction morphisms
to a morphism which is strict with respect to the weight filtrations.
\end{Cor}

\begin{Proof}
Our claim is clearly stable under
passage to direct factors. 
As in the proof of Theorem~\ref{7G}, we may therefore assume $M$ to equal
$f_* \one_S$, 
for a projective morphism $f: S \to X$ with nilregular source $S$.
Furthermore, the base change
$S \times_\BB \Spec k$ is quasi-projective over $k$. 
The functors $R_\bullet$ commute with $f_*$, $i^!$ and $i^*$;
therefore, $R_Z$ maps the morphism $i^! M \to i^* M$ to
the composition 
\[
i^! f_* (R_S \one_S) \longto i^* f_* (R_S \one_S)
\]
of the adjunction maps in $D^b_c(Z(\BC),F)$.
But this composition lies in the image of the forgetful functor
\[
D^b ( \MHM_\BQ (Z \times_\BB \Spec \BC) ) \longto
D^b_c(Z(\BC),F) \; .
\]
Altogether, this establishes that for all $m \in \BZ$, the functor
$H^m \circ R_Z$ maps the composition 
\[
i^! f_* \one_S \longto i^* f_* \one_S
\]
of the two adjunction morphisms
to a morphism which is strict with respect to the weight filtrations
of Hodge theory.

It remains to establish the analogous statement for the weight filtrations
of \cite[Def.~2.1.1]{Bo}. First, we apply Theorem~\ref{7G}~(a), which tells us
that  
\[
R_X (f_* \one_S) \cong \bigoplus_{n \in \BZ} H^n(R_X (f_* \one_S))[-n] \; .
\]
Thus,
\[
R_Z (i^! f_* \one_S) = i^! R_X (f_* \one_S) 
\cong \bigoplus_{n \in \BZ} i^! H^n(R_X (f_* \one_S))[-n] \; ,
\]
\[
R_Z (i^* f_* \one_S) = i^* R_X (f_* \one_S) 
\cong \bigoplus_{n \in \BZ} i^* H^n(R_X (f_* \one_S))[-n] \; ,
\]
and the morphism
\[
H^m \circ R_Z (i^! f_* \one_S \longto i^* f_* \one_S)
\]
corresponds to the direct sum over $n \in \BZ$ of the morphisms
\[
H^{m-n} \circ \bigl( i^! H^n(R_X (f_* \one_S)) \longto i^* H^n(R_X (f_* \one_S)) \bigr) \; .
\]
Second, observe that $H^n(R_X (f_* \one_S))$ being a perverse sheaf on $X(\BC)$,
the complexes $i^! H^n(R_X (f_* \one_S))$ and $i^* H^n(R_X (f_* \one_S))$
are concentrated in degrees $\ge 0$
and $\le 0$, respectively \cite[Prop.~1.4.16~(a)]{BBD}. 
Thus, we may assume that $m=n$, and need to study the morphism 
\[
H^0 i^! H^n(R_X (f_* \one_S)) \longto H^0 i^* H^n(R_X (f_* \one_S)) \; .
\]

that the weight filtrations of Hodge theory on the
\[
H^m \circ r_Z (i^! f_* \one_S) \quad \text{and} \quad
H^m \circ r_Z (i^* f_* \one_S) \; , \; m \in \BZ \; , 
\]
coincide with the weight filtrations of the functor $H^m \circ r_Z$.
By duality for motives \cite[Cor.~14.3.31~(d)]{CD}
and for Hodge modules \cite[\?]{Sa}, it suffices to prove the claim
for $H^m \circ r_Z (i^* f_* \one_S)$ (since $f$ is projective,
hence proper). By proper base change \cite[Thm.~2.4.50~(4)]{CD},
\[
i^* f_* \one_S = f'_* \one_{S_Z} \; ,
\]
where $f': S_Z \to Z$ is the base change of $f$ \emph{via} $i$. 
Now apply Proposition~\ref{7Fa}.
\end{Proof}
}

Let us come back to the morphism $\pi: S(\FS) \to Y(\Phi)$.
The prece\-ding results allow us to establish the key technical ingredient 
of the proof of Theorem~\ref{7C}, as far as gluing is concerned. 

\begin{Cor} \label{7I}
Let $\pi: S(\FS) \to Y(\Phi)$ be a 
proper morphism of good stratifications of schemes $S(\FS)$ and $Y(\Phi)$.
In the context of the Betti rea\-li\-zation,
assume that the base change $Y(\Phi)_k$ is quasi-projective.
Let $\Phi_U$ be an open subset of $\Phi$.
Denote by $\FS_U \subset \FS$ the pre-image of $\Phi_U$ under $\pi$,
and by $\FS_Z$ and $\Phi_Z$ the respective complements.
Denote by $i: Y(\Phi_Z) \into Y(\Phi)$ the closed immersion. 
Suppose the validity of Assumption~\ref{As1}. Suppose also that 
the restriction of the functor
\[
R_{Y(\Phi_Z)} : \DBcFYZM \longto D^b_c(Y(\Phi_Z)_k)_F
\]
to the sub-category $\pi_* DMT_{\FS_Z}(S(\FS_Z))_{F,w = 0}^\natural$ of $\DBcFYZM$
maps the ra\-di\-cal of $\pi_* DMT_{\FS_Z}(S(\FS_Z))_{F,w = 0}^\natural$ to 
the kernel of all perverse cohomology functors 
\[
H^m: D^b_c(Y(\Phi_Z)_k)_F \longto \Perv_c (Y(\Phi_Z)_k,F) \; , 
\; m \in \BZ \; .
\]
(a)~Let $M \in \pi_* DMT_{\FS}(S(\FS))_{F,w = 0}^\natural$ without non-zero factors in 
the image of $i_*$.
Then for any integer $m$, both 
\[
H^0 i^! H^m (R_{Y(\Phi)} M) \quad \text{and} \quad 
H^0 i^* H^m (R_{Y(\Phi)} M)
\]
are zero. \\[0.1cm]
(b)~If both 
\[
R_{Y(\Phi_Z),w=0}: \pi_* DMT_{\FS_Z}(S(\FS_Z))_{F,w = 0}^\natural \longto 
D^b_c(Y(\Phi_Z)_k)_F
\] 
and 
\[
R_{Y(\Phi_U),w=0}: \pi_* DMT_{\FS_U}(S(\FS_U))_{F,w = 0}^\natural \longto 
D^b_c(Y(\Phi_U)_k)_F
\]
map the radical to the kernel 
of all perverse cohomology functors $H^m$, then so does
\[
R_{Y(\Phi),w=0}: \pi_* DMT_{\FS}(S(\FS))_{F,w = 0}^\natural \longto 
D^b_c(Y(\Phi)_k)_F \; .
\]
\end{Cor}

\begin{Proof}
According to Theorem~\ref{7G}~(b), 
\[
H^0 i^! H^m (R_{Y(\Phi)} M) \isoto H^0 i^* H^m (R_{Y(\Phi)} M) \; .
\]
Thanks to Assumption~\ref{As1},
Main Theorem~\ref{5A} can be applied; therefore,
the category $\pi_* DMT_{\FS_Z}(S(\FS_Z))_{F,w = 0}^\natural$ is semi-primary.

We claim that the assumptions of Lemma~\ref{7F} are met, in other words, that 
$H^m \circ R_{Y(\Phi_Z)}$ maps the composition 
$i^! M \to i^* M$
of the two adjunction morphisms
to a morphism which is strict with respect to the weight filtrations.
For the Betti realization, this is the content of Corollary~\ref{7H}.
For the $\ell$-adic realization, 
we have a much more general statement: \emph{any} morphism in 
the image of $H^m \circ R_{\ell,Y(\Phi_Z)}$ is strict
with respect to the weight filtrations
\cite[Thm.~2.5.4~(II)~(1), Prop.~1.3.2~(II)~(2)]{Bo4}.
Note that $Y(\Phi_Z)_k$ 
is of finite type over $k$, hence \emph{very reasonsable}
in the sense of \cite[Def.~2.1.1~(4)]{Bo4}.

Thus, according to Lemma~\ref{7F}, the morphism
\[
H^m i^! (R_{Y(\Phi)} M) \longto H^m i^* (R_{Y(\Phi)} M)
\]
is zero. Thanks to Theorem~\ref{7G}~(a), this implies the isomorphism
\[
H^0 i^! H^m (R_{Y(\Phi)} M) \isoto H^0 i^* H^m (R_{Y(\Phi)} M) 
\]
is zero, too, which shows claim~(a). 

As for claim~(b), note that the assumptions of
Proposition~\ref{7D} are all satisfied. 
\end{Proof}

The main point that remains to be addressed, concerns the hypothesis
``$R_{Y(\Phi_\bullet),w=0}(\rad) \subset \ker H^m$''.
The following is a reformulation of the main result from \cite{Li}.

\begin{Thm} \label{7K}
Let $p_1: B_1 \to X$ and $p_2: B_2 \to X$ be proper, smooth morphisms
of schemes. Assume that the pull-back of 
\[
p_{r,k} := p_r \times_\BB \Spec k: B_{r,k} \longto X_k \; ,
\; r = 1,2 \; ,
\] 
to any geometric
point of $X_k$ lying over a generic point
is isomorphic to a finite disjoint union of Abelian varieties. 
In the context of the Betti rea\-li\-zation,
assume that the base change $X_k$ is quasi-projective. Then 
the restriction of the functor
\[
R_X : \DBcFXM \longto D^b_c(X_k)_F
\]
to $\CHFXM$ maps the ra\-di\-cal 
\[
\rad_{\CHFXM} \bigl( p_{1,*} \one_{B_1} , p_{2,*} \one_{B_2} \bigr) \subset
\Hom_{\CHFXM} \bigl( p_{1,*} \one_{B_1} , p_{2,*} \one_{B_2} \bigr)
\]
to the kernel of all perverse cohomology functors 
\[
H^m: D^b_c(X_k)_F \longto \Perv_c (X_k,F) \; , 
\; m \in \BZ \; .
\]
\end{Thm}

\begin{Proof}
Given that the $p_r$ are proper and smooth,
the perverse cohomology objects of the 
$R_X (p_{r,*} \one_{B_r})$ are, up to a shift, local
systems (in the context of the Betti realization) and lisse $\ell$-adic sheaves
(in the context of the $\ell$-adic realization), respectively. A morphism of such 
is trivial as soon as it is generically trivial.
Therefore, we may assume that $X = \Spec k$.
In particular (Example~\ref{5E}~(c)),
the smooth Chow motives $p_{r,*} \one_{B_r}$ are finite
dimensional.

The radical
\[
\rad_{\CHFXM} \bigl( p_{1,*} \one_{B_1} , p_{2,*} \one_{B_2} \bigr)
\]
then coincides with the maximal tensor ideal 
\[
\CN_X \bigl( p_{1,*} \one_{B_1} , p_{2,*} \one_{B_2} \bigr)
\]
\cite[Thm.~9.2.2]{AK}. It thus
consists of the classes of numerically trivial cycles on 
$B_1 \times_X B_2$ \cite[Ex.~7.1.2]{AK}. 
According to \cite[Thm.~4]{Li}, numerical and
homological equivalence coincide on Abelian varieties
over a field of characteristic zero. Given our hypotheses, 
\[
R_X \bigl( \rad_{\CHFXM} \bigl( p_{1,*} \one_{B_1} , 
                 p_{2,*} \one_{B_2} \bigr) \bigr)
\]   
therefore consists of morphisms which are zero on cohomology.
\end{Proof}

\begin{Cor} \label{7L}
Let $\pi: S(\FS) \to Y(\Phi)$ be a 
proper morphism of good stratifications $S(\FS)$ and $Y(\Phi)$.
Assume that all strata $\Yp$ and $\Ss$, $\varphi \in \Phi$, $\sigma \in \FS$,
are nilregular, and that for any 
$\iat: \St \into \bSs$, the functor $i_\tau^!$
maps $\one_{\bSs}$ to a Tate motive over $S_\tau$. 
Furthermore, assume 
that for all $\sigma \in \FS$
such that $\Ss$ is a stratum
of $\pi^{-1}(\Yp)$, the morphism
$\pis : \Ss \to \Yp$ can be factorized,
\[
\pis = \pios \circ \pits : \Ss \stackrel{\pits}{\longto} \Bs 
\stackrel{\pios}{\longto} \Yp \; ,
\]
such that conditions~\ref{5A}~$(1)_\sigma$ and $(2)_\sigma$ are satisfied,
and such that the pull-back of the base change 
\[
\pikos: \Bks \longto \Ykp  
\] 
to any geometric
point of $\Ykp$ lying over a generic point
is isomorphic to a finite disjoint union of Abelian varieties.
In the context of the Betti rea\-li\-zation,
assume that the base change $Y(\Phi)_k$ is quasi-projective.
Then the functor
\[
R_{Y(\Phi),w=0}: \pi_* DMT_{\FS}(S(\FS))_{F,w = 0}^\natural \longto 
D^b_c(Y(\Phi)_k)_F
\]
maps the ra\-di\-cal of $\pi_* DMT_{\FS}(S(\FS))_{F,w = 0}^\natural$ to 
the kernel of all perverse cohomology functors 
\[
H^m: D^b_c(Y(\Phi)_k)_F \longto \Perv_c (Y(\Phi)_k,F) \; , 
\; m \in \BZ \; .
\]
\end{Cor}

A proper morphism $\pi$ is \emph{of Abelian type} 
if it satisfies the hypotheses of Corollary~\ref{7L}. A
\emph{Chow motive of Abelian type over $Y(\Phi)$} is a motive isomorphic to
an object of $\pi_* DMT_{\FS}(S(\FS))_{F,w = 0}^\natural \,$, for a proper morphism 
$\pi: S(\FS) \to Y(\Phi)$ of Abelian type.

\medskip

\begin{Proofof}{Corollary~\ref{7L}}
According to Corollary~\ref{4I}~(b), 
the category $\pi_* \! \DFTSnatSM$ is obtained by gluing. 
Thanks to Corollary~\ref{7I}~(b), we are reduced to showing that each 
\[
R_{\Yp,w=0}: \pi_* DMT_{\FS}(\pi^{-1}(\Yp))_{F,w=0}^\natural \longto
D^b_c(\Ykp)_F
\]
maps the ra\-di\-cal of $\pi_* DMT_{\FS}(\pi^{-1}(\Yp))_{F,w = 0}^\natural$ to 
the kernel of all perverse cohomology functors 
$H^m$. By Proposition~\ref{5AL}, all objects of the source of $R_{\Yp,w=0}$
are direct factors of a finite direct sum of objects isomorphic to
$\pi'_{\sigma,*} \one_{\Bs}(p)[2p]$, for $p \in \BZ$ and $\sigma \in \FS$,
such that $\Ss$ is a stratum of $\pi^{-1}(\Yp)$. 
Our claim thus follows from our assumption on $\pios$,
the commutation of $R_{\Yp,w=0}$ with Tate twists,
and from Theorem~\ref{7K}.
\end{Proofof}

\begin{Proofof}{Theorem~\ref{7C}}
When applied to $\pi_Z: S(\FS_Z) \to Y(\Phi_Z)$,
Corollary~\ref{7L} tells us that the restriction of the functor
\[
R_{Y(\Phi_Z)} : \DBcFYZM \longto D^b_c(Y(\Phi_Z)_k)_F
\]
to the sub-category $\pi_* DMT_{\FS_Z}(S(\FS_Z))_{F,w = 0}^\natural$ of $\DBcFYZM$
maps the ra\-di\-cal of $\pi_* DMT_{\FS_Z}(S(\FS_Z))_{F,w = 0}^\natural$ to 
the kernel of all perverse cohomology functors 
\[
H^m: D^b_c(Y(\Phi_Z)_k)_F \longto \Perv_c (Y(\Phi_Z)_k,F) \; .
\]
(a): Recall (Definition~\ref{1Ea}) that the category
$\pi_* DMT_{\FS}(S(\FS))_{F,w=0}^{\natural,u}$
is the quotient $\pi_* DMT_{\FS}(S(\FS))_{F,w=0}^\natural / \Fg$, where
$\Fg$ is the two-sided ideal generated by
\[
\Hom_{\pi_* DMT_{\FS}(S(\FS))_{F,w=0}^\natural} (M,i_*N) \quad \text{and} \quad 
\Hom_{\pi_* DMT_{\FS}(S(\FS))_{F,w=0}^\natural} (i_*N,M) \; ,
\]
for all objects $(M,N)$ of 
$\pi_* DMT_{\FS}(S(\FS))_{F,w=0}^\natural \times \pi_* DMT_{\FS_Z}(S(\FS_Z))_{F,w=0}^\natural \,$, 
such that $M$ admits no non-zero direct factor in the image of $i_*$. 
In order to prove the claim, it thus suffices to show that for all
morphisms $f: M \to i_* N$ and $g: i_* N \to M$, for $M$ and $N$ as above,
the morphisms $R_{Y(\Phi)} f$ and $R_{Y(\Phi)} g$ are zero on cohomology.
But this follows from Corollary~\ref{7I}~(a). 

(b): Let $M_U \in \pi_* DMT_{\FS_U}(S(\FS_U))_{F,w = 0}^\natural$. Then
$j_{!*} M_U$ is an extension
of $M_U$ without non-zero direct factors in the image of $i_*$.
According to Corollary~\ref{7I}~(a), for any integer $m$, 
the perverse sheaves
\[
H^0 i^! H^m (R_{Y(\Phi)} j_{!*} M_U) \quad \text{and} \quad 
H^0 i^* H^m (R_{Y(\Phi)} j_{!*} M_U)
\]
are zero. Therefore \cite[Cor.~1.4.24]{BBD}, 
\[
H^m (R_{Y(\Phi)} j_{!*} M_U) = j_{!*} j^* H^m (R_{Y(\Phi)} j_{!*} M_U)
= j_{!*} H^m (R_{Y(\Phi_U)} M_U) \; .
\]  
This shows commutation of the diagram on the level of objects.
Now let $f$ be a morphism in $\pi_* DMT_{\FS_U}(S(\FS_U))_{F,w = 0}^\natural$.
In order to show the relation
\[
(H^m \circ R_{Y(\Phi)}) (\ujast f) 
= j_{!*} \circ H^m (R_{Y(\Phi_U)} f) \; ,
\] 
in $\Perv_c(Y(\Phi)_k,F)$, observe that the relation holds
trivially after restriction to $Y(\Phi_U)_k$. Thus, both 
$(H^m \circ R_{Y(\Phi)}) (\ujast f)$ and 
$j_{!*} \circ H^m (R_{Y(\Phi_U)} f)$ extend the same morphism.
But their source and target are intermediate extensions; therefore,
morphisms between their restrictions extend uniquely.
\end{Proofof}

\begin{Rem} \label{7M}
(a)~In the context of the Betti realization,
a more conceptual, and much less complex proof of Theorem~\ref{7C} could be
given if the Hodge theoretical realization 
\[
R_{{\bf H} , X} : \DBcFXM \longto 
D^b \bigl( \MHM_\BQ  (X \times_\BB \Spec \BC) \otimes_\BQ F \bigr) 
\] 
were known to exist, to yield the Betti realization $R_{\eta,X}$ after composition
with the forgetful functor, and to be compatible with the formalism of
six operations. Indeed, the target of $R_{{\bf H} , X}$ 
carries a weight structure (Example~\ref{2D}~(b)).
Proposition~\ref{6E} would then replace everything leading up to 
our main gluing result Corollary~\ref{7I} 
(Proposition~\ref{7D} -- Theorem~\ref{7G}). 
Given Corollary~\ref{7L}, the functor $R_{{\bf H} , X}$
can be expected to be radicial on Chow motives, 
given that the radical of its target consists precisely of the morphisms,
which are zero on cohomology. \emph{Vice versa}, given the interpretation
of \cite[Thm.~4]{Li} in the proof of Theorem~\ref{7K}, to conjecture
$R_{{\bf H} , X}$ to be radicial on the heart of the motivic weight structure
seems to be the correct generalization
of the conjecture ``numerical equivalence equals (Betti)
homological equivalence''
to Chow motives over arbitrary bases. \\[0.1cm]
(b)~The Hodge theoretical realization $R_{{\bf H} , X}$ left its ``shadow'' 
in the present proof of Theorem~\ref{7C}, \emph{via} (a)~strictness
of the Hodge theoretical weight filtration on perverse sheaves
of geometric origin (proof of Corollary~\ref{7H}), (b)~the Decomposition
Theorem (proof of Corollary~\ref{7I}), which should be consi\-de\-red as a 
Hodge theoretical phenomenon \cite[(4.5.4)]{Sa}. \\[0.1cm]
(c)~As far as the $\ell$-adic realization is concerned, we are unable
to say much in positive characteristic. In fact, a partial analogue of
\cite[Thm.~4]{Li} has been proved over algebraic
closures of finite fields \cite[Thm.~1]{C} (see also \cite{D}).
This result could be used to prove a version of Theorem~\ref{7C}
only if Assumption~\ref{As2}~(b) were replaced by the following
much more restrictive hypothesis: for all $\varphi \in \Phi_Z$, 
the strata $\Yp$ are spectra of finite products of fields, which
are algebraic over some $\BF_p$, for $p \ne \ell$. \\[0.1cm]
(d)~Assume the base scheme $\BB$ to be arithmetic.
Let $\pi: S(\FS) \to Y(\Phi)$ be a 
proper morphism of good stratifications satisfying Assumption~\ref{As2},
and such that the generic fiber $Y(\Phi) \otimes_\BZ \BQ$ is proper over $\BB \otimes_\BZ \BQ$.
Assume compatibi\-li\-ty of the intermediate extension with the $p$-adic 
realization. One may then employ \cite[Thm.~2~2)]{KM}
(\emph{e.g.} as in the proof of \cite[Thm.~1.2.4~(ii)]{Sch}), 
to see that for primes $p$ of good reduction, the characteristic polynomials
of a Frobenius on $p$-adic intersection cohomology of $Y(\Phi) \otimes_\BZ \BQ$
and on the associated $\phi$-filtered module coincide. 
This principle will be exploited elsewhere. \\[0.1cm]
(e)~At this point, the reader will have noted that the title of the present work
is somewhat misleading. Indeed, our main results on the existence and the properties
of the motivic intermediate extension (Corollary~\ref{5Aa}, Theorems~\ref{6G}, \ref{6H}
and \ref{7C}) only require $\pi_Z: S(\FS_Z) \to Y(\Phi_Z)$ to be of Abelian type.
Our theory thus works for Chow motives ``admitting a degeneration of Abelian type''.
In particular, it may be expected to apply to compactifications of moduli spaces
of $K3$-surfaces. 
\end{Rem}


\bigskip
%
%

\section{Applications}
\label{8}



Our main application of the material developed in the preceding sections
concerns Shimura varieties and their compactifications. We shall
recall the data necessary for their construction, and simultaneously check
the validity of Assumptions~\ref{As1}, \ref{As} and \ref{As2} 
(Lemmata~\ref{8A}--\ref{8C}). \\

Our notation is identical to the one used in \cite{P1,P2,W1,BW}.  
Let $(P, \FX)$ be \emph{mixed Shimura data} \cite[Def.~2.1]{P1}. 
In particular, $P$ is a connected algebraic linear group over $\BQ$, and
$P(\BR)$ acts on the complex ma\-ni\-fold $\FX$ by analytic automorphisms.
Denote by $W$ the unipotent radical of $P$. If $P$ is reductive, \emph{i.e.}, 
if $W=0$, then $(P, \FX)$ is called \emph{pure}.
The \emph{Shimura varieties} associated to $(P,\FX)$ are indexed by the 
open compact subgroups of $P (\BA_f)$,
where $\BA_f$ denotes the ring of finite ad\`eles over $\BQ$.
If $K$ is one such, then the analytic space of $\BC$-valued points of the
corresponding variety is given as
\[
M^K (\BC) := P (\BQ) \backslash ( \FX \times P (\BA_f) / K ) \; .
\]
According to
Pink's gene\-ral\-ization to mixed Shimura varieties of the Algebraization
Theorem of Baily and Borel
\cite[Prop.~9.24]{P1}, there exist canonical structures of normal quasi-projective 
varieties on the $M^K(\BC)$, which we denote as
\[
M^K_{\BC} := M^K (P, \FX)_{\BC} \; .
\]
According to
\cite[Thm.~2.18]{M} and \cite[Thm.~11.18]{P1}, there is a \emph{canonical model}
of $M^K_{\BC}$, which we denote as
\[
M^K := M^K (P, \FX) \; .
\]
It is defined over the \emph{reflex field} $E (P , \FX)$ of $(P , \FX)$ 
\cite[11.1]{P1}. \\

Any \emph{admissible parabolic subgroup} \cite[Def.~4.5]{P1} $Q$ of $P$
has a canonical normal subgroup $P_1$ \cite[4.7]{P1}.
There is a finite collection of \emph{rational boun\-dary components} 
$(P_1 , \FX_1)$, indexed by the $P_1 (\BR)$-orbits in $\pi_0 (\FX)$ 
\cite[4.11]{P1}. The $(P_1 , \FX_1)$ are themselves mixed Shimura data. 
Consider the following condition on $(P, \FX)$ (cmp.\ \cite[Condition~(3.1.5)]{P2}).
\begin{enumerate}
\item [$(+)$] If $G$ denotes the maximal reductive quotient of $P$, 
then the neutral connected component $Z (G)^0$ of the center $Z (G)$ of 
$G$ is, up to isogeny, a direct product of a $\BQ$-split torus with a torus 
$T$ of compact type (\emph{i.e.}, $T(\BR)$ is compact) defined over $\BQ$.
\end{enumerate}
If $(P, \FX)$ satisfies $(+)$, 
then so does any rational boundary component $(P_1 , \FX_1)$
\cite[proof of Cor.~4.10]{P1}.
Similarly, the reflex field does not change when passing from $(P , \FX)$ to a 
rational boundary component \cite[Prop.~12.1]{P1}. \\

Fix a connected component $\FX^0$ of $\FX$, and denote by
$(P_1 , \FX_1)$ the associated rational boundary component. 
Denote by $U_1 \subset P_1$ the ``weight $-2$'' part of $P_1$.
It is Abelian, normal in $Q$, and central in the unipotent radical 
$W_1$ of $P_1$.
\cite[4.11, 4.14, Prop.~4.15]{P1} give the construction of an open convex cone
\[
C (\FX^0 , P_1) \subset U_1 (\BR) (-1)
\]
canonically associated to $P_1$ and $\FX^0$. \\

Set-theoretically, the \emph{co\-ni\-cal complex} associated to $(P,\FX)$ is defined as
\[
\CC (P , \FX) := \coprod_{(\FX^0 , P_1)} C (\FX^0 , P_1) \; .
\]
By \cite[4.24]{P1}, the conical complex is naturally equipped with a 
topology (usually different from the coproduct topology). 
The closure $C^{\ast} (\FX^0 , P_1)$ of $C (\FX^0 , P_1)$ 
inside $\CC (P , \FX)$
can still be considered as a convex cone in $U_1 (\BR) (-1)$, 
with the induced topology. \\

The \emph{toroidal compactifications} of 
$M^K$ are para\-me\-terized by \emph{$K$-admissible complete cone 
decompositions}, which are collections of subsets of
\[
\CC (P,\FX) \times P (\BA_f)
\]
satisfying the axioms \cite[6.4~(i'), (ii), (iii), (iv), (v)]{P1}. 
If $\FS$ is one such, then in particular any member of $\FS$ is of the shape
\[
\sigma \times \{ p \} \; ,
\]
$p \in P (\BA_f)$, $\sigma \subset C^{\ast} (\FX^0 , P_1)$ a
convex rational polyhedral cone in 
the vector space $U_1 (\BR) (-1)$ 
not containing any non-trivial linear subspace. \\

Let $M^K (\FS) := M^K (P , \FX , \FS)$ be the associated compactification;
we refer to \cite[Thm.~9.21, 9.27, Prop.~9.28]{P1}
for criteria sufficient to guarantee its existence.
It comes equipped with a natural stratification into locally closed strata.
Any such stratum is obtained as follows: 
fix a pair $(\FX^0 , P_1)$ as above, $p \in P (\BA_f)$ and
\[
\sigma \times \{ p \} \in \FS
\]
such that $\sigma \subset C^{\ast} (\FX^0 , P_1)$,
$\sigma \cap C (\FX^0 , P_1) \neq \emptyset$.
To $\sigma$, one associates Shimura data
\[
\left( \Pes , \Xes \right)
\]
\cite[7.1]{P1}, whose underlying group 
$\Pes$ is the quotient of $P_1$ by the algebraic subgroup
\[
\langle \sigma \rangle \subset U_1
\]
satisfying 
$\BR \cdot \sigma = \frac{1}{2 \pi i} \cdot \langle \sigma \rangle (\BR)$. Set 
\[
K_1  :=  P_1 (\BA_f) \cap p \cdot K \cdot p^{-1} \; , \quad
\pi_{[\sigma]}  :  P_1 \longonto \Pes \; .
\]
According to \cite[7.3]{P1}, there is a canonical map
\[
\isp (\BC) : \Sps (\Pes , \Xes) (\BC) \longrightarrow 
M^K (\FS) (\BC) := M^K (P , \FX , \FS) (\BC)
\]
whose image is locally closed. 
The latter is disjoint from $M^K (\BC)$
if and only if the admissible parabolic subgroup $Q$ giving rise to $P_1$ is \emph{proper},
i.e., unequal to $P$. 
$\isp (\BC)$
comes from a morphism of schemes over $E (P , \FX)$
\cite[Thm.~12.4~(a), (c)]{P1}, denoted 
\[
\isp : \Sps := \Sps (\Pes , \Xes) \longto M^K (\FS) \; .
\]
Its image equals
the stratum associated to $\sigma \times \{ p \}$.
Letting $\sigma \times \{ p \}$ vary, the natural stratification 
\[
M^K (\FS) = \coprod_{\sigma \times \{ p \}} \isp (\Sps)
\]
is indexed by a quotient of $\FS$, which is finite \cite[7.3]{P1}. 
By abuse of notation, it will be denoted by the same letter $\FS$.
Similarly, we shall write $\sigma$ for the class of $\sigma \times \{ p \}$. 
If $(P, \FX)$ satisfies $(+)$, and $K$ is \emph{neat} 
(see e.g.\ \cite[0.6]{P1}), then $\isp$ is an immersion \cite[Prop.~1.6]{W1}, \emph{i.e.}, it 
identifies $\Sps$ with a locally closed sub-scheme of $M^K (\FS)$. 
In that case, the natural stratification takes the form
\[
M^K (\FS) = \coprod_{\sigma \in \FS} \Sps \; .
\]

\begin{Lem} \label{8A}
(a)~The stratification $M^K (\FS) = \coprod_{\sigma \in \FS} \isp (\Sps)$ is good,
\emph{i.e.}, the closure $\overline{\isp (\Sps)}$ of any stratum
is a union of strata. \\[0.1cm]
(b)~Assume that $(P, \FX)$ satisfies $(+)$, and that $K$ is neat. 
Then modulo a possible replacement of $\FS$ by a suitable refinement, 
$M^K (\FS)$ is a projective variety over $E(P,\FX)$, and 
the closures $\bSps$ of strata, $\sigma \in \FS$ are smooth over $E (P , \FX)$.
\end{Lem}

\begin{Proof}
By \cite[7.11]{P1}, 
\[
\isp (\BC) : \Sps (\BC) \longrightarrow 
M^K (\FS) (\BC) 
\]
extends to a continuous map
\[
\overline{\isp} (\BC) : \Sps (\FS_{1,[\sigma]}) (\BC) \longrightarrow 
M^K (\FS) (\BC) 
\]
of toroidal compactifications.
Here, $\FS_{1,[\sigma]}$ is the cone decomposition 
defined in \cite[7.7]{P1}. Since $\FS$ is assumed to be complete, the image of 
$\overline{\isp} (\BC)$ equals the closure of $\isp (\Sps) (\BC)$ \cite[Prop.~6.27, Prop.~7.8]{P1}.

By construction \cite[7.11]{P1}, the extension
$\overline{\isp} (\BC)$ maps any stratum of $\Sps (\FS_{1,[\sigma]}) (\BC)$
onto a stratum of $M^K (\FS) (\BC)$. This shows part~(a), since the claim can be checked 
on $\BC$-valued points.

Let us turn to the proof of part~(b). 
According to \cite[proof of Thm.~9.21]{P1}, a suitable refinement $\FS'$ of 
$\FS$ is smooth with respect to $K$ \cite[6.4]{P1}, satisfies condition~\cite[7.12~($\argast$)]{P1},
and the associated toroidal compactification $M^K(\FS')$ (exists and)
is projective over $E(P,\FX)$.
Let us change notation, and write $\FS$ instead of $\FS'$.  

Together with condition $(+)$ and neatness of $K$,
smoothness of $\FS$ implies that the cone decomposition $\FS_{1,[\sigma]}$
is smooth, too. Indeed, the groups $\Gamma_U \subset U_2 (\BQ)$ defined in \cite[6.4]{P1},
for any boundary component $(P_2 , \FX_2)$ of $(P , \FX)$, and any element $p \in P (\BA_f)$,
are then equal to $U_2 (\BQ) \cap pKp^{-1}$, \emph{i.e.}, they do not change when passing
from $\FS$ to the cone decomposition denoted 
\[
\FS_1 := \bigl( [\ \cdot p]^* \FS \bigr) _{(P_1 , \FX_1)}
\]
in \cite[7.3]{P1}. 

By \cite[Prop.~6.26]{P1}, $\Sps (\FS_{1,[\sigma]}) (\BC)$ is smooth.

Thanks to condition~\cite[7.12~($\argast$)]{P1}, we may apply
\cite[Cor.~7.17~(a)]{P1}. It tells us that the image of $\overline{\isp} (\BC)$
equals the quotient of $\Sps (\FS_{1,[\sigma]}) (\BC)$ by the action of a certain group denoted
$\Stab_{\Delta_1} ([\sigma])$. But according to \cite[Lemma~1.7, Rem.~1.8]{W1},
this group is trivial, thanks to condition $(+)$ and neatness of $K$.
This establishes the claim, since smoothness of the $\bSps$ can be checked 
on $\BC$-valued points.
\end{Proof}

Now let $(G, \FH)$ be the pure Shimura data underlying $(P, \FX)$, \emph{i.e.}, let
$(G, \FH) := (P, \FX) /W$ be the quotient of $(P, \FX)$ by $W$ \cite[Prop.~2.9]{P1}. 
Fix an open compact subgroup $L \subset G (\BA_f)$.
Let $(M^L)^*$ denote the \emph{Baily--Borel compactification}
of $M^L := M^L(G,\FH)$. 
It is a projective variety over the reflex field $E(G, \FH) = E(P, \FX)$ \cite[Main Theorem~12.3~(a), (b)]{P1}.
It comes equipped with a natural stratification into locally closed strata.
Any stratum disjoint from $M^L$ is obtained as follows: 
fix a 
proper admissible parabolic subgroup $Q$ of $G$ with
associated normal subgroup $P_1$. Fix a
rational boundary
component $(P_1, \FX_1)$ of $(G, \FH)$,
and an element $g \in G(\BA_f)$.
Define
\[
L_1  :=  P_1 (\BA_f) \cap g \cdot L \cdot g^{-1} \; .
\]
Denote by $W_1$
the unipotent radical of $P_1$, and by 
\[
\pi_1: (P_1,\FX_1) \longonto (G_1,\FH_1) := (P_1,\FX_1)/W_1
\]
the quotient
of $(P_1 , \FX_1)$ by $W_1$.
According to \cite[7.6, Main Theorem~12.3~(c)]{P1}, there is
a canonical morphism
\[
\ig : \Spi := \Spi (G_1, \FH_1) \longto 
(M^L)^* 
\]
whose image is locally closed, and identical to the stratum 
associated to $(P_1 , \FX_1)$ and $g$. \\

Letting $(P_1 , \FX_1)$ and $g$ vary, the natural stratification
\[
(M^L)^* = \coprod_{(P_1 , \FX_1),g} \ig (\Spi)
\]
is indexed by a quotient $\Phi$ of $\{ (P_1 , \FX_1) \} \times G(\BA_f)$, which
is finite \cite[6.3]{P1}.

\begin{Lem} \label{8B}
(a)~The stratification $(M^L)^* = \coprod_\Phi \ig (\Spi)$ is good. \\[0.1cm]
(b)~Assume that $(G, \FH)$ satisfies $(+)$, and that $L$ is neat. 
Then all strata $\ig (\Spi)$, $((P_1,\FX_1),g) \in \Phi$ are smooth over $E (G , \FH)$.
\end{Lem}

\begin{Proof}
As in the proof of Lemma~\ref{8A}, the claims can be checked on 
the level of $\BC$-valued points.
By \cite[7.6]{P1}, 
\[
\ig (\BC): \Spi (\BC) \longto (M^L)^* (\BC)
\]
extends to a continuous map
\[
\overline{\ig} (\BC) : (\Spi)^* (\BC) \longto (M^L)^* (\BC)
\]
of Baily--Borel compactifications. Since $(\Spi)^* (\BC)$ is compact, the image of
$\overline{\ig} (\BC)$ equals the closure of $\ig (\Spi) (\BC)$. By construction 
\cite[7.6]{P1}, the extension
$\overline{\ig} (\BC)$ maps any stratum of $(\Spi)^* (\BC)$
onto a stratum of $(M^L)^* (\BC)$. This shows part~(a).

In order to prove part~(b), define the following group (cmp.\ \cite[Sect.~1]{BW}):
\[
H_Q := \Stab_{Q (\BQ)} (\FH_1) \cap P_1 (\BA_f) \cdot L' \; .
\]
$H_Q$ acts by analytic automorphisms on
$\FH_1 \times P_1 (\BA_f) / L_1$. Hence the group
$\Delta_1 := H_Q / P_1(\BQ)$ acts naturally on
\[
\Spi (\BC) =
P_1(\BQ) \backslash ( \FH_1 \times P_1 (\BA_f) / L_1 ) \; .
\]
By \cite[6.3]{P1}, the action of $\Delta_1$
factors through a finite quotient $\Delta$, and
the quotient by this action is precisely the image of $\ig (\BC)$.
By \cite[Prop.~1.1~(d)]{BW}, the action of $\Delta$ on $\Spi (\BC)$
is free provided that $(G, \FH)$ satisfies $(+)$, and $L$ is neat. 
But $\Spi (\BC)$ is smooth: indeed, $L$ being neat, so is $\pi (L_1)$.
Now apply \cite[Prop.~3.3~(b)]{P1}.
\end{Proof}

We need to analyze morphisms 
associated to change of the level structure $L$,
and their extension to the Baily-Borel compactifications. Thus, let $L'$ be a second
open compact subgroup of $G(\BA_f)$, and $h \in G(\BA_f)$ such that
$L \subset h \cdot L' \cdot h^{-1}$. Then there is a finite, surjective morphism
\[
[ \ \cdot h] : M^{L} \longto M^{L'}
\]
\cite[3.4~(a), Def.~11.5]{P1}. It extends to a finite, surjective morphism
\[
(M^{L})^* \longto (M^{L'})^*
\]
\cite[Main Theorem~12.3~(b)]{P1}, which we still denote by $[ \ \cdot h]$. 
By the description from \cite[6.3]{P1}, 
one sees that $[ \ \cdot h]: (M^{L})^* \to (M^{L'})^*$ is a morphism
of stratifications with respect to the natural stratifications 
\[
(M^{L})^* = \coprod_\Phi \ig (\Spi) \quad \text{and} \quad
(M^{L'})^* = \coprod_\Phi \ig (\Spti)
\]
from above. Fix $\Spi$ and $\Spti$ occurring in these stratifications.

\begin{Lem} \label{8Ba}
Assume that $(G, \FH)$ satisfies $(+)$, that ${L'}$ is neat (hence so is $L$),
and that $\ig(\Spi)$ is a stratum of $[ \ \cdot h]^{-1} (\ig (\Spti)$. 
Then the morphism induced by $[ \ \cdot h]$ on $\ig(\Spi)$,
\[
[ \ \cdot h]_g : \ig(\Spi) \longto \ig(\Spti)
\]
is finite and \'etale.
\end{Lem}

\begin{Proof}
The varieties $\Spi$ and $\Spti$ are associated to the same Shimura data
$(G_1,\FH_1)$. From the description of the morphisms $\ig$ on complex points
\cite[6.3]{P1}, it follows that $[ \ \cdot h]_g$ fits into a commutative diagram
\[
\vcenter{\xymatrix@R-10pt{
        M^{\tilde{L}_1}:= M^{\tilde{L}_1}(G_1,\FH_1) \ar[d]_{[\ \cdot h_1]} \ar@{>>}[r]^-{b'} &
        \Delta \backslash M^{\tilde{L}_1} \cong \ig(\Spi) \ar[d]^{[ \ \cdot h]_g} \\
        M^{\tilde{{L'}}_1}:= M^{\tilde{{L'}}_1}(G_1,\FH_1) \ar@{>>}[r]^-{b} &
        \Delta' \backslash M^{\tilde{{L'}}_1} \cong \ig(\Spti)
\\}}
\]
for suitable open compact subgroups $\tilde{L}_1$ and $\tilde{L'}_1$ of $G_1(\BA_f)$,
which are neat, an element $h_1$ of $G_1(\BA_f)$ such that 
$\tilde{L}_1 \subset h_1 \cdot \tilde{{L'}}_1 \cdot h_1^{-1}$, 
and finite groups $\Delta$ and $\Delta'$ acting on $M^{\tilde{L}_1}$ and $M^{\tilde{L'}_1}$, respectively (see the proof of Lemma~\ref{8B}). Their actions being free, the surjective 
morphisms $b$ and $b'$ are finite and \'etale. But so is the morphism $[\ \cdot h_1]$.
\end{Proof}

Let us now connect the two types of compactifications. Let $(P,\FX)$ be mixed Shimura
data as before, 
\[
\pi : (P,\FX) \longonto (G,\FH) = (P,\FX)/W \; ,
\]
$K$ an open compact subgroup of $P(\BA_f)$. 
Put $L := \pi(K)$. The morphism of Shimura data $\pi$ gives rise to a surjective morphism
\[
M^K = M^K (P,\FX) \longto M^L = M^L (G,\FH) 
\]
\cite[3.4~(b), Prop.~11.10]{P1}, equally denoted $\pi$. 
Let $\FS$ be a $K$-admissible complete cone decomposition. Then $\pi$ extends to a proper,
surjective morphism
\[
M^K(\FS) \longto (M^L)^*
\]
\cite[6.24, Main Theorem~12.4~(b)]{P1}, still denoted $\pi$.
From the description given in \cite[7.3]{P1}, one sees that $\pi: M^K(\FS) \to (M^L)^*$ is a morphism
of stratifications with respect to the natural stratifications 
\[
M^K (\FS) = \coprod_{\sigma \in \FS} \isp (\Sps) \quad \text{and} \quad
(M^L)^* = \coprod_\Phi \ig (\Spi)
\]
from above. Fix $\Sps$ and $\Spi$ occurring in these stratifications.

\begin{Lem} \label{8C}
Assume that $(P, \FX)$ satisfies $(+)$, that $K$ is neat, and that
$\Sps$ is a stratum of $\pi^{-1}(\ig (\Spi))$. 
Then the morphism induced by $\pi$ on $\Sps$,
\[
\pis : \Sps \longto \ig (\Spi)
\] 
can be factorized
\[
\pis = \pios \circ \pits : \Sps \stackrel{\pits}{\longto} \Bs 
\stackrel{\pios}{\longto} \ig (\Spi) \; ,
\]
such that the morphisms $\pios$ and $\pits$ satisfy the following conditions.
\begin{enumerate}
\item[$(1)$] $\pios$ is proper and smooth, and its pull-back to any geometric
point of $\ig (\Spi)$ is isomorphic to a finite disjoint union of Abelian varieties.
\item[$(2)$] The motive 
\[
\pi''_{\sigma,*} \one_{\Sps} \in \DBcFBsM
\]
belongs to the category $\DFTBsM$ of Tate motives over $\Bs$.
\end{enumerate}
\end{Lem} 

\begin{Proof}
The varieties $\Spi$ and $\Sps$ are associated to Shimura data $(G_1,\FH_1)$ and $(\Pes , \Xes)$,
respectively. Since $\Sps$ belongs to
$\pi^{-1}(\ig (\Spi))$, the Shimura data $(G_1,\FH_1)$ are
in fact the pure Shimura data underlying $(\Pes , \Xes)$. 
From the description of the morphisms $\ig$ and $\isp$ on complex points
\cite[6.3, 7.3]{P1}, it follows that $\pis$ is isomorphic to the composition
\[
M^{K_{1, [\sigma]}}:= M^{K_{1, [\sigma]}}(\Pes , \Xes) 
\stackrel{\pi_1}{\longto} M^{\pi_1 (K_{1, [\sigma]})} 
\stackrel{b}{\longonto} \Delta \backslash M^{\pi_1 (K_{1, [\sigma]})} \; ,
\]
where we let $\pi_1 : (\Pes , \Xes) \onto (G_1,\FH_1)$ denote the quotient map,
for a suitable open compact subgroup $K_{1, [\sigma]}$
of $\Pes (\BA_f)$, which is neat, and for a finite
group $\Delta$ acting on $M^{\pi_1 (K_{1, [\sigma]})}:= M^{\pi_1 (K_{1, [\sigma]})}(G_1,\FH_1)$
(see the proof of Lemma~\ref{8B}). 
The action of $\Delta$ being free,
the morphism $b$ is finite and \'etale.
The factorization of $\pi_1: (\Pes,\Xes) \onto (G_1,\FH_1)$
corresponding to the weight filtration $1 \subset \Ues \subset \Wes$ of the unipotent radical 
$\Wes$ of $\Pes$ gives the following:
\[
\vcenter{\xymatrix@R-10pt{
(\Pes,\Xes) \ar@{->>}[r]^-{\pi_t} \ar@/_2pc/[rr]^-{\pi_1}&
(P'_1,\FX'_1) := (\Pes,\Xes)/\Ues \ar@{->>}[r]^-{\pi_a} &
(G_1,\FH_1)
\\}}
\]
On the level of Shimura varieties, we get:
\[
\vcenter{\xymatrix@R-10pt{
M^{K_{1, [\sigma]}} \ar@{->>}[r]^-{\pi_t} \ar@/_2pc/[rr]^-{\pi_1} &
M^{\pi_t (K_{1, [\sigma]})}(P'_1,\FX'_1) \ar@{->>}[r]^-{\pi_a} &
M^{\pi_1 (K_{1, [\sigma]})}
\\}}
\]
By \cite[3.12--3.22~(a), Prop.~11.10]{P1}, 
$\pi_a$ is in a natural way a torsor under an
Abelian scheme, while $\pi_t$ is a torsor under a split torus $T$. Now put 
\[
B_\sigma:= M^{\pi_t (K_{1, [\sigma]})}(P'_1,\FX'_1) \quad , \quad 
\pits := \pi_t \quad \text{and} \quad \pios := b \circ \pi_a \; .
\]
\end{Proof}

The verification of Assumptions~\ref{As1}, \ref{As} and \ref{As2}
is now complete. The theory of the intermediate extension therefore
applies. In order to make it explicit,
assume that $(P, \FX)$ satisfies $(+)$, and that $K$ is neat.
As before, denote by $\pi$ the morphism $M^K(\FS) \to (M^L)^*$.
Let $U \subset (M^L)^*$ be an open union of strata, 
and denote by $\FS_U \subset \FS$ the set of strata of $M^K(\FS)$
mapping to $U$ under $\pi$. Assume that $\FS$ is fine enough for
Lemma~\ref{8A}~(b) to hold.

\begin{Thm} \label{8D} 
The hypotheses
of Theorem~\ref{2K}~(a) are fulfilled for $\CC(X) = \pi_* \! \DFTSMM$ and 
$\CC(U) = \pi_* \! \DFTSUM$. 
In particular, the intermediate extension
\[
\ujast : \pi_* DMT_{\FS_U} \bigl( \pi^{-1}(U) \bigr)_{F,w = 0}^\natural
\longinto \pi_* DMT_\FS \bigl( M^K(\FS) \bigr)_{F,w=0}^{\natural,u} 
\]
is defined, and it satisfies the properties listed in Summary~\ref{2M}. 
\end{Thm}

\begin{Proof}
The morphisms $\pi$ is proper and respects the stratifications,
which are good according to Lemmata~\ref{8A}~(a) and \ref{8B}~(a).
Assumption~\ref{As1}~(a) holds thanks to Lemma~\ref{8A}~(b).
So do Assumption~\ref{As1}~(b) (by Lemma~\ref{8B}~(b))
and Assumption~\ref{As1}~(c) (by Lemma~\ref{8C} and Example~\ref{5E}~(c)).

Now apply Corollary~\ref{5Aa}.
\end{Proof}

\medskip

\begin{Proofof}{Theorem~\ref{0A}}
This is the special case of a trivial ``weight $-2$'' part of $P$,
and of $U = M^L \subset (M^L)^*$. According to Theorem~\ref{8D},
\[
\ujast : \pi_* DMT \bigl( M^K \bigr)_{F,w = 0}^\natural
\longinto \pi_* DMT_\FS \bigl( M^K(\FS) \bigr)_{F,w=0}^{\natural,u} 
\]
is defined for a suitable choice of $\FS$. 
Theorem~\ref{0A} concerns the Tate motive 
$\one_{M^K}$, and direct factors
$N$ of its direct image
\[
\pi_* \one_{M^K} \in \pi_* DMT \bigl( M^K \bigr)_{F,w = 0}^\natural \; .
\] 
\end{Proofof}

\begin{Thm} \label{8E} 
The intermediate extension is compatible with the Betti and $\ell$-adic realizations;
more precisely, the following holds. \\
(a)~For any integer $m$, the restriction of the composition
\[
H^m \circ R_{(M^L)^*}: \DBcM \bigl( (M^L)^* \bigr)_F 
\longto \Perv_c \bigl( (M^L)^*_{E(G, \FH)},F \bigr) 
\]
to $\pi_* DMT_{\FS}(M^K(\FS))_{F,w =0}^\natural$
factors over $\pi_* DMT_{\FS}(M^K(\FS))_{F,w=0}^{\natural,u}$. \\[0.1cm]
(b)~For any integer $m$, the diagram
\[
\vcenter{\xymatrix@R-10pt{
    \pi_* DMT_{\FS_U} \bigl( \pi^{-1}(U) \bigr)_{F,w = 0}^\natural 
              \ar[d]_{H^m \circ R_U} \ar@{^{ (}->}[r]^-{\ujast} &
    \pi_* DMT_{\FS} \bigl( M^K(\FS) \bigr)_{F,w=0}^{\natural,u} 
              \ar[d]^{H^m \circ R_{(M^L)^*}} \\
    \Perv_c(U_{E(G, \FH)},F)
                    \ar@{^{ (}->}[r]^-{j_{!*}} &
    \Perv_c \bigl( (M^L)^*_{E(G, \FH)},F \bigr)
\\}}
\]
commutes. 
\end{Thm}

\begin{Proof}
Assumption~\ref{As2}~(a) and (b) hold thanks to 
Lemmata~\ref{8A} and \ref{8B}. So does
Assumption~\ref{As2}~(c) (by Lemma~\ref{8C}).

Now apply Theorem~\ref{7C}.
\end{Proof}

\medskip

\begin{Proofof}{Theorem~\ref{0B}}
This is a special case of Theorem~\ref{8E}.
\end{Proofof}

Now let $L'$ be a second
open compact subgroup of $G(\BA_f)$, $h \in G(\BA_f)$ such that
$L \subset h \cdot L' \cdot h^{-1}$, and 
\[
[ \ \cdot h] : (M^L)^* \longto (M^{L'})^*
\]
as before. Assume that the open subscheme $U$ of $(M^L)^*$ equals 
$[ \ \cdot h]^{-1}(U')$, for some open subscheme $U'$ of $(M^{L'})^*$.
Denote by $Z$ the complement of $U$ in $(M^L)^*$.

\begin{Thm} \label{8F}
The intermediate extension is compatible with the direct image $[ \ \cdot h]_*$;
more precisely, the following holds. \\[0.1cm]
(a)~The functor
\[
[ \ \cdot h]_* : \pi_* DMT_{\FS} \bigl( \pi^{-1}(Z) \bigr)_{F,w = 0}^\natural \longto 
\bigl( [ \ \cdot h] \circ \pi \bigr)_* DMT_{\FS} \bigl( \pi^{-1}(Z) \bigr)_{F,w = 0}^\natural
\]
is radicial. \\[0.1cm] 
(b)~The diagram
\[
\vcenter{\xymatrix@R-10pt{
    \pi_* DMT_{\FS_U} \bigl( \pi^{-1}(U) \bigr)_{F,w = 0}^\natural 
                    \ar[d]_{[ \ \cdot h]_*} \ar@{^{ (}->}[r]^-{\ujast} &
    \pi_* DMT_{\FS} \bigl( M^K(\FS) \bigr)_{F,w=0}^{\natural,u} \ar[d]_{[ \ \cdot h]_*^u} \\
    \bigl( [ \ \cdot h] \circ \pi \bigr)_* 
                                DMT_{\FS_U} \bigl( \pi^{-1}(U) \bigr)_{F,w = 0}^\natural 
                    \ar@{^{ (}->}[r]^-{\ujast} &
    \bigl( [ \ \cdot h] \circ \pi \bigr)_* 
                                DMT_{\FS} \bigl( M^K(\FS) \bigr)_{F,w=0}^{\natural,u}
\\}}
\]
commutes.
\end{Thm}

\begin{Proof}
Assumption~\ref{As}~(a), (b) and (d) hold thanks to 
Lemmata~\ref{8A}, \ref{8B} and \ref{8C}. So does
Assumption~\ref{As}~(c) (by Lemma~\ref{8Ba}).

Now apply Theorem~\ref{6G}.
\end{Proof}

\medskip

\begin{Proofof}{Theorem~\ref{0C}}
This is a special case of Theorem~\ref{8F}.
\end{Proofof}

Let us conclude with a discussion of Hecke operators. 
We keep the situation of Theorems~\ref{8D}--\ref{8F}.
In particular, we continue to assume that $(P, \FX)$ satisfies 
$(+)$, and that $K$ is neat. In addition, we impose two more conditions.
First, we assume that
the ``weight $-2$''par of $P$ is trivial, \emph{i.e.}, that
the unipotent radical $W$ is pure of weight $-1$. 
Second, we fix a Levi section $i: G \into P$ 
of the projection from $P$ to $G$, and
suppose that $K$ contains the image of $L$ under $i$.
In other words, $K$ is the semi-direct product of $K \cap W(\BA_f)$ and $i(L)$. 
According to
\cite[3.12--3.22~(a), Prop.~11.10]{P1}, the morphism $\pi: M^K \to M^L$
is then equipped with the structure of an Abelian scheme.
We are thus in the situation
of the Introduction. As there, we fix a direct factor $N$ of 
\[
\pi_* \one_{M^K} \in \CHFMLMs \; .
\]
Let us recall the action of the \emph{Hecke algebra} $R(L,G(\BA_f))$
associated to the Shimura variety $M^L$. Elements of $R(L,G(\BA_f))$
are formal linear combinations of double cosets $LxL$, for $x \in G(\BA_f)$.
Fix such an $x$, and define $L_x := L \cap x^{-1}Lx$. According to
what was recalled further above, there are
two finite, \'etale morphisms 
\[
g_1:= [\ \cdot 1] \; , \; g_2:= [\ \cdot x^{-1}]: M^{L_x} \longto M^L \; .
\]
In order to construct the action of $LxL$, we need three ingredients:
(1)~the adjunction morphism $adj_1: N \to g_{1,*} g_1^* N$,
(2)~the adjunction morphism $adj_2: g_{2,*} g_2^* N = g_{2,!} g_2^! N \to N$,
(3)~a canonical morphism $\varphi_x: g_1^* N \to g_2^* N$
of smooth Chow motives over $L_x$, whose construction
we recall now. By smooth (or proper) base change \cite[Thm.~2.4.50~(4)]{CD},
\[
g_r^* \pi_* \one_{M^K} = \pi_{r,*} \one_{M^{K_r}} \; , \; r= 1,2 \; ,
\]
where $\pi_r : M^{K_r}:= M^{K_r}(P,\FX) \to M^{L_x}$, 
$r=1,2$ are the projections associated
to $K_1 := K \cap \pi^{-1}(L_x)$ and $K_2 := y^{-1}Ky \cap \pi^{-1}(L_x)$,
for $y := i(x)$ (recall that $i$ is the Levi section fixed above).
Then $\pi_1$ and $\pi_2$ are simultaneously dominated by 
$\pi_{1,2} : M^{K_{1,2}}:= M^{K_{1,2}}(P,\FX) \to M^{L_x}$, for 
$K_{1,2} := K \cap y^{-1}Ky \cap \pi^{-1}(L_x)$; more precisely,
$\pi_{1,2}$ factors over both $M^{K_r}$ \emph{via} finite, \'etale morphisms
$M^{K_{1,2}} \to M^{K_r}$, $r = 1,2$. Consider
the composition of the adjunctions associated to these morphisms
\[
\psi_x: g_1^* \pi_* \one_{M^K} = \pi_{1,*} \one_{M^{K_1}} 
\longto \pi_{1,2,*} \one_{M^{K_{1,2}}} 
\longto \pi_{2,*} \one_{M^{K_2}} = g_2^* \pi_* \one_{M^K} \; .
\]
The desired morphism $\varphi_x$ is the composition of the
inclusion of $g_1^* N$ into $g_1^* \pi_* \one_{M^K}$, of $\psi_x$, and
of the projection of $g_2^* \pi_* \one_{M^K}$ onto $g_2^* N$. \\

Denote by $m^\circ$ the structure morphism $M^L \to \Spec E$,
for $E:= E (G , \FH)$.
Given (1), (2), and (3), we may now define
the endomorphism $LxL$ of $m^\circ_* N$. Namely, apply $m^\circ_*$ to
(1) and (2),
\[
m^\circ_* (adj_1): m^\circ_* N \longto (m^\circ \circ g_1)_* g_1^* N \; , \;
m^\circ_* (adj_2): (m^\circ \circ g_2)_* g_2^* N \longto m^\circ_* N \; ,
\]
and $(m^\circ \circ g_1)_* = (m^\circ \circ g_2)_*$ to (3),
\[
(m^\circ \circ g_1)_* (\varphi_x): 
(m^\circ \circ g_1)_* g_1^* N \longto (m^\circ \circ g_2)_*g_2^* N \; .
\]
Set
\[
LxL := 
m^\circ_* (adj_2) \circ (m^\circ \circ g_1)_* (\varphi_x) \circ m^\circ_* (adj_1): 
m^\circ_* N \longto m^\circ_* N \; .
\]
In the same way, one gets an endomorphism of $m^\circ_! N$, which
is denoted by the same symbol $LxL$. 
Denote the structure morphism $(M^L)^* \to \Spec E$ by $m$;
thus, 
\[
m^\circ_* N = m_* j_* N \quad \text{and} \quad 
m^\circ_! N = m_* j_! N \; .
\]
Denote the extension of $g_r$ to a morphism
$(M^{L_x})^* \to (M^L)^*$ \cite[Main Theorem~12.3~(b)]{P1}
by the same letter $g_r$, $r = 1,2$.

\begin{Cor} \label{8G}
(a)~The action of $LxL$ can be extended to $m_* \ujast N$, 
in a way compatible with the actions on $m_* j_! N$ and on $m_* j_* N$: 
there is a commutative diagram
\[
\vcenter{\xymatrix@R-10pt{
        m_* j_! N \ar[d]_{LxL} \ar[r] &
        m_* \ujast N \ar[d]_{LxL} \ar[r] &
        m_* j_* N \ar[d]^{LxL} \\
        m_* j_! N \ar[r] &
        m_* \ujast N \ar[r] &
        m_* j_* N 
\\}}
\]
(b)~Assume that one of the following additional conditions (i), (ii) are satisfied: 
(i)~the radicals
\[
\rad_{\CHFMLastM}(\ujast N, g_{1,*} \ujast g_2^* N) 
\subset \Hom_{\CHFMLastM}(\ujast N, g_{1,*} \ujast g_2^* N)
\]
and 
\[
\rad_{\CHFMLastM}(g_{2,*} \ujast g_2^* N, \ujast N) 
\subset \Hom_{\CHFMLastM}(g_{2,*} \ujast g_2^* N, \ujast N)
\]
are trivial, (ii)~the radicals
\[
\rad_{\CHFMLastM}(\ujast N, g_{1,*} \ujast g_1^* N) 
\subset \Hom_{\CHFMLastM}(\ujast N, g_{1,*} \ujast g_1^* N)
\]
and
\[
\rad_{\CHFMLastM}(g_{2,*} \ujast g_1^* N, \ujast N) 
\subset \Hom_{\CHFMLastM}(g_{2,*} \ujast g_1^* N, \ujast N)
\]
are trivial. Then the action of $LxL$ on $m_* \ujast N$
is canonical. \\[0.1cm]
(c)~Assume that the realization
$R_{M^L}(N)$ is concentrated in a single perverse degree. 
Then the endomorphism
\[
R_{\Spec E} (LxL) : R_{\Spec E}(m_* \ujast N) \longto R_{\Spec E}(m_* \ujast N)
\]
coincides with the Hecke operator defined on the complex computing
intersection cohomology \emph{via} the isomorphism
\[
R_{\Spec E} (m_* \ujast N) \cong m_* \ujast R_{M^L} (N)
\]
from Corollary~\ref{0Bcor}.
\end{Cor} 

\begin{Rem}
(a)~One may expect the additional hypothesis of \ref{8G}~(b) to be
met as soon as $N$ is indecomposable (cmp.~Summary~\ref{2M}~(c),
Conjecture~\ref{3C}). \\[0.1cm]
(b)~Recall that the motive
$\pi_* \one_{M^K}$ admits a Chow--K\"unneth decomposition \cite[Thm.~3.1]{DM}.
Any direct factor $N$ contained in a Chow--K\"unneth component 
(in particular, any indecomposable $N$) then satisfies 
the additional hypothesis of \ref{8G}~(c). \\
(c)~Assume that the Shimura data $(G,\FH)$ are  \emph{of PEL-type}.
According to the main result from \cite{Anc}, 
the isomorphism classes of indecomposable 
direct factors of $\pi_* \one_{M^K}$ are in bijective correspondence with certain of the
isomorphisms classes of irreducible algebraic representations of $G$. 
More precisely, the Chow motive $\pi_* \one_{M^K}$
corresponds to a direct sum $V$ of irreducible representations of $G$, and any
constituent $V'$ of $V$ induces a direct factor $N$ of $\pi_* \one_{M^K}$
\cite[Thm.~8.6]{Anc}. 
\end{Rem}

\begin{Proofof}{Corollary~\ref{8G}}
Fix a $K_{1,2}$-admissible complete cone decomposition 
$\FS$ satisfying Lemma~\ref{8A}~(b).
Applying $\ujast$ to (1) and (2), we get morphisms
\[
\ujast (adj_1): \ujast N \longto \ujast g_{1,*} g_1^* N \quad \text{and} \quad
\ujast (adj_2): \ujast g_{2,*} g_2^* N \longto \ujast N 
\]
in the quotient category 
$(g_1 \circ \pi_{1,2} \coprod g_2 \circ \pi_{1,2})_* 
DMT_{\FS} ( M^{K_{1,2}}(\FS) )_{F,w=0}^{\natural,u}$.
By Theo\-rem~\ref{8F}~(b),
\[
\ujast g_{r,*} g_r^* N = g_{r,*}^u \ujast g_r^* N \; , \; r = 1,2 \; .
\]
Applying $\ujast$ to (3), we get
\[
\ujast (\varphi_x): \ujast g_1^* N \longto \ujast g_2^* N
\]
in $\pi_{1,2,*} DMT_{\FS} ( M^K(\FS) )_{F,w=0}^{\natural,u}$, yielding
\[
g_{1,*}^u \ujast (\varphi_x): 
g_{1,*}^u \ujast g_1^* N \longto g_{1,*}^u \ujast g_2^* N
\]
in $(g_1 \circ \pi_{1,2} \coprod g_2 \circ \pi_{1,2})_* 
DMT_{\FS} ( M^{K_{1,2}}(\FS) )_{F,w=0}^{\natural,u}$. Composing the latter
with $\ujast (adj_1)$, we get
\[
g_{1,*}^u \ujast (\varphi_x) \circ \ujast (adj_1): 
\ujast N \longto g_{1,*}^u \ujast g_2^* N \; .
\]
Lift both this latter morphism and $\ujast (adj_2)$ to actual morphisms
\[
\ujast N \longto g_{1,*} \ujast g_2^* N \quad \text{and} \quad
g_{2,*} \ujast g_2^* N \longto \ujast N
\]
of Chow motives over $(M^L)^*$. 
According to Theorems~\ref{8F}~(b) and \ref{1H}~(b), such lifts are unique
up to morphisms belonging to the radical.
Applying $m_*$, we get
\[
m_* \ujast N \longto (m \circ g_1)_* \ujast g_2^* N \quad \text{and} \quad
(m \circ g_2)_* \ujast g_2^* N \longto m_* \ujast N \; .
\]
But since $m \circ g_1 = m \circ g_2$, the two morphisms can be composed;
this is the desired endomorphism
\[
LxL : m_* \ujast N \longto m_* \ujast N \; .
\]
This proves parts~(a) and (b). Part~(c) follows from 
Theorem~\ref{8E}~(a) and Corollary~\ref{0Bcor}.
\end{Proofof}

\medskip

\begin{Proofof}{Theorem~\ref{0D}}
Theorem~\ref{0D} is identical to parts (a) and (c) of Corollary~\ref{8G}.
\end{Proofof}

\begin{Rem} \label{8H}
For specific Shimura varieties, certain of the above results,
or variants thereof, are already known. \\[0.1cm]
(a)~In \cite[Sect.~1]{Sch}, the construction of the intersection motive $m_* \ujast N$
of modular curves with coefficients in symmetric powers $N$ of the relative ``$h^1$''
of the universal elliptic curve is given. It carries an action of Hecke operators,
and satisfies the properties from Theorem~\ref{0D} \cite[Prop.~4.1.3]{Sch}.
Actually, all the ingredients to define the intermediate extension $\ujast N$ itself can be
found in \loccit --- except for the notion of Chow motive over $(M^L)^*$.
Scholl then proceeds and uses the Hecke action to construct 
(Grothendieck) motives associated to modular forms. \\[0.1cm]
(b)~The main result of \cite{GHM} applies to Hilbert--Blumenthal varieties $M^L_\BC$
over $\BC$. Indeed, according to \cite[Thm~I]{GHM}, for any Kuga family $\pi: M^K \to M^L$,
any direct factor $N$ of $\pi_* \one_{M^K_\BC}$ admits an extension to a Chow motive
on $(M^L)^*_\BC$ satisfying 
Theorem~\ref{0A}~(a)~(1); according to Remark~\ref{6I},
it is thus isomorphic to the base change of the intermediate extension of $N$   
to $\BC$. Compatibility with the Betti realization (Theorem~\ref{0B})
is also established in \cite[Thm~I]{GHM}. \\[0.1cm]
(c)~The intersection motive of a surface $S$ (or more generally, of a variety admitting a
\emph{semismall resolution}) over $\BC$, with constant coefficients,
is constructed in \cite[Thm.~4.0.4]{CM}. It satisfies Theorem~\ref{0B} \cite[Rem.~4.0.6]{CM}.
For surfaces, all the ingredients to define the intermediate extension 
itself can actually be found in \loccit . This construction even works over arbitrary
fields. For details, we refer to \cite[Thm.~3.11~(b)]{W5}. \\[0.1cm]
(d)~By \cite[Thm.~1]{MS}, any direct factor
$N$ of $\pi_* \one_{M^K_\BC}$ admits an extension to a Chow motive
on $(M^L)^*_\BC$ satisfying 
Theorem~\ref{0B}, as soon as $M^K$ is of dimension $3$. \\[0.1cm]
(e)~Let $M^L$ be a Hilbert--Blumenthal variety, and $N$ a direct factor
of $\pi_* \one_{M^K}$ not containing any shift of a Tate twist.
Then according to \cite[Cor.~2.7]{W4}, the \emph{interior motive}
of $M^L$ with coefficients in $N$ (which in this geometric setting
coincides with the intersection motive) exists. 
Theorem~\ref{0D} is proved in \cite[Cor.~2.8]{W4},
and Theorem~\ref{0B} follows from \cite[Thm.~4.7]{W3}. \\[0.1cm]
(f)~Let $M^L$ be a Picard surface, and $N$ a direct factor
of $\pi_* \one_{M^K}$, which under \cite[Thm.~8.6]{Anc} 
is induced by an irreducible
representation of $G$, whose highest weight is \emph{regular}.
Then according to \cite[Cor.~3.9]{W6}, the interior motive
of $M^L$ with coefficients in $N$ (which in this geometric setting, too, 
coincides with the intersection motive) exists. 
Theo\-rem~\ref{0D} is proved in \cite[Cor.~3.10]{W6},
and Theorem~\ref{0B} follows from \cite[Thm.~4.7]{W3}. 
We refer to \cite{Clt} for the generalization of these results to Picard varieties
of arbitrary dimension. \\[0.1cm]
(g)~To the best of the author's knowledge, the first non-trivial example
of the motivic intermediate extension over a Baily--Borel compactification
of dimension at least $3$, 
and defined over the reflex field of the Shimura data, appears in \cite{V}. 
There, the author shows the existence of a
motive $EM_X$ \cite[Def.~3.2.12]{V}, which is a pre-image of 
the \emph{weight truncated complex} $EC_X$ \cite[Def.~3.1.15]{NV1}
under the Betti realization \cite[Prop.~4.2.4]{V}. He also stu\-dies the relation
of $EC_X$ and the intersection complex $IC_X$ \cite[Thm.~4.1.2]{V};
in particular, both motives coincide if $X$ is the Baily--Borel compactification 
of a Hilbert--Blumenthal variety. 
This concerns the case of constant coefficients 
that was left open in \cite{W4}.  \\[0.1cm]
(h)~In \cite{NV2}, the authors consider the case of Siegel threefolds $M^L$,
and the universal family $\pi: M^K \to M^L$ of Abelian surfaces with 
an $n$-structure, for some $n \ge 3$. 
\cite[Thm.~I, Thm.~II]{NV2} establish Theorem~\ref{0A}~(a) and (d),
for $N = \pi_* \one_{M^K}$.
Contrary to our present approach, the construction of \cite{NV2} is very
explicit. For a particular choice of $\FS$, it is shown that the action of
the $n$-torsion extends to $M^K(\FS)$. Writing $\pi: M^K(\FS) \to (M^L)^*$
as before, it is then established that the direct factor of $\pi_* \one_{M^K(\FS)}$
on which the $n$-torsion acts trivially, has no components supported in the
boundary of $(M^L)^*$. Therefore, it is an intermediate extension of 
$\pi_* \one_{M^K}$. The claim from Theorem~\ref{0A}~(a)~(2) is verified by using
intersection theory of cycles on $M^K(\FS)$.
\end{Rem}

\begin{Rem} \label{8I}
Given the recent results from \cite{La1,La2}, it appears likely that analogues 
of the results of the present section, and hence of Theorems~\ref{0A}--\ref{0D},
hold for \emph{integral models} of Baily--Borel compactifications
of arbitrary (pure) PEL-type Shimura varieties
over primes of good reduction, and for Chow motives
occurring as direct factors in the relative motive of integral models of
Kuga families over such Shimura varieties. More precisely, the analogue
of Lemma~\ref{8A} holds in this context (\cite[Thm.~2.15~(1), Sect.~2C~(1)]{La2},
together with \cite[Thm.~6.4.1.1~2. and 3.]{La1}), and so does Lemma~\ref{8B}~(a)
\cite[Thm.~7.2.4.1~4.]{La1}. As in the proof of Lemma~\ref{8B}, the strata of 
the model of the Baily--Borel compactification are quotients, denoted $[M_{\CH}^{Z_{\CH}}]$
in \cite[Thm.~7.2.4.1~4.]{La1}, of smooth objects $M_{\CH}^{Z_{\CH}}$; given the
limitations of our present understanding of $\cite{La1}$, it does not appear
obvious that the quotient map $M_{\CH}^{Z_{\CH}} \to [M_{\CH}^{Z_{\CH}}]$ is  
finite and \'etale (which would imply the analogue of Lemma~\ref{8B}~(b)).
Given what is called the ``Hecke action'' on the models of Baily--Borel compactifications 
in \cite[Prop.~7.2.5.1]{La1}, the analogue of Lemma~\ref{8Ba} would require
the same information on $M_{\CH}^{Z_{\CH}} \to [M_{\CH}^{Z_{\CH}}]$;
in addition, as in its proof, we would need $[ \ \cdot h]$ to lift to 
a morphism $M_{\CH}^{Z_{\CH}} \to M_{\CH'}^{Z_{\CH'}}$, which is finite and \'etale.
In order to obtain a statement analogous to Lemma~\ref{8C}, one would need to analyze
the restriction of the morphism of models of toroidal compactifications from
\cite[Thm.~2.15~(2)]{La2} to individual strata. 
Finally, the proof of the analogue of Corollary~\ref{8G}
for integral models of arbitrary $PEL$-type Shimura varieties
would involve \cite[Thm.~2.15~(4)]{La2}. 
\end{Rem}

Note that in the context of integral models of \emph{Siegel varieties}, the missing geometric 
information is available \cite{FC};
let us finish with a more detailed discussion of that situation. 
Fix an integer $g \ge 1$, consider the 
reductive group $G_{2g} := CSp_{2g,\BQ}$, and the pure Shimura data 
$(G_{2g},\FH_{2g})$ of \cite[Ex.~2.7]{P1}. Fix a neat open compact
subgroup $L$ of $G_{2g} (\BA_f)$, and suppose that it lies between two 
\emph{principal} open compact subgroups,
\[
L_k \le L \le L_n \; , \; n \tei k \; ,
\]
defined as the kernels of the reduction maps modulo $k$ and $n$ on $G_{2g} (\widehat{\BZ})$,
respectively. We also suppose that $n \ge 3$. It follows from the Main Theorem of 
Shimura--Taniyama on complex multiplication (see \cite[Thm.~4.19, Thm.~4.21]{D0})
that $M^{L_k}$ equals the moduli space of principally polarized Abe\-lian varieties
of dimension $g$ with principal level-$k$ structure, denoted $\A_{g,k}$.
Given the modular interpretation of the latter, it admits a canonical extension to
$\Spec \BZ[\frac{1}{k}]$. The resulting scheme is denoted by the same symbol
$\A_{g,k}$; it is identical to the object from \cite[Def.~I.4.4]{FC}. 
The scheme $\A_{g,k}$ is quasi-projective and smooth over $\Spec \BZ[\frac{1}{k}]$, and the 
morphism $[ \ \cdot 1]: M^{L_k} \onto M^{L_n}$ extends to a finite, \'etale morphism
$\A_{g,k} \onto \A_{g,n} \otimes \Spec \BZ[\frac{1}{k}]$ \cite[Rem.~IV.6.2~(c)]{FC}.
In fact \cite[Rem.~IV.6.2~(d)]{FC}, this morphism is identical to the map to the
quotient by the natural (right) action of $L_n / L_k$ on $\A_{g,k}$; this action
is free since $n \ge 3$. An integral model of $M^L$ over 
$\Spec \BZ[\frac{1}{k}]$ can thus be defined as the partial quotient by
the action of $L / L_k$:
\[
\A_{g,L} := \A_{g,k} / (L / L_k) \; .
\] 
It is quasi-projective and smooth over $\Spec \BZ[\frac{1}{k}]$.
This construction can be employed to show that \emph{mutatis mutandis}, 
the main results from \cite{FC} on $\A_{g,k}$ are equally valid for $\A_{g,L}$. 
In particular, the Baily-Borel compactification $(M^L)^*$
admits a proper model $\A_{g,L}^*$ over $\Spec \BZ[\frac{1}{k}]$, defined as
\[
\A_{g,L}^* := \A_{g,k}^* / (L / L_k)
\]
\cite[Thm.~V.2.5~(3)]{FC}. \\

Now fix a second integer $s \ge 0$.
The mixed Shimura data $(P,\FX)$ is supposed
to equal the $s$-fold fibre product of the data $(P_{2g},\FX_{2g})$ of \cite[Ex.~2.25]{P1};
the underlying pure Shimura data thus equals $(G_{2g},\FH_{2g})$, and the unipotent
radical of $P$ equals the $s$-th power of the standard representation $V_{2g}$
of $G_{2g}$. Let $K$ be an open compact sub-group of $P_{2g} (\BA_f)$, whose image
under the projection $\pi : (P,\FX) \to (G_{2g},\FH_{2g})$ equals $L$, and which 
contains the image of $L$ under the canonical Levi section of $\pi$. As in
the pure case discussed above, the main results from \cite{FC} on the $s$-fold
fibre product $\B_{g,k} \to \A_{g,k}$ of the universal Abelian scheme over $\A_{g,k}$ 
(denoted $Y$ in \loccit ) generalize to the model $\pi: \B_{g,K} \to \A_{g,L}$ of $M^K$ over
$\Spec \BZ[\frac{1}{k}]$. In particular, for suitable choices of complete 
cone decomposition $\FS$, the toroidal compactification $M^K(\FS)$
extends to a proper, smooth model $\B_{g,K} (\FS)$ over $\Spec \BZ[\frac{1}{k}]$,
and the morphism $\pi: M^K(\FS) \to (M^L)^*$ extends to a proper morphism
$\B_{g,K} (\FS) \to \A_{g,L}^*$ \cite[Thm.~VI.1.1]{FC}, 
which we shall still denote by $\pi$. \\

Fix a direct factor
$N$ of $\pi_* \one_{\B_{g,K}}$, viewed as an object of the category $\CHFALM$;
note that according to \cite[Prop.~5.1.1]{O'S2}, the isomorphism classes of such
$N$ are in bijective correspondence with the direct factors of the generic fibre
$\pi_* \one_{M^K}$ considered in the Introduction. 
Denote by $j: \A_{g,L} \into \A_{g,L}^*$ the open immersion, and by 
$i: \partial \A_{g,L}^* \into \A_{g,L}^*$ the complement of $\A_{g,L}$.

\begin{Thm} \label{8K}
(a)~There is a Chow motive $\ujast N \in \CHFALastM$
extending $N$, and satisfying
the following properties.
\begin{enumerate}
\item[(1)] $\ujast N$ admits no non-zero direct
factor belonging to $i_* \CHFpA$.
\item[(2)] Any element of the kernel of
\[
j^* : \End_{\CHFALastM}(\ujast N) \longto \End_{\CHFALM}(N)
\]
is nilpotent. 
\end{enumerate}
(b)~The Chow motive $\ujast N$ satisfies analogues of properties~(b)--(d) from Theorem~\ref{0A}. \\[0.1cm]
(c)~The construction of $\ujast N$ is compatible with Theorem~\ref{0A}
under pull-back to the generic fibre $\Spec \BQ$ of $\Spec  \BZ[\frac{1}{k}]$. 
In particular, it is compatible with the Betti and $\ell$-adic realizations.  
\end{Thm}

Now let $[\ \cdot h]: M^L \to M^{L'}$ be a finite morphism
associated to change of the ``level'' $L$, with 
\[
L_{k'} \le L' \le L_{n'} \; , \; n' \tei k' 
\]
and $n' \ge 3$. Then $[\ \cdot h]$
extends to a finite morphism 
\[
\A_{g,L}^* \otimes \Spec \BZ[\frac{1}{k \cdot k'}] \longto 
\A_{g,L'}^* \otimes \Spec \BZ[\frac{1}{k \cdot k'}] \; ,
\]
denoted by the same symbol $[\ \cdot h]$. 

\begin{Thm} \label{8L}
The intermediate extension is compatible with $[\ \cdot h]_*$. More precisely,
the Chow motive $[\ \cdot h]_* \ujast N \in \CHFAtLastM$ satisfies
the analogues of the properties from Theorem~\ref{8K}. 
\end{Thm}

Write $m$ for the structure morphism $\A_{g,L} \to \Spec \BZ[\frac{1}{k}]$.
Fix $x \in G_{2g}(\BA_f)$, and let $k_x$ be an integer satisfying
\[
L_{k_x} \le L_x := L \cap x^{-1} L x \; .
\] 

\begin{Thm} \label{8M}
The Hecke operator $LxL$ acts on $m_* \ujast N \otimes \Spec \BZ[\frac{1}{k_x}]$, in a way 
compatible with its action on $m_* j_! N \otimes \Spec \BZ[\frac{1}{k_x}]$ and on 
$m_* j_* N \otimes \Spec \BZ[\frac{1}{k_x}]$,
and with the action from Theorem~\ref{0D} on the generic fibre. 
\end{Thm}

\begin{Proofof}{Theorems~\ref{8K}, \ref{8L} and \ref{8M}}
The integral versions of
Lemmata~\ref{8A} and \ref{8B}~(a) hold (see Remark~\ref{8I}).
As for Lemma~\ref{8B}~(b), the quotient map $M_{\CH}^{Z_{\CH}} \to [M_{\CH}^{Z_{\CH}}]$, 
is the identity \cite[Thm.~V.2.5~(4)]{FC},
as was explained in Remark~\ref{8I}; in particular, it is finite and \'etale.
The same result (and the modular interpretation of $[ \ \cdot h]$) shows that
$[ \ \cdot h]$ lifts to 
a morphism $M_{\CH}^{Z_{\CH}} \to M_{\CH'}^{Z_{\CH'}}$, which is finite and \'etale.
We thus get the integral version of Lemma~\ref{8Ba}, 
and also the extension of Hecke operators to integral
models. 
As for Lemma~\ref{8C}, we refer to \cite[Thm.~V.2.5~(5) 
and (6)]{FC} (if $s = 0$) and to 
\cite[Thm.~VI.1.1 and its proof]{FC} (if $s \ge 1$). 
\end{Proofof}  


\bigskip

%
%

\end{document}